\documentclass{elsarticle}

\usepackage{amsmath,amssymb,amsthm,mathrsfs,epsf,a4,color,graphicx,eucal}

\newtheorem{theorem}{Theorem}[section]
\newtheorem{lemma}[theorem]{Lemma}
\newtheorem{definition}[theorem]{Definition}

\newtheorem{remark}[theorem]{Remark}
\newtheorem{notation}[theorem]{Notation}
\newtheorem{problem}[theorem]{Problem}

\numberwithin{equation}{section}

\newcommand{\COL}[1]{{\color{black}#1}}
\newcommand{\COLL}[1]{{\color{black}#1}}
\newcommand{\COLLL}[1]{{\color{black}#1}}

\definecolor{ddmagenta}{rgb}{0.7,0,0.9}

\newcommand{\ITEM}[2]{\parbox[t]{.05\textwidth}{#1}\hfill\parbox[t]{.95\textwidth}{#2}\vspace*{.8mm}}

\newcommand{\DDD}[3]{\begin{array}[t]{c}#1\vspace*{-1em}\\_{#2}\vspace*{-.5em}\\_{#3}\end{array}}
\newcommand{\ddd}[3]{\DDD{\begin{array}[t]{c}\underbrace{#1}\vspace*{.6em}\end{array}}{\text{\footnotesize #2}}{\text{\footnotesize #3}}}

\newcommand\eps{\varepsilon}
\newcommand\R{\mathbb R}
\newcommand\N{\mathbb N}

\newcommand{\eq}[1]{(\ref{#1})}
\newcommand\DT[1]{\mathchoice
                 {{\buildrel{\hspace*{.1em}\text{\LARGE.}}\over{#1}}}
                 {{\buildrel{\hspace*{.1em}\text{\Large.}}\over{#1}}}
                 {{\buildrel{\hspace*{.1em}\text{\large.}}\over{#1}}}
                 {{\buildrel{\hspace*{.1em}\text{\large.}}\over{#1}}}}
\newcommand\DDT[1]{\mathchoice
   {{\buildrel{\hspace*{.1em}\text{\LARGE.\hspace*{-.1em}.}}\over{#1}}}
   {{\buildrel{\hspace*{.1em}\text{\Large.\hspace*{-.1em}.}}\over{#1}}}
   {{\buildrel{\hspace*{.1em}\text{\large.\hspace*{-.1em}.}}\over{#1}}}
   {{\buildrel{\hspace*{.1em}\text{\large.\hspace*{-.1em}.}}\over{#1}}}}
\newcommand{\pl}{\partial}
\newcommand\JUMP[2]{\mathchoice
                   {\big[\hspace*{-.3em}\big[#1\big]\hspace*{-.3em}\big]_{#2}}
                   {[\hspace*{-.15em}[#1]\hspace*{-.15em}]_{#2}}
                   {[\![#1]\!]_{#2}}
                   {[\![#1]\!]_{#2}}}
\newcommand\bbC{\mathbb C}
\newcommand\bbD{\mathbb D}
\newcommand\bbE{\mathbb E}
\newcommand\bbK{\mathbb K}
\newcommand\bbB{\mathbb B}
\renewcommand\d{\mathrm d}
\newcommand\w{\vartheta}
\newcommand\ent{h}
\renewcommand\pl{\partial}
\newcommand\GC{\Gamma_{\mbox{\tiny\rm C}}}
\newcommand\SC{\Sigma_{\mbox{\tiny\rm C}}}
\newcommand{\dt}{\mathrm{D}_t}

\newcommand{\FRM}{F}
\newcommand{\GRM}{G}
\newcommand{\fRM}{f}
\newcommand{\gRM}{g}

\newcommand{\IRM}{\mathrm{I}}
\newcommand{\calD}{\mathcal{R}}

\newcommand{\aein}{\text{a.e. in}}
\newcommand{\foraa}{\text{for a.a.}}

\newcommand{\testu}{v}
\newcommand{\testw}{w}

\newcommand{\minus}{-}
\newcommand{\Gdir}{\Gamma_{\mbox{\tiny\rm D}}}
\newcommand{\Gnew}{\Gamma_{\mbox{\tiny\rm N}}}
\newcommand{\Sdir}{\Sigma_{\mbox{\tiny\rm D}}}
\newcommand{\Snew}{\Sigma_{\mbox{\tiny\rm N}}}

\newcommand{\dela}{\kappa}
\newcommand{\delam}{\kappa}
\newcommand{\het}{\eta}
\newcommand{\ind}{I}
\newcommand{\dd}{\mathrm{d}}

\newcommand{\piecewiseConstant}[2]{\overline{#1}_{\kern-1pt#2}}
\newcommand{\pwc}{\piecewiseConstant}
\newcommand{\underpiecewiseConstant}[2]{\underline{#1}_{\kern-1pt#2}}
\newcommand{\upwc}{\underpiecewiseConstant}

\newcommand{\piecewiseLinear}[2]{#1_{\kern-1pt#2}}
\newcommand{\pwl}{\piecewiseLinear}

\newcommand{\weaksto}{\stackrel{*}{\rightharpoonup}}
\newcommand{\weakto}{{\rightharpoonup}\,}
\newcommand{\GE}{\succeq}
\newcommand{\GEstar}{\succeq\hspace{-.8em}^{^*}\hspace{.3em}}
\newcommand{\BV}{\mathrm{BV}}
\newcommand{\pairing}[4]{ \sideset{_{#1 }}{_{ #2}}  {\mathop{\langle #3 , #4  \rangle}}}

\newcommand{\ut}{u_{\eps\tau}}
\newcommand{\wt}{\w_{\eps\tau}}
\newcommand{\zt}{z_{\eps\tau}}
\newcommand{\ue}{u_{\eps}}
\newcommand{\we}{\w_{\eps}}
\newcommand{\ze}{z_{\eps}}
\newcommand{\uej}{u_{\eps}}
\newcommand{\wej}{\w_{\eps}}
\newcommand{\zej}{z_{\eps}}
\newcommand{\ude}{u}
\newcommand{\wde}{\w}
\newcommand{\zde}{z}
\newcommand{\udelj}{u_{\delam}}

\newcommand{\zdelj}{z_{\delam}}
\newcommand{\uk}{u_{\eps\tau}^k}
\newcommand{\wk}{\w_{\eps\tau}^k}
\newcommand{\zk}{z_{\eps\tau}^k}

\newcommand{\yosd}{(I_{\cone}^\varepsilon)'}
\newcommand{\yosdej}{(I_{\cone}^{\varepsilon_j})'}
\newcommand{\maxad}{\mathcal{A}}

\newcommand{\dom}{\mathrm{dom}}
\newcommand{\cone}{K}

\newcommand{\indj}{J}

\newcommand{\Omegaone}{\Omega_+}
\newcommand{\Omegatwo}{\Omega_-}
\newcommand{\indabs}{\partial \abs}
\newcommand{\abs}{\mathcal{I}_{\cone}}
\newcommand{\norm}{\nu^{\pm}}
\newcommand{\varmea}{\xi_{\DT{z}}^{\mathrm{surf}}}
\newcommand{\varmeaps}{\xi_{\DT{z}_\eps}^{\mathrm{surf}}}
\newcommand{\jum}{J}
\newcommand{\yosapp}{I_{\cone}^{\eps}}
\newcommand{\yosappej}{I_{\cone}^{\eps}}

\journal{Nonlinear Analysis: Theory, Methods \& Applications}

\begin{document}

\begin{frontmatter}

\title{Thermodynamics and analysis of rate-independent\\[.5mm]
       adhesive contact at small strains}

\author[brescia]{Riccarda Rossi}
\ead{riccarda.rossi@ing.unibs.it}
\author[cuni,cas]{Tom\'a\v s Roub\'\i\v cek}
\ead{tomas.roubicek@mff.cuni.cz}

\address[brescia]{Dipartimento di Matematica, Universit\`a di
Brescia,  Via Valotti 9, I--25133 Brescia, Italy}

\address[cuni]{Mathematical Institute, Charles University, Sokolovsk\'a 83,
CZ--186~75~Praha~8,  Czech Republic}

\address[cas]{Institute of Thermomechanics of the ASCR, Dolej\v skova~5,
CZ--182 00 Praha 8, Czech Republic}

\begin{abstract}
We address a model for adhesive unilateral frictionless Signorini-type
contact between bodies of heat-conductive viscoelastic material, in
the linear Kelvin-Voigt rheology, undergoing thermal expansion. The
flow-rule for debonding the adhesion is considered rate-independent
and unidirectional, and a thermodynamically consistent model is
derived and analysed as far as the existence of a weak solution is
concerned.
\end{abstract}

\begin{keyword}
 Adhesive contact \sep nonlinear heat
equation \sep rate-independence \sep energetic solution \sep
existence \MSC 35K85 \sep 49S05 \sep 74A15 \sep 74M15 \sep 80A17
\end{keyword}

\end{frontmatter}

\section{Introduction}\label{sec-intro}
 We are interested in the modelling of elastic bodies glued
together by an adhesive, which can undergo an inelastic process of
so-called delamination (sometimes also called debonding).
``Microscopically'' speaking, some macromolecules in the adhesive
may break upon loading and we assume that they can never be glued
back, i.e., no ``healing'' is possible. This makes the process {\it
unidirectional}; sometimes it is  also referred to as  {\it
irreversible}, although this adjective has  an alternative
thermodynamical meaning as dissipative in general. On the glued
surface, we consider the {\it delamination process as
rate-independent} and, in the bulk, we also consider rate-dependent
{\it inertial, viscous-like}, and {\it thermal-expansion effects}.
Moreover, we confine ourselves to {\it small strains} and, just for
the sake of notational simplicity,  we restrict the analysis to the
case of  two bodies $\Omegaone$ and $\Omegatwo$ glued together along
the {\it contact surface} $\GC$. The material in the bulk is
considered as heat conductive, and thus the system is completed by
the nonlinear heat equation in a  thermodynamically consistent way.
The contact surface is considered infinitesimally thin, so that the
thermal capacity of the adhesive is naturally neglected. The
coupling of the mechanical and thermal effects  thus results from
thermal expansion, dissipative/adiabatic heat
production/consumption, and here also from the
possible dependence of the heat-transfer through the contact surface
$\GC$ on the
delamination itself, and on the possible slot between the bodies if
the contact is debonded.

We consider an elastic response of the adhesive, and then one speaks
about  {\it adhesive contact} (in contrast to  {\it brittle
contact}, see Remark~\ref{rem-brittle}). Within the realm of the
literature on (frictionless adhesive) contact, in the isothermal
case  we  refer  e.g. to~\cite{{KoMiRo}} in the framework of
rate-independent problems.  For \emph{rate-dependent} models, we
mention~\cite{BBR1,BBR2,Fre82,Fre87,Point,Raous} (cf. the
monograph~\cite{SoHaSh06AACP} for further references). The
anisothermal rate-dependent case has been recently addressed
in~\cite{BBR3,BBR4}. The present paper extends the analysis
in~\cite{KoMiRo} of rate-independent adhesive contact, to encompass
inertial, viscous, and thermal effects.

The elastic response in the adhesive will be considered linear,
determined by the scalar elastic modulus $\kappa>0$;
cf.~Remark~\ref{rem-kappa} for a generalization. At a current time,
the ``volume fraction'' of debonded molecular links will be
``macroscopically'' described by the {\it scalar delamination
parameter} $z:\GC\to[0,1]$. The state $z(x)=1$ means that the
adhesive is still $100\%$ undestroyed and thus fully effective,
while the intermediate state $0<z(x)<1$ means that there are some
molecular links which have been broken but the remaining ones are
effective, and eventually $z(x)=0$ means that the surface is already
completely debonded at $x\in\GC$. As already pointed out
in~\cite{KoMiRo},  one needs a specific energy to break the
macromolecular structure of the adhesive, independently of the rate
of this process. Thus,  delamination is a \emph{rate-independent}
and activated phenomenon, governed by the maximum dissipation
principle, and we shall accordingly consider a rate-independent flow
rule for $z$. Activating the delamination process in the adhesive
contact at a given point $x\in\GC$ again needs the
(phenomenologically prescribed) energy $a(x)$.

In the thermodynamical context, the energy $a(x)$ needed for
delamination is dissipated by the system in two ways: one part $a_1$
is spent to the chaotic vibration of the atomic lattice of both
sides of the delaminating surface $\GC$, which leads
``macroscopically'' to heat production (cf.~also \cite[Remark
4.2]{tr-LS-CZ}), while another part $a_0$ is spent to create a new
delaminated surface (or, ``microscopically'' speaking, to break the
macromolecules of the adhesive). Thus $a(x)=a_0(x)+a_1(x)$.

The mathematical difficulties, arising both from the proper
thermodynamical coupling and from hosting a rate-independent process
on $\GC$, have been already revealed for other inelastic processes
in the bulk in \cite{tr1}. The essential ingredient is the
satisfaction of the energy balance and, for this, the mentioned
concept of energetic solutions to rate-independent systems recently
developed in
\cite{Miel05ERIS,MieThe99MMRI,MieThe04RIHM,MiThLe02VFRI}, and
adapted to systems with inertia and viscosity in \cite{tr2}, appears
truly essential.


In Section~\ref{sec2}, we set up our model and, in
Sect.~\ref{sec-th}, discuss its thermodynamics and various
modifications. After making a suitable transformation of the problem
using an enthalpy variable instead of  the temperature, and
introducing a suitable weak formulation in
Sect.~\ref{sec-ent-trans}, the main existence results are presented
in Sect.~\ref{sec-main-res}, and proved throughout
Sections~\ref{s:4}--\ref{s:5}. For this, in Sect.~\ref{s:4} we set
up procedures of regularization of the Signorini-type unilateral
contact. As we shall observe in Sect.~\ref{s:4}, such a regularized
problem has its own interest. We further approximate it by
convexifying  some nonlinear terms, and setting up a
time-discretization procedure in Sect.~\ref{s:new}. Hence, we prove
fine a-priori estimates. Ultimately, a careful passage to the limit
is executed in two consecutive steps in Sections~\ref{ss:4.4}
and~\ref{s:5}.

\section{The model}\label{sec2}
Hereafter, we  suppose that the elastic body occupies a reference
domain
$$
\text{$\, \Omega\subset\R^d\,$, $d=2$ or $3$,  bounded and with a
Lipschitz boundary $\partial\Omega$.}
$$
We assume that
\[
\Omega = \Omegaone \cup \GC \cup \Omegatwo\,,
\]
 with $\Omegaone$ and $\Omegatwo$
disjoint Lipschitz subdomains and $\GC$ their common boundary, which
represents a prescribed delamination $(d{-}1)$-dimensional surface.
We  denote by $\nu$ the outward unit normal to $\partial \Omega$,
and by $\norm$
 the unit normal to $\GC$, which we consider
oriented from $\Omegatwo$ to $\Omegaone$. Moreover, given  $v \in
W^{1,2} (\Omega {\setminus} \GC)$, $v^+$ (respectively, $v^-$) shall
signify the restriction of $v$ to $\Omegaone$ (to $\Omegatwo$,
resp.).
 We further suppose  that
\[
\partial \Omega = \Gdir\cup \Gnew\,,
\]
with $\Gdir$ and $\Gnew $ open subsets in the relative topology of
$\partial\Omega$,  disjoint one from each other and each of them
with a smooth boundary.

As {\it state} variables, inside $\Omega$ we have  the {\it
displacement} $u:\Omega{\setminus}\GC \to\R^d$ and the absolute
temperature $\theta:\Omega{\setminus}\GC\to(0,+\infty)$, while on
the contact boundary we consider a {\it delamination variable}
$z:\GC\to[0,1]$, having the meaning of the integrity fraction of the
adhesive. Namely, $z=1$ (respectively $z=0$) means that the adhesive
has full (resp.~no) integrity. We denote by
\[
\text{$\JUMP{u}{}= u^+|_{\GC} - u^-|_{\GC}\ $= the jump of $u$ across
$\GC$}.
\]
Furthermore, we shall denote by $T=T(u,v,\theta)$  the traction
stress on some $(d{-}1)$-dimensional surface $\Gamma$
 (later, we shall take either $\Gamma= \GC$ or
$\Gamma=\Gnew$), i.e.
\begin{align}\label{T-stress}
T(u,v,\theta):=\sigma\big|_{\Gamma}\nu\,,\quad \text{with} \quad
\sigma:=\bbD e(\DT{u})+\bbC\big(e(u){-}\bbE\theta\big),
\end{align}
where of course we take as $\nu$ the unit normal  $\norm$ to $\GC$,
if $\Gamma=\GC$. In~\eqref{T-stress}, $\sigma$ is the stress
(assuming Kelvin-Voigt's rheology and {\it thermal expansion},
see~\eqref{kelvin-voigt} later on).

To describe various general situations in a unified and simple way,
we introduce
\[
\text{ a closed, convex cone $\cone(x)\subset\R^d$, possibly
depending on $x\in\GC$,}
\]
and assume the boundary conditions
on $\GC$ in the complementarity form as
\begin{align}\label{b.c.}
\left.\begin{array}{rll}
\JUMP{u}{}&\GE&0,\ \ \\[0em]
T(u,\DT{u},\theta)&\GEstar&0,\ \ \\[0em]
T(u,\DT{u},\theta){\cdot}\JUMP{u}{}&=&0
\end{array}\right\}
\ \text{ on }\GC.
\end{align}
In \eqref{b.c.}, $\GE $ is the ordering induced by the multivalued,
cone-valued mapping $\cone: \GC \rightrightarrows \R^d$,  in the
sense that, for $v_1,v_2: \GC \to \R^d$,
\begin{equation}\label{e:ordering}
\text{$v_1\GE v_2\ $ if and only if $v_1(x){-}v_2(x)\in \cone(x)$
for a.a.~$x\in\GC$.}
\end{equation}
Likewise,  $\GEstar$ is the dual ordering induced by the negative
polar cone to $\cone$, in the sense that,
for $\zeta_1,\zeta_2: \GC \to \R^d$,
\[
\text{$\zeta_1 \GEstar\zeta_2\ $ if
and only if $\zeta_1(x){\cdot}v\ge\zeta_2(x){\cdot}v$ for all
$v\in \cone(x)$, for a.a.~$x\in\GC$.}
\]

Possible  choices for the cone-valued  mapping $\cone:\GC
\rightrightarrows \R^d$ are
\begin{subequations}\label{CC}\begin{align}\label{CC1}
&&&&&&&&&\cone(x)=\cone=\R^d \quad \foraa\, x \in \GC,&&\text{or
}&&&&&&&&
\\\label{CC2}
&&&&&&&&&\cone(x)=\{v\in\R^d;\ v{\cdot}\norm(x)\ge0\}\quad \foraa\, x
\in \GC,&&\text{or }&&&&&&&&
\\\label{CC3}
&&&&&&&&&\cone(x)=\{v\in\R^d;\ v{\cdot}\norm(x)=0\}\quad \foraa\, x
\in \GC.&&&&&&&&
\end{align}\end{subequations}
In the first case \eqref{CC1}, the second of boundary conditions
\eqref{b.c.} translates into $T(u,\DT{u},\theta) =0$ on $\GC$, while
no constraint on $\JUMP{u}{}$ is imposed. Thus, \eqref{CC1}
 allows for no interaction of
the  bodies $\Omegaone$ and $\Omegatwo$ after a complete
delamination. In fact, this model is very simplified because
it does not prevent
possible interpenetration and delamination can be thus triggered,
rather unphysically,   by mere compression. Nevertheless, a model
like this  may be feasible in some situations. In this connection,
let us point out that the interpenetration after developed cracks is
neglected in several crack models used in mathematical literature
(as e.g.~\cite{ChGiPo08CIBM,DMFraToa05,GiaPon06GCAS}), too. The case
\eqref{CC2} yields the standard model of unilateral frictionless
Signorini contact in the normal displacement at $x\in\GC$. The last
case \eqref{CC3} prescribes  the normal jump of the displacement,
variable at $x\in\GC$, to zero. Thus, it only allows   for a
tangential slip along $\GC$. This may be a relevant model under high
pressure, when no cavity of $\GC$ can be expected anyhow. Such a
situation occurs, e.g., on lithospheric faults deep under the earth
surface. Note that, both in \eqref{CC1} and in \eqref{CC3},
$\cone(x)$ is a linear manifold for a.a. $x \in \GC$. As we shall
see later, this feature may allow for some special benefits.
\paragraph{Classical formulation of the adhesive contact problem}
Beside the force equilibrium, coupled with the heat equation inside
$\Omega{\setminus}\GC$ and  supplemented  with standard boundary
conditions, we have two complementarity problems on $\GC$, namely
\begin{subequations}\label{eq6:adhes-class-form}\begin{align}
\label{eq6:adhes-class-form1}& \varrho\DDT{u} -\mathrm{div}\big(\bbD
e(\DT{u})+\bbC\big(e(u){-}\bbE\theta\big)\big)= \FRM &\text{in  }
Q{\setminus}\SC,
\\
\label{eq6:adhes-class-form1bis}
 & c_{\rm v}(\theta)\DT{\theta}
-\mathrm{div}\big(\bbK(e(u),\theta)\nabla\theta\big)= \bbD
e(\DT{u}){:} e(\DT{u}) +\theta\bbC\bbE{:} e(\DT{u})+ \GRM
 &\text{in }Q{\setminus}\SC,
\\
\label{eq6:adhes-class-form2} &u=0 &\text{on }\Sdir,\hspace{1.2em}
\\\label{eq6:adhes-class-form3-bis}
&T(u,\DT{u},\theta)=\fRM&\text{on }\Snew,\hspace{1.2em}
\\\label{eq6:adhes-class-form3}
&(\bbK(e(u),\theta)\nabla\theta)\nu=\gRM&\text{on }\Sigma,\hspace{1.7em}
\\\label{adhes-form-d1}
& \JUMP{\bbD e(\DT{u})+\bbC (e(u){-}\bbE\theta)}{}\norm=0 &\text{on
}\SC,\hspace{1.2em}
\\\label{adhes-form-d2}
&
\JUMP{u}{}\GE0  &\text{on }\SC,\hspace{1.2em}
\\\label{adhes-form-d3}
&
T(u,\DT{u},\theta)+\kappa z\JUMP{u}{}\GEstar0
&\text{on }\SC,\hspace{1.2em}
\\
\label{adhes-form-d4}
&
\big(T(u,\DT{u},\theta)+\kappa z\JUMP{u}{}\big){\cdot}\JUMP{u}{}=0
&\text{on }\SC,\hspace{1.2em}
\end{align}
\begin{align}
\label{adhes-form-d6}
&\DT{z}\le0&\text{on }\SC,\hspace{1.2em}
\\\label{adhes-form-d7}
& d\le a_1+a_0 &\text{on }\SC,\hspace{1.2em}
\\\label{adhes-form-d8-bis}
& \DT{z} \left( d - a_0-a_1\right) =0
 &\text{on }\SC,\hspace{1.2em}
\\\label{adhes-form-d8}
& d\in N_{[0,1]}(z)+
\mbox{$\frac12$}\kappa\big|\JUMP{u}{}\big|^2
 &\text{on }\SC,\hspace{1.2em}
\\\label{adhes-form-d9}
& \frac12\big(\bbK(e(u),\theta)\nabla\theta|_{\GC}^+
+\bbK(e(u),\theta)\nabla\theta|_{\GC}^-\big){\cdot}\norm
+\het(\JUMP{u}{},z)\JUMP{\theta}{}=0
 &\text{on }\SC,\hspace{1.2em}
\\
\label{adhes-form-d10} &
\JUMP{\bbK(e(u),\theta)\nabla\theta}{}{\cdot}\norm ={-} a_1\DT{z}
&\text{on }\SC,\hspace{1.2em}
\end{align}
\end{subequations}
   where   we have used the
notation
\begin{displaymath}
Q:=(0,T)  \times \Omega, \quad \Sigma : = (0,T)  \times \partial
\Omega, \quad \SC:= (0,T)  \times\GC , \quad \Sdir:= (0,T)
\times \Gdir , \quad \Snew:= (0,T) \times\Gnew,
\end{displaymath}
$T>0$ being a fixed time horizon. In  \eqref{eq6:adhes-class-form},
$\FRM:Q\to\R^d$ is the applied bulk force, $w_{\rm D}:\Sdir\to\R^d$
the prescribed time-dependent boundary displacement,
$\fRM:\Snew\to\R^d$ the applied traction, while
$\GRM:Q\to\R$ and $\gRM:\Sigma\to\R$ are some external heat sources.
In addition,
\begin{equation}\label{posit} \bbC,\,\bbD:\R_{\mathrm{sym}}^{d\times d}\to
\R_{\mathrm{sym}}^{d\times d}\quad\text{ are $4$th-order
positive definite and symmetric tensors,}
\end{equation}
(i.e. $\bbC_{ijkl}= \bbC_{jikl}= \bbC_{klij}$, and the same for
$\bbD$), $\bbK=\bbK(e,\theta)$ is the positive definite matrix of the
heat conduction coefficients, and $\bbE\in\R^{d\times d}$ is a
matrix of thermal-expansion coefficients. Furthermore, the constant
$\kappa>0$ phenomenologically describes  the elastic response of the
adhesive. The complementarity
problem \eqref{adhes-form-d2}--\eqref{adhes-form-d4} describes
general, possibly unilateral (depending on the choice of the mapping
$\cone:\GC \rightrightarrows \R^d$) contact, whereas the
 adhesive contact results from  the complementarity
conditions~\eqref{adhes-form-d6}--\eqref{adhes-form-d8}.
In~\eqref{adhes-form-d7}
 the coefficient $a_0$ (resp.~$a_1$) is the
phenomenological specific energy (per area) which is stored
(resp.~dissipated) by disintegrating the adhesive. The overall
activation energy to trigger the debonding process in the adhesive
is then $a_0+a_1$. Note that the term $\kappa\big|\JUMP{u}{}\big|^2$
in~\eqref{adhes-form-d8} is in fact a penalization of the
 delamination condition $z\JUMP{u}{}=0$,
cf. with the brittle delamination model
\eqref{eq6:delam-class-form}.
 Moreover,
 $N_{[0,1]}$ denotes the normal cone to the interval $[0,1]$,
i.e. the subdifferential in
 the sense of convex analysis of the indicator function $\ind_{[0,1]}$ of
 $[0,1]$.
Finally, $\het=\het(x,\JUMP{u}{},z)\ge0$ is a phenomenological
\emph{heat-transfer coefficient}, determining the linear heat
convection through $\GC$. We shall suppose that $\het$ depends
affinely  on the delamination variable $z$, cf.~\eqref{eta-affine}
below.

\section{Thermodynamics of the model and various remarks}\label{sec-th}
         ~~~~~~~~~~~~~~~~~~~~~~~~~~~~~~~~~~~~~~~~~~~~~~~~~~

Let us briefly present the {\it thermodynamics} of the boundary-value problem
\eqref{eq6:adhes-class-form}.
The underlying overall Helmholtz {\it free energy}
$\Psi:\R^d\times\R\times (0,+\infty)\to\R$ has a bulk and a surface
part, i.e.
\begin{subequations}
\begin{align}\label{Phi-epos}
\Psi(u,z,\theta)= \int_{\Omega {\setminus} \GC} \psi^{\mathrm{bulk}}
(e(u),\theta)\, \d x + \int_{\GC} \psi^{\mathrm{surf}}(\JUMP{u}{},
z)\, \d S\,,
\end{align}
with $\psi^{\mathrm{bulk}}$ and $\psi^{\mathrm{surf}}$, respectively
being the bulk and the contact-surface contributions to the specific
Helmholtz energy. One can identify
\begin{align}\nonumber
\psi^{\mathrm{bulk}} (e,\theta)&=
\frac12\bbC(e{-}\bbE\theta){:} (e{-}\bbE\theta)
-\frac{\theta^2}2\bbB{:} \bbE-\psi_0(\theta)
\\\label{psi-bulk}
&=\frac12\bbC e{:} e-\theta \bbB{:}e-\psi_0(\theta)
\quad\text{ with }\ \bbB:=\bbC\bbE.
\end{align}
Here $\frac12\bbC e{:} e$ is the {\it mechanical part of the
internal  energy in the bulk}, while $ -\psi_0(\theta)$ is the {\it
thermal part of the free energy}. Hereafter, we shall assume that
\begin{equation}
\label{psizero} \text{$\psi_0: (0,+\infty)\to\R$ a strictly convex
function.}
\end{equation}
The specific \emph{contact surface} energy
$\psi^{\mathrm{surf}}(\JUMP{u}{},z)$ is then
\begin{equation}
\label{psi-surf} \psi^{\mathrm{surf}} (v, z)=
\begin{cases}
\displaystyle{\int_{\GC}\left(\frac{\dela}2
z\big|\JUMP{u}{}\big|^2-a_0z\right)\,\d S}  & \mbox{if
$\JUMP{v}{}\GE0$ and $\ 0{\le}z{\le}1$
a.e.~on $\GC$},\\
+\infty&\text{otherwise}.\end{cases}
\end{equation}
\end{subequations}
The other underlying ingredient of the model is
the overall {\it dissipation rate} $\xi$, which also has   bulk and
 surface contributions and $\xi^{\mathrm{bulk}}$ and
$\xi^{\mathrm{surf}}$, namely:
\begin{align}\label{dissip-rate}
\Xi(\DT{e},\DT{z}):=\int_{\overline{\Omega}}\big[\xi(\DT{e},\DT{z})\big](\d x)=
\int_{\Omega{\setminus}\GC}\xi^{\mathrm{bulk}}(\DT{e})\, \d x
+\int_{\GC}\!\xi^{\mathrm{surf}}(\DT{z})\d S
\end{align}
where the specific dissipation rate $\xi(\DT{e},\DT{z})=
\xi^{\mathrm{bulk}}(\DT{e})\d x+\xi^{\mathrm{surf}}(\DT{z})\d S$ is
a measure in general, with absolutely continuous part determined by
the (pseudo)potential of viscous-type dissipative forces in the
bulk, and a possibly concentrating part, supported on $\GC$, i.e.
\[
\xi^{\mathrm{bulk}}(\DT{e})=2\zeta_2(\DT{e}),\qquad
\zeta_2(\DT{e}):=\frac12 \bbD\DT{e} {:} \DT{e}, \qquad
\xi^{\mathrm{surf}}(\DT{z})=\zeta_1(\DT{z})
:=\begin{cases} \displaystyle{a_1\big|\DT{z}\big|}&
\mbox{ if $\DT{z}\le0$ a.e. in $\GC$},\\
+\infty&\text{ otherwise},\end{cases}
\]
the latter term  representing the potential (and also the specific
dissipation rate) of the rate-independent delamination process on
the contact boundary $\GC$.

Standardly, one then defines the \emph{specific entropy}
$s=s(\theta,e)$ by the so-called Gibbs' relation $ \langle
s,\tilde\theta\rangle=-\Psi'_\theta(u,z,\theta;\tilde\theta)$ where
$\Psi'_\theta(u,z,\theta;\tilde\theta) $ is the directional
derivative of $\Psi$ at $(u,z,\theta)$ in the direction
$\tilde\theta$. This yields the entropy in the bulk as
\begin{align}\label{spec-entropy}
s=-\frac{\partial\psi^{\mathrm{bulk}}}{\partial\theta}(e(u),\theta)
=\bbB{:} e(u)+\psi_0'(\theta).
\end{align}
Further, we shall use the so-called {\it entropy equation}
\begin{align}\label{8-3-*}
\theta\DT{s}=\xi^{\mathrm{bulk}}(e(\DT{u}))-{\rm div}(j)+\GRM.
\end{align}
Substituting $\DT{s} =\bbB {:}
e(\DT{u})-\psi_0''(\theta)\DT{\theta}$, cf.~\eqref{spec-entropy},
into the entropy equation \eqref{8-3-*} yields
the {\it heat equation}
\begin{align}\label{heat-eq}
c_{\rm v}(\theta)\DT{\theta}+{\rm div}(j)=
2\zeta_2(e(\DT{u}))+\theta\bbB{:} e(\DT{u})+\GRM
\end{align}
with the {\it heat capacity}
\begin{align}\label{cap}
 c_{\rm v}(\theta)=\theta\psi_0''(\theta).
\end{align}
Hence, postulating the constitutive relation for the heat flux
\begin{align}\label{Fourier-law}
j:=-\bbK(e(u),\theta)\nabla\theta,
\end{align}
i.e.~{\it Fourier's law} in an anisotropic medium, one obtains to
the heat equation in the form \eqref{eq6:adhes-class-form1bis}.

Similar, but simpler thermodynamics can be seen also on the contact
boundary by involving $\psi^{\mathrm{surf}}$ and
$\xi^{\mathrm{surf}}$. As \eqref{psi-surf} is independent of
temperature, the ``boundary
entropy''=$-\frac{\partial}{\partial\theta}\psi^{\mathrm{surf}}$ is
simply zero, and the corresponding entropy equation reduces to
$0=\xi^{\mathrm{surf}}(\DT{z})-\JUMP{j}{}$ (as an analog of
\eqref{8-3-*}), which then results in \eqref{adhes-form-d10}.
Incorporating the analog of the phenomenological law
\eqref{Fourier-law}, we arrive at~\eqref{adhes-form-d9}.


\noindent \textbf{Momentum equation.}
 As in  {\it Kelvin-Voigt rheology}, the total stress
$\sigma$ is postulated as
\begin{align}
\label{kelvin-voigt} \sigma=
\frac{\partial}{\partial e} \zeta_2\big(e(\DT{u})\big)
+\frac{\partial}{\partial
e}\psi^{\mathrm{bulk}}\big(e(u),\theta\big),
\end{align}
which just gives $\sigma$ from \eqref{T-stress}. From {\it
Hamilton's principle}, generalized for dissipative systems as in
\cite{bedford} and with the specific kinetic energy $\frac12\rho
|\DT{u}|^2$, one then obtains the equilibrium
equation~\eqref{eq6:adhes-class-form1}.
 For later use,
we introduce the indicator functional $\ind_{\cone}$ associated with
$\cone : \GC \rightrightarrows \R^d$, defined on $L^2 (\GC;\R^d)$ by
\begin{equation}
\label{ind_D} \ind_{\cone}(v) = \int_{\GC}\ind_{\cone(x)}(v(x))\,\dd
S(x) \quad \text{for all $v \in L^2 (\GC;\R^d)$.}
\end{equation}
We point out that the complementarity conditions
\eqref{adhes-form-d2}--\eqref{adhes-form-d4}
 may be reformulated as
the subdifferential inclusion
\begin{equation}
\label{subdiff-adhes} \partial \ind_{\cone}(\JUMP{u}{}) +
T(u,\DT{u},\theta) +\kappa z \JUMP{u}{}\ni 0 \quad \text{in $\SC$,}
\end{equation}
featuring  the (convex analysis) subdifferential $\partial
\ind_{\cone}: L^2(\GC;\R^d) \rightrightarrows L^2(\GC;\R^d)$ of the
indicator functional $\ind_\cone$ introduced in \eqref{ind_D}.
\par\noindent\textbf{Evolution of the delamination parameter (a flow rule).}
Finally,
 we consider the following differential inclusion
for the inelastic evolution of the parameter $z$
\[
\partial \zeta_1 \left(\DT{z}\right) + \partial_z \psi^{\mathrm{surf}}
(u,z) \ni 0  \quad \text{in $\SC$,}
\]
which in the  adhesive case
 results in
 \begin{equation}
\label{e:reactivation-2}
\partial \ind_{(-\infty,0]}(\DT{z}) + \partial\ind_{[0,1]}(z) +
\mbox{$\frac12$}\kappa\big|\JUMP{u}{}\big|^2 -a_0-a_1 \ni 0 \quad
\text{in $\SC$.}
\end{equation}
It is immediate to check that \eqref{e:reactivation-2} is a
reformulation of~\eqref{adhes-form-d6}--\eqref{adhes-form-d8}.

The entropy equation (\ref{8-3-*}) is designed to balance the total
energy, i.e.~the sum of the kinetic energy integrated over $\Omega
{\setminus} \GC$ with the overall dissipated energy  (i.e., $\Xi$
from~\eqref{dissip-rate} integrated in time), and with the bulk {\it
internal energy}
\begin{align}\nonumber
e^{\mathrm{bulk}}(e,\theta):=\psi^{\mathrm{bulk}}+\theta s&=
\frac12\bbC e(u) {:} e(u)-\psi_0(\theta)
-\theta\bbB{:}\bbE+\theta(\bbB{:}\bbE+\psi_0'(\theta))
\\\label{ie}&=h(\theta) +\frac12\bbC e(u) {:} e(u)
\qquad\text{ with }\quad
h(\theta):=\theta\psi_0'(\theta)-\psi_0(\theta);
\end{align}
we convene to refer to $h$ as the \emph{enthalpy}, see
also~\cite[Sect.~2]{tr1}.
One can then derive the \emph{total energy balance}:
\begin{align}\nonumber
  &\frac{\d}{\d t}\!\!\!\!\ddd{\int_{\Omega{\setminus}\GC}\!\!
\varrho|\DT{u}|^2+\frac12\bbC e(u){:}e(u)+h(\theta)\,\d
x}{kinetic, elastic, and thermal energies}{}
\\&\label{total-energy}
\hspace{6em}
+ \!\!\! \!\!\!\ddd{\int_{\GC}\!\! \frac{\dela}2
z\big|\JUMP{u}{}\big|^2-a_0 z\,\d S}{mechanical energy}{in the
adhesive}\!\!\!
   =\!\!\!\!\ddd{\int_{\Omega}\!\GRM+\FRM{\cdot}\DT{u}\,\d x}{power of
bulk heat}{and mechanical load}\!\!\!\!
+\!\!\!\ddd{\int_{\Gamma}\gRM\,\d
S+\int_{\Gnew}\!\!\fRM{\cdot}\DT{u}\,\d S}{power of surface
heat}{and mechanical load}\!\!\!.
\end{align}

Assuming $\theta_0>0$, $\GRM\ge0$ a.e. in $\Omega$, and $\gRM\ge0$ a.e.~in
$\partial \Omega$, we can
rely on the fact that  $\theta>0$
a.e.~in $\Omega$ (proved later in Theorem~\ref{th:3.0})
and, using (\ref{8-3-*}), we derive the {\it Clausius-Duhem
inequality}:
\begin{equation}
\label{e:clausius-duhem}
\begin{aligned}
\frac{\d}{\d t}\int_\Omega\!s\,\d x= \int_\Omega\! \left(\frac{{\rm
div}\left(\!\bbK \nabla\theta\right)}{\theta}+
\frac{\GRM}{\theta}\right)\,\d x =\int_\Omega\!
\left(\frac{\bbK\nabla\theta{\cdot}\nabla\theta}{\theta^2}+\frac{\GRM}{\theta}\right)\,\d
x +\int_{\partial \Omega}\frac { \gRM }\theta\,\d S \ge0.
\end{aligned}
\end{equation}

\begin{remark}[Partly linearized ansatz]\upshape
An important feature is that, as a consequence of the partly
linearized ansatz \eq{psi-bulk}, the mechanical and thermal
variables are additively separated in \eq{spec-entropy}, which makes
$c_{\rm v}$ in \eq{heat-eq} independent of $u$, and thus makes
mathematical analysis easier.
\end{remark}

\begin{remark}[Non-homogeneous boundary conditions]\label{rem:Robin} \upshape
We could supplement \eqref{eq6:adhes-class-form1}  with
non-ho\-mo\-ge\-neous,
 Dirichlet boundary conditions on $\Gdir$,
i.e. impose
\begin{equation}
\label{non-homog}
 u=\omega_{\rm D}\quad  \text{on }\Sdir,
 \end{equation}
for some prescribed time-dependent loading $\omega_{\rm D}:[0,T] \to
H^{1/2}(\Gdir) $. The analysis we are going to perform in the case
of homogeneous Dirichlet conditions can be carried over to the case
of \eqref{non-homog} by arguing as in~\cite{tr-LS-CZ}, and thus
recurring to  the additive split $u(t)=\tilde{u}(t) +
u_{\mathrm{D}}(t)$ for almost all $t \in (0,T)$, with $\tilde{u}
:[0,T] \to W_{\Gdir}^{1,2}(\Omega {\setminus} \GC;\R^d)$ and
$u_{\mathrm{D}}:[0,T] \to  W^{1,2}(\Omega {\setminus} \GC;\R^d)$ an
extension of $\omega_{\mathrm{D}}$ to $\Omega$. It was observed
in~\cite{tr-LS-CZ} that, if $\overline{\GC}\cap
\overline{\Gdir}=\emptyset$, one can assume that
$\JUMP{u_{\mathrm{D}}(t)}{} = 0$, whence $\JUMP{u(t)}{}=
\JUMP{\tilde{u}(t)}{}$ for almost all $t \in (0,T)$. This allows for
a reformulation of the problem in terms of the unknown $\tilde{u}$,
hence   reducing the analysis to the case  with homogeneous
Dirichlet conditions.
\end{remark}

\begin{remark}[Heat-transfer contact conditions]\upshape
Note that the transient conditions~\eqref{adhes-form-d9}--\eqref{adhes-form-d10} on $\GC$ for the heat
equation can equivalently be written as two Robin-type conditions
\begin{subequations}\begin{align}
&\label{robina} \bbK(e(u),\theta)\nabla\theta|_{\GC}^+{\cdot}\norm+
\het(\JUMP{u}{},z)\theta|_{\GC}^+ =\het(\JUMP{u}{},z)\theta|_{\GC}^-
\minus \frac12 a_1\DT{z} \qquad \text{on $\SC$},
\\
\label{robinb}
&\bbK(e(u),\theta)\nabla\theta|_{\GC}^-{\cdot}\nu^{\mp}+
\het(\JUMP{u}{},z)\theta|_{\GC}^- =\het(\JUMP{u}{},z)\theta|_{\GC}^+
\minus \frac12 a_1\DT{z}\qquad \text{on $\SC$},
\end{align}\end{subequations}
where we have highlighted the unit normals $\norm$ from $\Omegatwo$
to $\Omegaone$ and $\nu^{\mp}$ from $\Omegaone$ to $\Omegatwo$. This
reveals that the heat generated by delamination is distributed  with
proportions $\frac12$ and $\frac12$ into the two subdomains adjacent
to $\GC$.
In principle, we could also consider a contribution from the
Stefan-Boltzmann radiation, which would then result in the condition
\begin{align}
\! \! \! \! \! \!\frac12\big(\bbK(e(u),\theta)\nabla\theta|_{\GC}^+
+\bbK(e(u),\theta)\nabla\theta|_{\GC}^-\big){\cdot}\norm
+\het_0(\JUMP{u}{},z)
\JUMP{\theta}{}+\het_1(\JUMP{u}{},z)\JUMP{\theta^4}{}=0\ \  \text{on
$\SC$},
\end{align}
with $\het_0(\JUMP{u}{},\cdot),\,\het_1(\JUMP{u}{},\cdot) \ge0$
affine. In fact, this would lead to a modification of the present
analysis  which is quite routine. Thus, we shall not scrutinize this
generalization here. Let us also remark that, in alternative to the
dependence of $\het$ on $\JUMP{u}{}$, a dependence on the normal
stress is sometimes considered, cf.~\cite{AKRS02ODTC}. However, it
seems difficult to adapt the present multidimensional analysis to
that case.
\end{remark}

\begin{remark}[Elastic response in adhesive]\label{rem-kappa}\upshape
One can easily imagine a positive-definite $d{\times d}$-matrix in
place of $\kappa$, which would more properly describe the
phenomenological elastic response of the adhesive. The related
analysis would be just a standard modification of the presented one.
\end{remark}
\begin{remark}[Griffith concept]\label{rem-brittle}\upshape
The classical concept of delamination is based on the {\it Griffith
criterion} \cite{Grif20PRFS}, phenomenologically prescribing  the
amount of energy $a$ (in J/m$^2$, in 3-dimensional situations)
needed to de\-la\-mi\-na\-te the surface, independently  of the rate
of the process. The classical Griffith-type approach considers the
adhesive inelastic and one speaks of a {\it brittle delamination}.
Our adhesive contact problem can be viewed as a regularization of
this brittle delamination, and in fact makes mathematical analysis
and numerical implementation easier. It has its own interpretation
and many applications, and it is thus often considered as the
original problem,
cf.~\cite{BBR1,BBR2,BBR3,BBR4,KoMiRo,Raous}. 
 In the
quasistatic isothermal case, it has been proved in \cite{tr-LS-CZ}
that the adhesive contact approximates the brittle delamination as
the elastic modulus $\kappa\to\infty$ in the framework of the
so-called {\it energetic solution} concept. Furthermore, in
\cite{MiRoTh??DDNE}, any energetic solution to the brittle
delamination has  been proved to  be of ``Griffith-type'' in the
sense that $z$ indeed takes either the  value $1$ or the value $0$.
\end{remark}

\begin{remark}[Engineering models]\upshape
In the engineering literature, the Griffith-type delamination on a
prescribed so-called ``weak surface'' is a quite accepted concept
(for example in the framework of the so-called Finite Fracture
Mechanics), although  it is often combined with the heuristically
devised stress criterion, which in some situations seems to provide
a better understanding of the initiation of the delamination
process, cf.~e.g.~\cite{Legu02STCC,Mant08ICOC}. The initiation of
the delamination process can sometimes be triggered by another crack
approaching the weak surface, according to the classical so-called
Cook-Gordon mechanism \cite{CooGor64MCCP}, which has been confirmed
experimentally. The present form  of the activation energy $a(x)$
may typically correspond to crack growth in a pure fracture mode
(e.g.~Mode I). Nonetheless, it is believed that the present approach
can be extended to a generalization of this form, in order to cover
more complex phenomenological engineering models, working with the
so-called ``fracture mode mixity" which reflects the character of
the load (the ratio of its shear and normal components) on the crack
tip.
\end{remark}

\begin{remark}[Constant heat capacity]\upshape
The special case $\psi_0(\theta)=c_0\theta{\rm
ln}(\theta/\theta_0)$, with $c_0>0$ and $\theta_0>0$ constant, would
give $c_{\rm v}(\theta)=c_0$ in (\ref{cap}).  However, this case is
not within the scope of our analysis, since   $c_{\rm v}(\cdot)$
does not  have a \emph{compatible
 growth}, cf.~\eqref{30b}.
\end{remark}


\begin{remark}[Brittle delamination model]\label{rem:brittle}\upshape
Let us now briefly comment on the  model for brittle
delamination with thermal effects which would result from the above
derivation. As in the case of adhesive contact, we focus on its
classical formulation, which couples the momentum equilibrium
equation \eqref{eq6:adhes-class-form1}, the heat equation
\eqref{eq6:adhes-class-form1bis}, the boundary conditions
\eqref{eq6:adhes-class-form2}--\eqref{adhes-form-d1} and
\eqref{adhes-form-d9}--\eqref{adhes-form-d10} with the   two
following complementarity problems on $\GC$:
\begin{subequations}\label{eq6:delam-class-form}\begin{align}
\label{class-form-d2} & \JUMP{u}{}\GE0  &\text{on }\SC,
\\\label{class-form-d3}
&
T(u,\DT{u},\theta)\GEstar0\qquad\text{wherever $z(\cdot)=0$}
&\text{on }\SC,
\\\label{class-form-d4}
&
T(u,\DT{u},\theta){\cdot}\JUMP{u}{}=0 &\text{on }\SC,
\\\label{class-form-d5}
& z\JUMP{u}{}=0 &\text{on }\SC,
\\\label{class-form-d6}
&\DT{z}\le0&\text{on }\SC,
\\\label{class-form-d7}
& d\le a_1+a_0 &\text{on }\SC,
\\
\label{class-form-d8-bis} & \DT{z} (d-a_0-a_1) =0 &  \text{on }\SC,
\\\label{class-form-d8}
& d\in N_{[0,1]}(z)+
\partial_z \indj \big(\JUMP{u}{},z\big)
 &\text{on }\SC.
   \end{align}
   \end{subequations}
Indeed, \eqref{class-form-d2}--\eqref{class-form-d4} and
\eqref{class-form-d5}--\eqref{class-form-d8} respectively correspond
to (possibly unilateral) contact and activated delamination.
 Note that the penalization terms $\delam
\JUMP{u}{}$ in \eqref{adhes-form-d3}--\eqref{adhes-form-d4}  and
$\delam/2 |\JUMP{u}{}|^2$ in
 \eqref{adhes-form-d8} are
no longer present, and  the delamination constraint
\eqref{class-form-d5} is enforced by the second subdifferential
operator in \eqref{class-form-d8},
 featuring  the indicator function
\begin{equation}
\label{ind-noncvx}
 \indj(v,z)=\ind_{\{vz=0\}}, \ \  \text{i.e.} \ \
 \indj(v,z)=\left\{ \begin{array}{lll}
 0 & \text{if $vz=0$,}
 \\
 +\infty & \text{otherwise.}
 \end{array}
 \right.
 \end{equation}
  In  fact, as function of the two variables $v$ and $z$   $\indj$ is nonconvex, but
separately convex. Hence, the subdifferentials of the convex
functions $\indj(\cdot,z)$ and $\indj(v,\cdot)$ are well-defined,
and in particular $\partial_zJ(\JUMP{u}{},z)$ is given by
\[
\partial_zJ\big(\JUMP{u}{},z\big)
=\left\{
\begin{array}{cl}\emptyset&\text{ if $z \neq 0$ and}\
\JUMP{u}{}\ne0,
\\
0 &\text{ if $z=0$  and $\JUMP{u}{}\ne0$},
\\
\R &\mbox{ if $z\neq 0$}\ .
\end{array}
\right.
\]

As we have already mentioned,  existence for the (global) energetic
formulation of the brittle delamination problem in the isothermal
quasistatic case has been proved in \cite{tr-LS-CZ}.  In contrast,
the analysis of the corresponding thermomechanical model given by
(\ref{eq6:adhes-class-form}a-f,n-o)--(\ref{eq6:delam-class-form}) is
for the moment being an open problem.  The main difficulties
attached to this problem are  related to the presence of two
multivalued operators in \eqref{class-form-d8}, and in particular to
the \emph{essentially nonconvex} character of the nonlinearity
\eqref{ind-noncvx}.

However, taking into account \eqref{class-form-d8}, we clearly
identify a drawback of the differential
formulation~(\ref{eq6:delam-class-form}d-h) of brittle delamination.
Indeed, in this framework
 any driving tendency towards delamination is
smeared  out if $0<z<1$,  because then the driving force is
$d=0<a_0{+}a_1$,  and necessarily, by \eqref{class-form-d8-bis}, we
have $\DT z=0$.  The adhesive contact problem shows a similar
behaviour if $\kappa\to\infty$.
\end{remark}

\section{Enthalpy transformation and energetic solution}\label{sec-ent-trans}

Throughout the paper,  we shall adopt the notation
\[
\begin{aligned}
&  W_{\Gdir}^{1,2}(\Omega {\setminus} \GC ;\R^d)=\big\{ \testu \in
W^{1,2}(\Omega {\setminus} \GC ;\R^d)\, : \ \ \testu =0 \ \ \text{on
$\Gdir$}\big\}\,,
\\
& W_{\GC}^{1,2}(\Omega {\setminus} \GC ;\R^d)=\big\{ \testu \in
W^{1,2}(\Omega {\setminus} \GC ;\R^d)\, : \ \ \testu =0 \ \ \text{on
$\GC$}\big\}\,.
\end{aligned}
\]
Furthermore, in the case  $\cone(x)$ is  a linear subspace of $\R^d$
for almost all $x \in \GC$, we shall use the notation
\[
\begin{aligned}
 W_{\cone}^{1,2}(\Omega {\setminus} \GC ;\R^d)=\big\{ \testu \in
W^{1,2}(\Omega {\setminus} \GC;\R^d)\, : \ \ \JUMP{\testu(x)}{} \in
\cone(x)\ \ \foraa\, x \in \GC\big\}.
\end{aligned}
\]
 We shall also
extensively exploit that, for  $d\le3$,
\begin{equation}
\label{e:contemb}
\begin{aligned}
& W^{1,2}(\Omega) \subset L^p(\Omega) \text{ continuously for} \  1
\leq p \leq 6,
\\
 &  u\mapsto u|_\Gamma:W^{1,2}(\Omega) \to H^{1/2}
(\Gamma) \subset L^m (\Gamma)\
\begin{cases}
\text{ continuously for} \  1 \leq m \leq 4,
\\
\text{ compactly for} \  1 \leq m < 4,
\end{cases}
\end{aligned}
\end{equation}
with  $\Gamma = \partial \Omega$, or $\Gamma = \GC$, or $\Gamma =
\Gnew$.  The same embeddings hold for  the  Sobolev space
$W^{1,2}(\Omega;\R^d) $ of vector-valued functions.
 Finally, we
shall denote by   $\pairing{}{}{\cdot}{\cdot}$ the duality pairing
between  the spaces $W_{\Gdir}^{1,2}(\Omega{\setminus} \GC;\R^d)^*$
and
$W_{\Gdir}^{1,2}(\Omega{\setminus} \GC;\R^d)$. 

The analysis of the nonlinear heat
equation~\eqref{eq6:adhes-class-form1bis}, featuring the quadratic
coupling terms with the momentum balance
equation~\eqref{eq6:adhes-class-form1}, calls for rather
sophisticated techniques and suitable working assumptions. In
particular, one may impose some  conditions either on the growth of
$\mathbb K(e,\cdot)$ (cf., e.g., \cite{fpr09} for the analysis of a
similar nonlinear heat equation in some phase transition model), or
on the growth of $c_{\rm v}$ (cf., e.g., \cite{tr0,tr1} and, more
specifically, \cite[Sect.5.4.2]{eck} for contact problems in
thermo-viscoelasticity). Under the latter kind of  assumptions, the
Galerkin approximation method for proving existence of solutions
could serve quite effectively, cf.~\cite{tr0}.

On the other hand, system~\eqref{eq6:adhes-class-form} hosts the
delamination rate-independent process on $\GC$. Hence, the Rothe
method (i.e.~the implicit discretization in time) seems more natural
for the analysis, see e.g.~\cite{MieThe04RIHM,Miel05ERIS,MieFra06}.
In turn, the nonlinearity  $c_{\rm v}(\cdot)$ makes it technically
difficult to implement such a discretization method. This problem
can be circumvented by rewriting the original PDE
system~\eqref{eq6:adhes-class-form} in terms of the enthalpy,
instead of the temperature, as e.g.~in \cite{tr0}.

Namely, we introduce the  so-called \emph{enthalpy transformation}, setting
\begin{align}\label{hat-c}
\w=\ent_0(\theta):=\int_0^\theta\!\!c_{\rm v}(r)\,\d r.
\end{align}
Thus, $\ent_0$ is a primitive function of $c_{\mathrm v}$,
normalized in such a way that $\ent_0(0)=0$. In view of  \eqref{cap}
and~\eqref{ie}, we have
\begin{align}
\ent'(\theta)=(\theta\psi_0'(\theta)-\psi_0(\theta))'=\theta\psi_0''(\theta)
+\psi_0'(\theta)-\psi_0'(\theta)=\theta\psi_0''(\theta)=c_{\mathrm
v}(\theta) =\ent_0'(\theta),
\end{align}
hence $\ent_0$ differs from $h$ just by a constant, namely
$\psi_0(0)$. Furthermore,
 thanks to~\eqref{psizero},
 $c_{\rm v}$ is strictly positive and hence $h_0$ is strictly
increasing. Thus, we are entitled to
 define
\begin{align}\label{K-T}
\Theta(\w):=
\left\{\begin{array}{ll}
\ent_0^{-1}(\w) &\text{if }\w\ge0,
\\
0              &\text{if }\w<0,\end{array}\right. \qquad\quad
\mathcal{K}(e,\w):= \frac{\mathbb K(e,\Theta(\w))}{c_{\mathrm
v}(\Theta(\w))},
\end{align}
where $\ent_0^{-1}$ here denotes the inverse function to $\ent$.

Taking into account \eqref{K-T}, as well as the \emph{subdifferential
reformulations} \eqref{subdiff-adhes} and \eqref{e:reactivation-2}
of the complementarity problems
\eqref{adhes-form-d2}--\eqref{adhes-form-d4} and
\eqref{adhes-form-d6}--\eqref{adhes-form-d8}, respectively,  the PDE
system~\eqref{eq6:adhes-class-form} turns into
\begin{subequations}\label{eqsystem}\begin{align}
\label{eq6:delam-class-trans1} &\left.\begin{array}{ll}
&\hspace{-1.6em} \varrho\DDT{u} -\mathrm{div}\big(\bbD
e(\DT{u})+\bbC e(u)-\bbB\Theta(\w)\big)=\FRM
\\[.3em]
 &\hspace{-1.6em}
\DT{\w}-\mathrm{div}\big(\mathcal{K}(e(u),\w)\nabla\w\big)= \bbD
e(\DT{u}){:} e(\DT{u})+\Theta(\w)\bbB {:} e(\DT{u})+\GRM
\end{array}\hspace{2.9em}\right\}
&&\text{in }Q{{\setminus}}\SC,
\\ \label{eq6:delam-class-trans2}
&u=0
&&\text{on }\Sdir,&&
\\
\label{eq6:delam-class-trans3} &(\mathcal{K}(e(u),\w)\nabla\w)\nu=f
&&\text{on }\Snew,&&
\\\label{eq6:delam-class-trans3-bis} &\mathcal{T}(u,\DT{u},\w)=g
&&\text{on }\Sigma,&&
\\\label{class-form-d-trans}
&\left.\begin{array}{ll}
&\hspace{-1.6em}
  \JUMP{\bbD e(\DT{u}){+}\bbC e(u){-}\bbB\Theta(\w)}{}\norm=0
\\[.3em]
&\hspace{-1.6em}
\partial \ind_{\cone}(\JUMP{u}{}) +
T(u,\DT{u},\theta) + \delam z \JUMP{u}{}\ni 0
\\[.3em]
&\hspace{-1.6em}
\partial \ind_{(-\infty,0]}(\DT{z}) +
\partial\ind_{[0,1]}(z) +
\mbox{$\frac12$}\kappa\big|\JUMP{u}{}\big|^2 -a_0-a_1 \ni 0
\\[.3em]
&\hspace{-1.6em} \frac12\big(\mathcal{K}(e(u),\w)\nabla\w|_{\GC}^+
\!+\mathcal{K}(e(u),\w)\nabla\w|_{\GC}^-\big){\cdot}\norm +
\eta(\JUMP{u}{},z)\JUMP{\Theta(\w)}{}=0
\\[.3em]
&\hspace{-1.6em} \JUMP{\mathcal{K}(e(u),\w)\nabla\w}{}{\cdot}\norm
=\minus a_1\DT{z}
  \end{array}\hspace{-.3em}\right\}
&&\text{on }\SC,&&
\intertext{where}
&\mathcal{T}(u,v,\w):= T(u,v,\Theta(\w))= \big[\bbD
e(v)+\bbC e(u)-\bbB\Theta(\w)\big]\big|_{\Gamma}\nu
&&&&
\end{align}\end{subequations}
where again  we take as $\nu$ the unit normal to $\GC$ $\norm$, if
$\Gamma=\GC$.
%
%
\par\noindent
\textbf{Data qualification.}  Hereafter, the problem data $\FRM$,
$\GRM$, $\fRM$, and $\gRM$ shall be qualified by
\begin{subequations}\label{hypo-data}
\begin{align}
  \label{eFFe1}
&\FRM \in L^1 (0,T; L^2 (\Omega; \R^d));
\\
& \label{eFFe2}\fRM \in \ \begin{cases} W^{1,1} (0,T;
L^{4/3}(\Gnew;\R^3))  &  \text{if $d=3$,}
\\
 W^{1,1} (0,T;
L^{1+\epsilon}(\Gnew;\R^2))  & \text{for some $\epsilon >0$, if
$d=2$;}
\end{cases}
\\
& \label{posg1} \GRM \in L^1 (Q),  \quad \GRM \geq 0 \
 \aein \ Q;
\\
& \label{posg2} \gRM \in L^1 (\Sigma), \quad  \gRM \geq 0  \  \aein
\ \Sigma\,.
\end{align}
\end{subequations}

The energetic formulation
associated with system~\eqref{eqsystem} hinges on the following
energy functional $\Phi$ (which is in fact the mechanical part of
the internal energy~\eqref{Phi-epos}),
 and   on the  dissipation potential  $\calD$
\begin{align}
\label{8-1-k}
 &
\begin{aligned}
 \Phi(u,z):= \int_{\Omega{\setminus}\GC} \frac12\bbC
e(u){:} e(u)\,\d x   + \ind_{\cone} (\JUMP{u}{})   +
\int_{\GC}\!\!\left( \frac{\dela}2 z\big|\JUMP{u}{}\big|^2+
\ind_{[0,1]}(z) -  a_0 z\right)\,\d S
\end{aligned}
\\&
\label{DISS} \calD\big(\tilde{z}{-}z):=
\begin{cases}
\displaystyle{\int_{\GC}\!a_1|\tilde{z}{-}z|\,\d S}  & \text{if
$\tilde{z} \leq
 z$ a.e. in $\GC$,}
\\[-.2em] + \infty & \text{otherwise.}
\end{cases}
\end{align}
 For notational convenience, we
also set for all $v \in L^2 (\Omega)$
\begin{displaymath}
T_\mathrm{kin}^\varrho (v):=\frac{1}{2}\int_\Omega\varrho\,|v|^2\,\d x.
\end{displaymath}

We are now in the position of introducing the notion of weak
solution to system \eqref{eqsystem} which shall be analyzed
throughout this paper. The reader is referred to
\cite[Prop.~3.2]{tr1} for some justification of the energetic
solution concept in the framework of general thermomechanical
rate-independent processes, in particular for the proof of the fact
that energetic solutions are also conventional weak solutions
whenever $\DT{z}$ is absolutely continuous.

\begin{definition}[Energetic solution of the adhesive contact problem]\label{def4}
\upshape
 Given a quadruple of initial data
$(u_0,\DT{u}_0,z_0,\theta_0)$ satisfying suitable conditions
(cf.~\eqref{hyp-init} later on),
 we call a triple $(u,z,\w)$ an \emph{energetic solution} to the
Cauchy problem for (the enthalpy reformulation of)
system~\eqref{eqsystem}
 if
\begin{subequations}
\label{reguu}
\begin{align}
\label{reguu1} & u \in
W^{1,2}(0,T;W_{\Gdir}^{1,2}(\Omega{{\setminus}}\GC;\R^d)),
\\
\label{reguu2}  & u \in
 W^{1,\infty}(0,T;L^2(\Omega;\R^d)) \qquad \text{if $\varrho>0$,}
\\
\label{reguz} &  z \in L^\infty (\SC)\, \cap\,
\BV([0,T];L^1(\GC))\,, \\
 \label{reguw}  & \left.
\!\!\!\!\begin{array}{ll} \w \in
L^r(0,T;W^{1,r}(\Omega{{\setminus}}\GC))  \,\cap\,
L^\infty(0,T;L^1(\Omega)) \\
\w \in \BV ([0,T]; W^{1,r'}(\Omega{\setminus} \GC)^*)
\end{array}
\right\} \quad
  \mbox{ for every
$1\le r<\frac{d+2}{d+1}$},
\end{align}
\end{subequations}
with $r'$ denoting the conjugate exponent $\frac{r}{r-1}$ of $r$,
and the triple $(u,z,\w)$  complies with:\\
 \ITEM{(i)}{(weak formulation of the) momentum inclusion, i.e.:}
\begin{align}
\label{constraints-delam}
 &  \JUMP{u}{}\GE0 \  \ \text{on $\SC$,} \ \ \text{and}
\\\nonumber
 &\int_\Omega\varrho\DT{u}(T){\cdot}\big(\testu(T){-}u(T)\big)\,\d x
+\int_Q\big(\bbD e(\DT{u})+\bbC
e(u)-\bbB\Theta(\w)\big){:}e(\testu{-}u)
-\varrho\DT{u}{\cdot}\big(\DT{\testu}{-}\DT{u}\big)
 \,
\d x\d t\\
&
\ \ +\int_{\SC}\!\!\kappa z\JUMP{u}{}{\cdot}\JUMP{\testu{-}u}{}\d S\d t
\ge\int_\Omega\!\varrho\DT{u}_0{\cdot}\big(\testu(0){-}u(0)\big)\d x
+\int_Q\!\FRM{\cdot}(\testu{-}u)\,\d x\d t
+\int_{\Snew}\!\!\fRM{\cdot}(\testu{-}u)\,\d S\d t
\label{e:weak-momentum-variational}
\end{align}
\ITEM{}{for all $\testu$ in
$L^2(0,T;W_{\Gdir}^{1,2}(\Omega{\setminus}\GC;\R^d))$ with
$\JUMP{\testu}{}\GE0$  on~$\SC$ and, if $\varrho>0$, also in
$W^{1,2}(0,T;L^2(\Omega;\R^d))$,}
 \ITEM{(ii)}{total energy inequality}
\begin{align}\nonumber
 T_\mathrm{kin}^{\varrho}\big(\DT{u}(T)\big)
+\Phi\big(u(T),z(T)\big)+\int_\Omega\w(T)\,\d x   \leq
T_\mathrm{kin}^{\varrho}\big(\DT{u}_0\big) +\Phi\big(u_0,z_0\big)
+\int_\Omega\w_0\,\d x
\\
+\int_{Q}\FRM{\cdot}\DT{u}\,\d x \d t
+\int_{\Snew}\fRM{\cdot}\DT{u}\,\d S\d t
+\int_{Q}\GRM\,\d x \d t+\int_{\Sigma}
\gRM \, \d S  \d t
 \label{total-energy-brittle}
\end{align}
 \ITEM{(iii)}{
semistability for a.a.~$t\in (0,T)$}
\begin{align}
 \label{semistab} &\forall \tilde{z}\in
L^\infty(\GC):\qquad \Phi\big(u(t),z(t)\big)\le\Phi\big(u(t),\tilde
z\big) +\calD\big(\tilde z-z(t)\big)
\end{align}
 \ITEM{(iv)}{ (weak formulation of the) enthalpy equation:}
\begin{align}\nonumber
& \int_\Omega\!\w(T)  \testw(T)\,\d x +
\int_Q\mathcal{K}(e(u),\w)\nabla\w{\cdot}\nabla \testw-\w\DT{
\testw}\,\d x\d t +\int_{\SC} \het(\JUMP{u}{},z)
\JUMP{\Theta(\w)}{}\JUMP{ \testw}{}\,\d S\d t
\\&\nonumber
=\int_Q\!\big(\bbD e(\DT{u}){:} e(\DT{u})+\Theta(\w)\bbB {:}
e(\DT{u})\big) \testw\,\d x\d t
+\int_{\overline{\SC}}\!\!\frac{
\testw|_{\GC}^+{+} \testw|_{\GC}^-}2\ \varmea(\d S\d t)
\\
&\label{weak-heat}
+\int_{Q}\GRM \testw\,\d x \d t
+\int_{\Sigma}\gRM\testw\,\d S \d t
+\int_\Omega \w_0\testw(0)\,\d x
\end{align}
\ITEM{}{for all $\testw \in \mathrm{C}^0
([0,T];W^{1,r'}(\Omega{\setminus} \GC)) \cap
W^{1,r'}(0,T;L^{r'}(\Omega))$, where $\w_0:= h_0 (\theta_0)$ and
$\varmea$ is a measure (=heat produced by rate-independent
dissipation) defined by prescribing its values for every  closed set
of the type
 $A:=[t_1,t_2]{\times}  C \subset[0,T]\times\overline\GC$ as}
\begin{align}
\label{meash} \varmea(A): =\begin{cases}
\displaystyle{\int_{C}\!a_1\big|z(t_1,x){-}z(t_2,x)\big|\,\d S}
&\text{if $z(\cdot,x)$ nonincreasing on $[t_1,t_2]$ for a.a.~$x{\in}C$},\\
+\infty&\text{elsewhere},\end{cases}
\end{align}
\ITEM{(v)}{and the remaining initial conditions
(\COLLL{in addition to} $\DT{u}(0)=\dot u_0$, already involved in
\eqref{e:weak-momentum-variational}), i.e.}
\begin{equation}
\label{init} u(0)=u_0 \quad \aein \ \Omega, \qquad z(0)=z_0 \quad
\aein \ \GC, \qquad \w(0)=\w_0 \quad \aein \ \Omega.
\end{equation}
\end{definition}


\begin{remark}
\upshape \label{radon}
The inequality \eqref{total-energy-brittle} is the integrated
(inequality) version of the total energy balance~\eqref{total-energy}.
It is immediate to check that, for every  closed set of the type
$A:=[t_1,t_2]{\times}  C \subset[0,T]\times\overline\GC$,
$\varmea(A)$ coincides with
$\mathrm{Var}_{\mathcal{R}}(z|_C;[t_1,t_2])$. Now, relying, e.g.,
on~\cite[Prop.~1.3.10,Thm.~1.5.6]{fremlin3}, one can verify that
formula \eqref{meash} indeed defines a non-negative Radon measure on
$\overline{\SC}$. \COL{Subtracting \eqref{weak-heat} tested by 1
from \eqref{total-energy-brittle} reveals the mechanical energy inequality:
\begin{align}\nonumber
 T_\mathrm{kin}^{\varrho}\big(\DT{u}(T)\big) &
+\Phi\big(u(T),z(T)\big)
+\int_Q\!\bbD e(\DT{u}){:} e(\DT{u})\,\d x\d t
+\mathrm{Var}_{\mathcal{R}}(z;[0,T])
\\ &
\le T_\mathrm{kin}^{\varrho}\big(\DT{u}_0\big) +\Phi\big(u_0,z_0\big)
+\int_{Q}\FRM{\cdot}\DT{u}-\Theta(\w)\bbB {:}
e(\DT{u})\,\d x \d t
+\int_{\Snew}\fRM{\cdot}\DT{u}\,\d S\d t.
\label{mech-energy}
\end{align}
In particular, $z(\cdot,x)$ must be nonincreasing on $[0,T]$ for
a.a.~$x \in \GC$, otherwise
$\mathrm{Var}_{\mathcal{R}}(z;[0,T])=\infty$ and \eqref{mech-energy}
cannot hold. }
\end{remark}
\section{Main results}\label{sec-main-res}
\par\noindent
\textbf{Assumptions.}
Hereafter, we shall denote by the symbols $C$,
$C'$ most of the (positive) constants occurring in calculations and
estimates.
We suppose that
\begin{subequations}
\label{C-K}\begin{align}
 &
\label{30a}
 c_{\rm v}:[0,+\infty)\to\R^+\ \text{
continuous},
\\&
\label{30b} \exists\,\omega_1\ge\omega
>\frac{2d}{d{+}2},\ c_1\ge c_0>0\
\forall\theta\in\R^+:\quad c_0(1{+}\theta)^{\omega-1}\le c_{\rm
v}(\theta)\le c_1(1{+}\theta)^{\omega_1-1},
\\
\label{30c} &\mathbb K:\R^{d\times d}\times\R \to\R^{d\times d}\
\text{is  bounded, continuous, and }
\\&
\label{30dprimo}  \inf_{(e,\w,\xi)\in\R_{\rm sym}^{d\times
d}\times\R\times\R^d,\ |\xi|=1} \mathcal{K}(e,\w)\xi{:} \xi=
\mathsf{k}>0,
\end{align}
and that $\eta(x,v,\cdot)$ is a non-negative
 affine function of the
delamination parameter $z\in [0,1]$, i.e.
\begin{equation}
\label{eta-affine}
\begin{aligned}
\eta(x,v,z) &=\eta_1(x,v)z +\eta_0(x,v)
\ \text{ for $\eta_1,\eta_0:\GC{\times}\R^d\to\R^+$ Carath\'eodory
s.t.}
\\
& \exists\, C_\eta>0\,: \quad \forall\, (x,v) \in \GC{\times}\R^d \
\  |\eta_0(x,v)|+  |\eta_1(x,v)|  \leq C_\eta (|v|^{4/3} + 1);
\end{aligned}
\end{equation}
\end{subequations}
in fact, the above growth condition for the functions
$\eta_0(x,\cdot)$ and $\eta_1(x,\cdot)$ is not optimal and could be
slightly improved, as one can deduce from the proof of
Theorem~\ref{th:4.1} in Sect.~\ref{ss:4.4} later on.
 It is
immediate to deduce from~\eqref{30b} that
\begin{equation}
\label{growthTheta}\exists\, C_{\theta}^1,\, C_{\theta}^2>0\, \ \
\forall\, w \in \R^+\, : \quad C_{\theta}^1 (w^{1/\omega_1} -1) \leq
\Theta(w) \leq C_{\theta}^2 (w^{1/\omega} -1)\,.
\end{equation}
Moreover, it follows from~\eqref{30c}, \eqref{growthTheta},  and the
definition~\eqref{K-T} of $\mathcal{K}$ that
\begin{equation}
\label{growthKappa}\exists\, C_{\mathcal{K}}>0\, \ \ \forall\,\xi,\
\zeta \in \R^d\, : \quad \left| \mathcal{K}(e,\w)\xi{:}
\zeta\right| \leq C_{\mathcal{K}} |\xi||\zeta|\,.
\end{equation}
%
Finally,  we impose the following on the initial data
\begin{subequations}
\label{hyp-init}
\begin{align}
& \label{uzero} u_0 \in W_{\Gdir}^{1,2}(\Omega{\setminus}
\GC;\R^d)\,, \quad \JUMP{u_0}{} \GE 0 \ \ \text{on $\SC$,}
\\
& \label{vzero} \DT{u}_0 \in L^2(\Omega;\R^d) \quad\,\text{ if
$\varrho>0$}\,,
\\
& \label{zzero} z_0 \in L^\infty(\GC), \qquad 0 \leq z_0 \leq 1 \ \
\text{a.e. on}\, \GC\,,
\\
& \label{wzero} \theta_0 \in L^\omega (\Omega)\,,\qquad\ \ \theta_0
\geq 0 \ \ \aein\, \Omega\,.
\end{align}
\end{subequations}


\begin{theorem}[Existence for the adhesive contact problem]\label{th:3.0}
Let us assume \eqref{hypo-data}, \eqref{C-K},
\eqref{hyp-init} and\\
\ITEM{(i)}{if $\varrho=0$ (such a case is sometimes referred to as
\emph{quasistatic}), let also}
\begin{subequations}\label{hypo-data-bis}
\begin{align} \label{effebis} & \FRM \in  \ \begin{cases}
 W^{1,1} (0,T; L^{6/5}(\Omega;\R^3))   & \text{if
$d=3$,} \\   W^{1,1} (0,T; L^{1+\epsilon}(\Omega;\R^2))   &
\text{for  some $\epsilon >0$, if $d=2$};
\end{cases}
\\
\label{Dirichlet-part} & \mathscr{H}^{d-1} \left(\partial \Omegaone
\cap \Gdir \right)>0, \quad \mathscr{H}^{d-1} \left(\partial
\Omegatwo \cap \Gdir \right)>0,
\end{align}
\end{subequations}
\ITEM{}{with $\mathscr{H}^{d-1}$ denoting the $(d{-}1)$-dimensional
Hausdorff measure, or}
\ITEM{(ii)}{if $\varrho>0$, let also}
\begin{equation}
\label{varrho2}  \varrho>0 \ \ \text{and} \ \ \text{$\cone(x)$ is a
linear subspace of $\R^d$} \ \foraa\, x \in \GC.
\end{equation}

Then, there exists an energetic solution $(u,z,\w)$ to the adhesive
contact problem with the additional regularity
\begin{equation}\label{additional-cone}
u\in W^{2,2}(0,T;W_{\cone}^{1,2}(\Omega{\setminus} \GC;\R^d)^*)
\qquad \text{if $\varrho>0$.}
\end{equation}
Furthermore, in both cases $\varrho>0$ and $\varrho=0$,
the positivity of the initial temperature
\begin{align}\label{strict-pos}
\inf_{x\in\Omega}\theta_0=:\theta^*>0
\end{align}
implies $\inf_{(t,x)\in Q}\theta=\inf_{(t,x)\in
Q}\Theta(\w(t,x))>0$; in particular, $\theta$
is a.e.~positive on $Q$.
\end{theorem}
\noindent Theorem~\ref{th:3.0} shall be proved in Section~\ref{s:5}
by passing to the limit in some regularized problem (where the
contact conditions on $\GC$ are penalized), which we shall present
in Sect.~\ref{s:4}. In turn, existence for the latter problem shall
be proved in Section~\ref{ss:4.4} by passing to the limit in a
further approximation scheme, constructed in Section~\ref{s:new} by
a regularized semi-implicit time discretization.
\begin{remark}
\upshape \label{rem:math-diff} In Theorem~\ref{th:3.0} we
distinguish the cases $\varrho>0$ and $\varrho=0$, because in the
latter case we are able to prove existence for a far larger class of
cones yielding the unilateral constraint on the displacement, in
particular the Signorini conditions. This stems from the fact that
the analysis of the momentum equilibrium equation in which inertia
interacts with Signorini boundary conditions is remarkably
difficult. It has indeed been an open problem for a long time and
only very recently, in~\cite{petrov-schatzman2009,
petrov-schatzman2010}, existence results
 have emerged
for the
 dynamical viscoelastic equation with Signorini contact conditions,
 in the one- and three-dimensional case on unbounded domains. Such results have been proved
  with very sophisticated Fourier analysis techniques.
  In the one-dimensional framework of~\cite{petrov-schatzman2010}, it has  also been
   obtained that the solutions comply with
  the energy balance.
\end{remark}
\begin{remark}
\label{rem:relax-test-delam} \upshape
Under~\eqref{varrho2}, the qualification $\testu \in
W^{1,2}(0,T;L^2(\Omega;\R^d))$ for the test functions
in~\eqref{e:weak-momentum-variational} might be relaxed to
\begin{equation}
\label{relaxed-regu} \testu \in
W^{1,2}(0,T;W_{\cone}^{1,2}(\Omega{\setminus} \GC;\R^d)^*).
\end{equation}
Indeed, thanks to~\eqref{constraints-delam} and  to the linearity of
$\cone(x)$ for almost all $x \in \GC$, the function $u$
fulfilling~\eqref{e:weak-momentum-variational} is such that $ \DT{u}
\in L^2 (0,T;W_{\cone}^{1,2}(\Omega{\setminus} \GC;\R^d)).$ Now, the
spaces $L^2 (0,T;W_{\cone}^{1,2}(\Omega{\setminus} \GC;\R^d)^*)$ and
$L^2 (0,T;W_{\cone}^{1,2}(\Omega{\setminus} \GC;\R^d))$ are in
duality. Hence, \eqref{relaxed-regu} is sufficient to give meaning
to the term $\int_{\Omega} \DT{u}{\cdot}\DT{\testu} \, \d x$.
\COLLL{A similar extension will also apply  for the test functions
of \eqref{e:weak-momentum} below.}
\end{remark}


\section{Regularization}\label{s:4}


We shall approximate (the enthalpy-reformulation) of the adhesive contact
system~\eqref{eqsystem} by penalizing the contact condition $\JUMP{u}{}\GE0$.
This is a well-established routine in the analysis of contact problems, see
e.g.~\cite{eck}. We should emphasize that the penalized problems themselves
have their own practical usage because they allow, first, for combination of
inertia and the unilateral-type elastic contact condition and, second, for a
more physical interpretation of the coupling through the heat-transfer
coefficient, cf.~Remark~\ref{rem-heat-transfer} below.

Thus, we shall replace the subdifferential operator $\partial
\ind_{\cone}$ in the differential inclusion~\eqref{subdiff-adhes}
(equivalent to the complementarity problem
\eqref{adhes-form-d2}--\eqref{adhes-form-d4} on $\SC$), with
  its $\eps$-Yosida
regularization (see, e.g., \cite{attouch, brezis73, barbu76}). We
recall that the $\eps$-Yosida approximation of the indicator
functional $\ind_{\cone}$ is the lower semicontinuous, convex, and
Fr\'echet differentiable functional given by
\begin{align}\label{F-eps}
 \yosapp: L^2 (\GC;\R^d) \to [0,+\infty) \quad \text{given
 by} \quad
 \yosapp (v) = \frac1{2\eps} \min_{w \GE 0} \|
v-w\|_{L^2(\GC;\R^d)}^2;
\end{align}
cf. definition~\eqref{e:ordering} for the ordering $\GE$. We point
out that $\yosapp$ \textsc{Mosco}-converges to $\ind_\cone$ in the
 $L^2 (\GC;\R^d)$; see, e.g., \cite[\S~3.3]{attouch} for
the definition of \textsc{Mosco}-convergence
and~\cite[Thm.~3.66]{attouch} for the link with Yosida
regularizations. In particular, this entails that
\begin{equation}
\label{liminf-later} v_\eps \weakto v \quad \text{in $L^2
(\GC;\R^d)$} \ \Rightarrow \ \liminf_{\eps \to 0} \yosapp (v_\eps)
\geq \ind_{\cone}(v).
\end{equation}
The $\eps$-Yosida regularization of  $\partial \ind_{\cone}$ is the
Fr\'echet derivative  $\yosd: L^2(\GC;\R^d) \to L^2(\GC;\R^d)$ of
the functional  $\yosapp$.
 It is well known that
\begin{equation}
\label{yos-repre} \yosd =
\frac1\eps\left(\mathrm{Id}{-}\mathrm{P}_{\cone} \right),
\end{equation} where $\mathrm{Id} : L^2(\GC;\R^d) \to L^2(\GC;\R^d)$
is the identity operator and
 $\mathrm{P}_{\cone}: L^2(\GC;\R^d) \to L^2(\GC;\R^d)$ is the
projection  associated with the multivalued mapping $\cone: \GC
\rightrightarrows \R^d$.
For later use, we recall that, being
$\mathrm{P}_{\cone}$  a contraction on $L^2(\GC;\R^d) $, there holds
\begin{equation}
\label{e:contraction} \| \yosd(v)\|_{L^2 (\GC;\R^d)} \leq \frac2\eps
\| v \|_{L^2 (\GC;\R^d)} \quad \text{for all $v \in L^2
(\GC;\R^d)$.}
\end{equation}
\par\noindent
Hence, we shall consider the following regularized
conditions on $\GC$, where    \eqref{subdiff-adhes} is
approximated by {\it Yosida regularization}:
\begin{align}
\label{class-form-d-reg} &\left.\begin{array}{ll} &\hspace{-1.6em}
  \JUMP{\bbD e(\DT{u})+\bbC e(u){-}\bbB\Theta(\w)}{}\norm=0
\\[.3em]
&\hspace{-1.6em} \kappa z \JUMP{u}{} + \yosd(\JUMP{u}{})
  + T(u,\DT{u},\theta) =0
\\[.3em]
   &\hspace{-1.6em}
\partial \ind_{(-\infty,0]}(\DT{z}) + \partial\ind_{[0,1]}(z) +
\frac12\kappa\big|\JUMP{u}{}\big|^2 -a_0-a_1 \ni 0
\\[.3em]
&\hspace{-1.6em} \frac12\big(\mathcal{K}(e(u),\w)\nabla\w|_{\GC}^+
\!+\mathcal{K}(e(u),\w)\nabla\w|_{\GC}^-\big){\cdot}\norm +
\eta(\JUMP{u}{},z)\JUMP{\Theta(\w)}{}=0
\\[.3em]
&\hspace{-1.6em} \JUMP{\mathcal{K}(e(u),\w)\nabla\w}{}{\cdot}\norm =
\minus a_1\DT{z}
  \end{array}\hspace{-.5em}\right\} \quad \text{on $\SC$.}
\end{align}
The resulting regularized stored energy   is then
\begin{subequations}\label{Phi-eps}
\begin{align}
& \label{8-1eps}
\begin{aligned}
\Phi_{\eps}(u,z):= \int_{\Omega{\setminus}\GC} \frac12\bbC
e(u){:} e(u)\,\d x &  + \yosapp (\JUMP{u}{})
 + \int_{\GC}\!\!\left(\frac{\dela}2 z\big|\JUMP{u}{}\big|^2+
\ind_{[0,1]}(z) - a_0 z\right)\,\d S.
\end{aligned}
\end{align}
\end{subequations}
%
%
The main result of this section ensures the existence of
energetic solutions to the initial-boundary value problem for the
adhesive contact model supplemented with the regularized contact
conditions \eqref{class-form-d-reg}.

\begin{theorem}[Existence of energetic solutions to the
regularized problem]\label{th:4.1} Under assumptions
\eqref{hypo-data}, \eqref{hypo-data-bis},  \eqref{C-K} and
\eqref{hyp-init}, for every $\eps>0$ there exists a triple
$(u_\eps,z_\eps,\w_{\eps})$ as in~\eqref{reguu}, and such that, in
addition,
\[
u_\eps\in W^{2,2}(0,T; W_{\Gdir}^{1,2}(\Omega{\setminus} \GC;\R^d)^*)
\quad \text{if $\varrho>0$,}
\]
which solves the energetic formulation of the Cauchy problem for
system {\rm(\ref{eqsystem}a-d)} and~\eqref{class-form-d-reg}, namely
the initial conditions~\eqref{init} hold, as well as\\
\ITEM{(i)}{ the (weak formulation of the) momentum
equation:}
\begin{align}\nonumber
\!\!\!\!\!\!\!\!\!\!\!\!
\int_Q\!\big(\bbD e(\DT{u}_\eps){+}\bbC e(u_\eps){-}\bbB\Theta(\w_\eps)\big)
{:}e(\testu)-\varrho\DT{u}_\eps{\cdot}\DT{\testu}\,\d x\d t
+\int_{\SC}\!\!\left(\dela z_\eps\JUMP{u_\eps}{}{+}\yosd(\JUMP{u_\eps}{})
\right){\cdot}\JUMP{\testu}{} \, \d S\d t
\\\label{e:weak-momentum}
+\int_\Omega\varrho\DT u_\eps(T){\cdot} \testu(T)\,\d x
=\int_\Omega \varrho\DT u_0{\cdot} \testu(0)\,\d x + \int_{Q}\FRM{\cdot}
\testu \, \d x\d t +
\int_{\Snew}\!\fRM{\cdot}\testu\,\d S\d t
\end{align}
 \ITEM{}{ for all $\testu$ in $L^{2}
(0,T;W_{\Gdir}^{1,2}(\Omega{\setminus} \GC;\R^d))$ and, in the case
$\varrho>0$, in $W^{1,2}(0,T;L^2(\Omega;\R^d))$,}
 \ITEM{(ii)}{the total energy inequality~\eqref{total-energy-brittle} with
$\Phi_\eps$ {and $(u_\eps,z_\eps,\w_\eps)$} in place of $\Phi$ {and
$(u,z,\w)$},}
 \ITEM{(iii)}{ $z_\eps$ complies with
the semistability condition~\eqref{semistab}
for a.a.~$t\in (0,T)$,  again with $\Phi_{\eps}$ {and $(u_\eps,z_\eps)$,}
in place of $\Phi$ {and $(u,z)$}}
 \ITEM{(iv)}{ the weak formulation~\eqref{weak-heat} of the enthalpy equation.}
Furthermore, if~\eqref{strict-pos} holds, then
\begin{equation}\label{poswe}
\inf_{\eps>0,\,(t,x)\in Q}\w_\eps(t,x)>0\,.
\end{equation}
\end{theorem}

\begin{remark}[Signorini-contact case]\label{rem-heat-transfer}\upshape
Let us point out that, if $\cone:\GC \rightrightarrows \R^d$ is of
the type \eqref{CC2}, i.e.~corresponding to unilateral frictionless
Signorini contact on $\GC$, the $\eps$-Yosida regularization of
$\partial \ind_{\cone}$ is given  by $\yosd (u_\eps) = -\frac1{\eps}
[\hspace*{-.15em}[u_\eps]\hspace*{-.15em}]_{\mathrm n}^-$
with $[\hspace*{-.15em}[u_\eps]\hspace*{-.15em}]_{\mathrm n}
=\JUMP{u_\eps}{} {\cdot} \norm$ and with
$(\cdot)^-=-\min\{0,\cdot\}$, and the
second of \eqref{class-form-d-reg} reduces to
\[
 \kappa z \JUMP{u_\eps}{}-
\frac1{\eps}
\big[\hspace*{-.3em}\big[u_\eps\big]\hspace*{-.3em}\big]_{\mathrm n}
^{-}\norm+T(u_\eps,\DT{u}_\eps,\theta_\eps) = 0 \qquad \text{on $\SC$.}
\]
It follows from the above relation that, for fixed $\eps>0$, in
the case and $z(t,x)>0$
 one can express $\JUMP{u}{}$
as a function of the traction stress $T(u,\DT{u},\theta)$.
This also holds for
$[\hspace*{-.15em}[u]\hspace*{-.15em}]_{\mathrm n}$ in the case
 for $z(t,x)=0$
and $T_{\mathrm n}(u,\DT{u},\theta)>0$, while the tangential stress
$T_{\mathrm t}(u,\DT{u},\theta)=0$ because there is no friction.
\COLLL{Let us suppose} that the heat-transfer coefficient
$\eta(\cdot,z)$ vanishes if there is no contact, i.e.~on the set
$\{\JUMP{u}{};\ [\hspace*{-.15em}[u]\hspace*{-.15em}]_{\mathrm
n}<0\}$: then, by continuity, $\eta(\cdot,z)$ vanishes also on
$\{\JUMP{u}{};\ [\hspace*{-.15em}[u]\hspace*{-.15em}]_{\mathrm
n}\le0\}$, and hence by substitution one can express the
heat-transfer coefficient as a function of the normal stress and of
$z$.  \COLLL{We point out that the mentioned condition on
$\eta(\cdot,z)$ is,  to some extent,  a natural assumption,} also
advocated in the engineering literature,
cf.~e.g.~\cite{SadStu10MTCC}. Such an approach does not seem
mathematically amenable for the multidimensional non-penalized
Signorini problem; for $d=1$ we refer to \cite{AKRS02ODTC}.
\end{remark}
\par\noindent
\textbf{Scheme of the proof.}
 The proof of Theorem~\ref{th:4.1} shall be developed in the next
subsections by pursuing the  following scenario. First, in
Section~\ref{s:new}
 we shall devise a semi-implicit time discretization
(with a further regularization in the momentum equation), and prove
existence of   solutions to the time-discrete problem. Next, in
we shall  derive refined a priori estimates, enabling us to perform
the limit passage as the discretization time-step $\tau$ goes to $0$
in Sect.~\ref{ss:4.4}. In this way,  we
shall conclude the  proof of Theorem~\ref{th:4.1}. 

\section{Semi-implicit time discretization}\label{s:new}
We {perform a semi-implicit time-discretization using an equidistant
partition of $[0,T]$, with time-step $\tau>0$ and nodes $t_\tau^k:=k\tau$,
$k=0,\ldots,K_\tau$.

Hereafter, given any sequence $\{\phi^j\}_{j\ge 1}$, we shall use the following
notation for the {\it backward difference ope\-ra\-tor} and its iteration by,
respectively,
\begin{align}\label{not-interp0}
\dt\phi^k :=\frac{\phi^k{-}\phi^{k-1}\!\!\!}\tau, \qquad \dt^2\phi^k
:={\dt\big(\dt\phi^k\big)=}
\frac{\phi^k{-}2\phi^{k-1}{+}\phi^{k-2}\!\!\!}{\tau^2}.
\end{align}
 Secondly,  we recall the notion of piecewise constant  and
piecewise linear interpolants:
 for a given $K_\tau$-tuple $\{
b_\tau^k \}_{k=1}^{K_\tau} \subset \mathcal{B}$, $\left(\mathcal{B},
\| \cdot\|_{\mathcal{B}}\right)$ being some Banach space, the
left-continuous piecewise constant interpolant $\pwc b{\tau} : (0,T)
\to \mathcal{B}$, the right-continuous piecewise constant
interpolant $\upwc b{\tau} : (0,T) \to \mathcal{B}$,  and the
piecewise linear interpolant
 $\pwl b{\tau} : (0,T) \to \mathcal{B}$
 of the elements $\{ b_\tau^k
\}_{k=1}^{K_\tau}$ are the functions respectively  defined by
\begin{align}\label{not-interp}
\pwc b{\tau}(t) = b_\tau^k, \qquad \upwc b{\tau}(t) = b_\tau^{k-1},
\qquad \pwl b{\tau}(t) =\frac{t-t_\tau^{k-1}}{\tau} b_\tau^k +
\frac{t_\tau^k-t}{\tau}b_\tau^{k-1} \qquad \text{for $t \in
(t_\tau^{k-1}, t_\tau^k]$.}
\end{align}
Thirdly, we shall denote by $\overline{\mathsf{t}}_{\tau}$ and by
$\underline{\mathsf{t}}_{\tau}$ the left-continuous and
right-continuous piecewise constant interpolants associated with the
partition, i.e.
 $\bar{\mathsf{t}}_{\tau}(t) = t_\tau^k$ if $t_\tau^{k-1}<t \leq t_\tau^k $
and ${\underline{\mathsf{t}}}_{\tau}(t)= t_\tau^{k-1}$ if
$t_\tau^{k-1} \leq t < t_\tau^k $. For later use, we recall the
following elementary inequalities for all $t\in [0,T]$
\begin{align}
\label{elementary1} & \left\|\pwl b{\tau}(t)
\right\|_{\mathcal{B}}\leq \|\pwc b{\tau}(t)\|_{\mathcal{B}}+
\|\upwc b{\tau}(t)\|_{\mathcal{B}}=\|\pwc
b{\tau}(t)\|_{\mathcal{B}}+ \|\pwc
b{\tau}\left(t-{\underline{\mathsf{t}}}_{\tau}(t)\right)\|_{\mathcal{B}},
\\
 &
\label{elementary2}
  \left\|\pwc b{\tau}(t) - \pwl b{\tau}(t)
 \right\|_{\mathcal{B}}\leq
 \int_{{\underline{\mathsf{t}}}_{\tau}(t)}^{\bar{\mathsf{t}}_{\tau}(t)}\|\pwl
 {b'}{\tau}(s)\|_{\mathcal{B}}\, \d s\,.
\end{align}

\par\noindent
\textbf{Time-discrete problem.} We approximate the data $\FRM$,
$\fRM$ by local means, i.e. setting for all  $k=1,\ldots,K_{\tau}$
\[
\FRM_{\tau}^k:= \frac{1}{\tau}\int_{t_\tau^{k-1}}^{t_\tau^k}
\FRM(s)\, \d s\,,  \qquad  \fRM_{\tau}^k:=
\frac{1}{\tau}\int_{t_\tau^{k-1}}^{t_\tau^k} \fRM(s)\, \d s\,,
\]
and consider  the interpolants $\pwc \FRM{\tau}$, $\pwc \fRM{\tau}$,
and $\pwl \fRM{\tau}$  of the $K_\tau$-tuples $\{ \FRM_{\tau}^k
\}_{k=1}^{K_\tau}$, $\{\fRM_{\tau}^k \}_{k=1}^{K_\tau}$.
 In view of~\eqref{eFFe1}--\eqref{eFFe2},
the following estimates and
 strong convergences hold as $\tau \to 0$
\begin{subequations}
\label{data-converg}
\begin{align}
& \label{data-converg-1}
 \pwc \FRM{\tau} \to \FRM   \ \left\{
\begin{array}{ll}
\text{in $L^1(0,T;L^2(\Omega;\R^d))$}&\text{if $\varrho>0$,}
\\
\text{in $L^p(0,T;L^{6/5}(\Omega;\R^d))$ for all $1\le
p<\infty$}&\text{if $\varrho=0$;}
\end{array}
\right.
\\
& \label{data-converg-2}
\begin{aligned}
& \exists\, C>0 \ \  \forall\, \tau>0\, : \ \  \| \pwc
\fRM{\tau}\|_{L^\infty (0,T; L^{4/3}(\Gnew;\R^d))} \leq C \|
 \fRM\|_{L^\infty (0,T;L^{4/3}(\Gnew;\R^d))}\,,
 \\
& \pwc \fRM{\tau} \to \fRM \quad \text{in $ L^p
(0,T;L^{4/3}(\Gnew;\R^d)) $ for all $1\leq p <\infty$
 as $\tau \to
0$}\,,
\\
& \exists\, C>0 \ \  \forall\, \tau>0\, : \ \   \| \pwl
{\DT{\fRM}}{\tau}\|_{L^1 (0,T;L^{4/3}(\Gnew;\R^d))} \leq 2 \|
\DT{\fRM}\|_{L^1 (0,T;L^{4/3}(\Gnew;\R^d))}\,.
 \end{aligned}
\end{align}
\end{subequations}
Furthermore, we shall approximate  $\GRM$ and $\gRM$ with suitably
constructed discrete data $\{\GRM_{\tau}^k \}_{k=1}^{K_\tau}$,
$\{\gRM_{\tau}^k \}_{k=1}^{K_\tau}$ with
\begin{subequations}
\begin{equation}
\label{shall-use-later} \GRM_{\tau}^k \in W^{1,2}(\Omega)^*, \qquad
\gRM_{\tau}^k \in H^{1/2}(\partial \Omega)^*\qquad \text{ for all }
k=1,\ldots,K_\tau,
\end{equation}
 and  such that
\begin{align}
\label{data-converg-bis}
  \pwc \GRM{\tau} \to \GRM \quad \text{in  $L^1 (Q)$}, \qquad
 \pwc \gRM{\tau} \to \gRM \qquad \text{ in $L^1 (\Sigma)\ $ as $\ \tau\to0$}\,,
\end{align}
\end{subequations}
and approximate the initial datum $u_0$ with a
sequence $\{ u_{0,\tau}\} \subset  W_{\Gdir}^{1,\gamma}
(\Omega{\setminus} \GC; \R^d)$ (with
$\gamma>\max\{4,\frac{2\omega}{\omega-1}\}$ as assumed in
Problem~\ref{probk}) such that
\begin{equation}
\label{est-init-data}
\begin{gathered}
 \lim_{\tau \downarrow 0} \sqrt[\gamma]{\tau} \|
 e(u_{0,\tau})\|_{L^\gamma(\Omega;\R^d)}=0, \qquad  u_{0,\tau} \to
 u_0 \qquad\text{ in $W^{1,2}(\Omega;\R^d)\ $ as $\ \tau\to0$.}
 \end{gathered}
\end{equation}
We are now in the position of formulating the time-discrete problem,
which we again write in the classical formulation for notational
simplicity.

\begin{problem}
\label{probk} \upshape {Let
$\gamma>\max\{4,\frac{2\omega}{\omega-1}\}$.} Given
\begin{align}\label{IC2}
u_{\tau}^0=u_{0,\tau},\qquad u_{\tau}^{-1}=u_{0,\tau}-\tau\DT u_0,
\qquad z_{\tau}^0=z_{0},\qquad\w_{\tau}^0=\w_0,
\end{align}
find $\{(\uk,\wk,\zk) \}_{k=1}^{K_\tau}$ fulfilling
 for $k=1,...,K_\tau$ the equations in $ \Omega {\setminus} \GC $
\begin{subequations}\label{GMa}
\begin{align}
& \label{GM1a} \varrho \dt^2 \uk -{\rm div}\Big(\mathbb De\big(\dt
\uk\big) +\bbC e(\uk){-}\bbB\Theta(\wk)
+\tau\big|e(\uk)\big|^{\gamma-2}e(\uk)\Big)=\FRM_{\tau}^k,
\\&\label{GM3a}
\begin{aligned}
\dt \wk& -\mathrm{div}
\big(\mathcal{K}(\wk,e(\uk))\nabla\wk\big)
=\frac{2{-}\sqrt{\tau}}2\mathbb De\big(\dt
\uk\big){:}e\big(\dt \uk\big) +\Theta(\wk)\bbB{:}e\big(\dt
\uk\big) + \GRM_{\tau}^k,
\end{aligned}
\end{align}
\end{subequations}
 with the boundary conditions
\begin{subequations}\label{BC1}
\begin{align}
& \label{BC0} \uk = 0 &&\text{{on} $\Gdir$}\,,&&
\\& \label{BC1a}
\Big(\mathbb De\big(\dt \uk\big)
+\bbC e(\uk)-\Theta(\wk)\bbB
+\tau\big|e(\uk)\big|^{\gamma-2}e(\uk)\Big)\nu=\fRM_{\tau}^k\ \
&&\text{{on} $\Gnew$}\,,&&
\\
& \label{BC1b} \big(\mathcal{K}(\wk,e(\uk))\nabla\wk\big){\cdot}\nu
=\gRM_{\tau}^k &&\text{{on} $\partial \Omega$}\,,&&
\end{align}\end{subequations}
and the conditions on the contact boundary
\begin{subequations}\label{BC-cont}
\begin{align}
& \label{BC-conta}
 \left. \begin{array}{ll}
&\hspace{-1.6em}\partial\ind_{(-\infty,0]}\big(\dt \zk\big)
-a_0-a_1+\frac\dela{2}\big|\JUMP{\uk}{}\big|^2 +\tau^{\alpha} \zk +
r(\zk)\ni0 \vspace{1em}
\\
  &\hspace{-1.6em} r(\zk) \in \partial \ind_{[0,1]}(\zk)
  \end{array}
  \right\}
&&\text{{on} $\GC$}\,,&&
\\
&\label{BC-contb} \JUMP{\bbD e(\dt \uk)+
\bbC e(\uk){-}\Theta(\wk)\bbB
+\tau\big|e(\uk)\big|^{\gamma-2}e(\uk)}{}\norm=0 &&\text{{on}
$\GC$}\,,&&
\\\nonumber
 &
\dela{\zk}\JUMP{\uk}{} + \yosd(\JUMP{\uk}{})
 +\Big(\bbD e(\dt \uk)+\bbC e(\uk){-} \Theta(\wk)\bbB&&&&
\\\label{BC-contc}
&\qquad\qquad+\tau\big|e(\uk)\big|^{\gamma-2}e(\uk)\Big)\norm+
\tau^{\beta} \big(1{+}\big|\JUMP{\uk}{}\big|^2\big)^{\frac\mu2
-1}\JUMP{\uk}{} =0 &&\text{{on} $\GC$}\,,&&
\\\nonumber
&
\frac12\big(\mathcal{K}(\wk,e(\uk))\nabla\wk|_{\GC}^+ +
\mathcal{K}(\wk,e(\uk))\nabla\wk|_{\GC}^-\big){\cdot}\norm&&&&
\\\label{BC-contd}
&\qquad\qquad\qquad\qquad\qquad\qquad\qquad\quad
+\eta({\JUMP{u_{\eps\tau}^{k-1}}{},}\zk)\JUMP{\Theta(\wk)}{}=0
&&\text{{on} $\GC$}\,,&&
\\[.3em]
\label{BC-conte} & \JUMP{\mathcal{K}(\wk,e(\uk))\nabla\wk}{}\norm =
\minus a_1\dt{\zk} &&\text{{on} $\GC$}\,.&&
\end{align}\end{subequations}
\end{problem}

\begin{remark}[Semi-implicit discretization]\upshape
The value $u_{\eps\tau}^{k-1}$ at the level $k-1$ in
\eqref{BC-contd} makes the above scheme semi-implicit, not just
fully implicit as it would be in the case $u_{\eps\tau}^{k}$ were in
place of $u_{\eps\tau}^{k-1}$ in~\eqref{BC-contd}.
 This makes the proof of Lemma~\ref{lem-exist} easier.
\end{remark}

\begin{remark}[Regularization]\upshape
Like in \cite{tr1}, 
a regularizing term
$\tau|e(u)|^{\gamma -2} e(u)$
was added to the momentum equation in the bulk and to the
corresponding boundary/contact conditions, too.
Its role
is to compensate the
growth of the right-hand side terms in the momentum equation,
cf.~the proof of Lemma~\ref{lem-exist}.
 Moreover, with the aim of
obtaining  some suitable semiconvexity property of the approximate
stored energy (cf. Lemma~\ref{l:semiconvex} below), we have  also
introduced the monotone terms
\begin{equation}
\label{alpha-beta} \tau^{\alpha} z  \ \ \text{and}\ \  \tau^{\beta}
\big(1{+}\big|\JUMP{u}{}\big|^2\big)^{\frac\mu2 -1}\JUMP{u}{}, \ \
\text{with }  4<\mu<5 \ \text{and} \ \alpha,\beta \in (0,1)
\end{equation}
in the differential inclusion for
the delamination parameter and in the  boundary conditions for $u$
on $\GC$, respectively; see Lemma~\ref{l:semiconvex} for some further
specification of the exponents $\alpha$ and $\beta$.
\end{remark}

\begin{lemma}[Existence of weak solutions to Problem~\ref{probk}]
\label{lem-exist} Under the assumptions  of Theorem~\ref{th:4.1},
for every $k=1,...,K_{\tau}$ there exists a triple $ (\uk,\zk,
\wk)\in W_{\Gdir}^{1,\gamma}(\Omega{{\setminus}}\GC;\R^d) \times
L^\infty(\GC)\times W^{1,2}(\Omega{{\setminus}}\GC),$ fulfilling the
weak formulation of the boundary value problem
\eqref{GMa}--\eqref{BC-cont}. Moreover, $\wk\ge0 $ a.e. in
$\,\Omega$. If, in addition, \eqref{strict-pos} holds, then there
exists some constant $\chi^*>0$ (cf.~\eqref{pos} below)
such that, for sufficiently small $\tau$,
\begin{equation}
\label{posw} \wk\ge\chi^*>0 \quad \aein \ \Omega \ \ \text{for
every $k=1,...,K_{\tau}$.}
\end{equation}
\end{lemma}
\noindent{\it Proof.} The existence of a weak solution to the
boundary value problem \eqref{GMa}--\eqref{BC-cont}  can be proved
relying on the standard theory of pseudomonotone set-valued
operators (see e.g.~\cite[Chap.~2]{NPDE_roubicek}).  In particular,
we may apply Leray-Lions type theorems. Indeed, the strict
monotonicity of the main part of the operator which comes into play
in  the weak formulation of problem \eqref{GMa}--\eqref{BC-cont}
derives from the presence of the term $\tau |e(u)|^{\gamma-2}e(u)$.
The latter counteracts the quadratic nonlinearity in $e(u)$ of the
dissipative heat source in~\eqref{GM3a}.

We now show that this operator is coercive w.r.t.
the norm of $W_{\Gdir}^{1,\gamma}(\Omega{{\setminus}}\GC;\R^d)
\times L^\infty(\GC)\times W^{1,2}(\Omega{{\setminus}}\GC)$. To this
aim, first of all we test equation \eqref{GM1a} by $\uk$. Thus, with
elementary calculations, we find
\begin{align}
&\frac{\varrho}{2\tau^2} \|\uk \|_{L^2 (\Omega;\R^d)}^2 +
\frac{\mathsf{d}}{2\tau C_K^2} \|\uk \|_{W^{1,2} (\Omega;\R^d)}^2 +
\frac{\tau}{C_K^\gamma} \|\uk \|_{W^{1,\gamma} (\Omega;\R^d)}^\gamma
\nonumber
\\
& +\dela\int_{\GC} \zk\left|\JUMP{\uk}{} \right|^2\, \d S
 + \int_{\GC}\yosd (\JUMP{\uk}{}){\cdot} \JUMP{\uk}{}\, \d S
   +  \tau^{\beta}
\int_{\GC}\big(1{+}\big|\JUMP{\uk}{}\big|^2\big)^{\mu/2-1}\left|\JUMP{\uk}{}\right|^2\,
\d S \nonumber
\\
& \leq C \left( \varrho \|u_{\tau}^{k-1} \|_{W^{1,2}
(\Omega;\R^d)}^2 + \varrho \|u_{\tau}^{k-2} \|_{L^2 (\Omega;\R^d)}^2
+ \|\FRM_{\tau}^k\|_{W^{1,2} (\Omega;\R^d)^*}^2 +
\|\fRM_{\tau}^k\|_{H^{1/2}(\Gnew; \R^d)^*}^2 \right) \nonumber
\\
& \label{coerc1}
 + 4 \tau
|\bbC|^2 \| \Theta(u_{\tau}^k)\|_{L^2 (\Omega)}^2\,,
\end{align}
where \COL{
$\mathsf{d}:=\inf_{\xi \in \R^d, \, |\xi|=1} \bbD \xi {:} \xi>0$,
cf.~\eqref{posit}, and where}
we have also used Korn's inequality in the form
\begin{equation}
\label{e:korn}
\begin{aligned}
\exists\, C_K = C_K(\Omega)\, \ \ \forall\, \testu \in
W_{\Gdir}^{1,2}(\Omega;\R^d)\,: \qquad \| \testu
\|_{W^{1,2}(\Omega;\R^d)} \leq C_K \| e(\testu)
\|_{L^2(\Omega;\R^{d\times d})}.
\end{aligned}
\end{equation}
Also taking into account the monotonicity of the operator $\yosd$, we have
that the fourth, the fifth and the sixth term on the left-hand side
of~\eqref{coerc1} are non-negative. Secondly, we test \eqref{BC-conta}
by $\zk$. Supposing that $\zk \neq z_{\tau}^{k-1}$, so that $\partial \zeta_1
((\zk-z_{\tau}^{k-1})/\tau)=-a_1$, we obtain with trivial calculations
\begin{align}
\label{coerc2}
\frac\dela2\int_{\GC}\!\!\zk\left|\JUMP{\uk}{}\right|^2\,\d S
+ \tau^\alpha\|\zk\|_{L^2 (\GC)}^2 + \int_{\GC}\!\!r(\zk)\zk\,\d S
\le(2a_1{+}a_0)\|\zk\|_{L^1 (\GC)}\,.
\end{align}
Note that the third term on the left-hand side of~\eqref{coerc2} is
non-negative by monotonicity of the operator $\partial \ind_{[0,1]}.
$ Finally, we test \eqref{GM3a} by $\wk$, thus obtaining
\begin{align}
 \frac1{2\tau}\| \wk\|_{L^2 (\Omega)}^2  & +  \mathsf{K}
\int_{\Omega}|\nabla \wk|^2 \, \d x+ \int_{\GC}
\eta({\JUMP{u_{\eps\tau}^{k-1}}{},}\zk) \JUMP{\Theta (\wk)}{}
\JUMP{\wk}{}\, \d S \nonumber\\ & \leq \frac1{2\tau}\|
\w_{\eps\tau}^{k-1}\|_{L^2 (\Omega)}^2 + I_{1} + I_{2} + I_{3} +
I_{4} \label{coerc3}
\end{align}
where we have used that $\mathcal{K}$ is positive definite (cf.
with~\eqref{30dprimo}). As for the remaining terms $I_{j}$,
$j=1,\ldots,4$, we have
\begin{align}
I_1 &  = \frac{2-\sqrt{\tau}}{2\tau^2} \int_{\Omega{\setminus} \GC}
\bbD (e(\uk) - e(u_{\eps\tau}^{k-1}))
 {:}  (e(\uk) - e(u_{\eps\tau}^{k-1}))  \wk\, \d x
 \nonumber
\\ &
  \leq \frac1{4\tau}\| \wk\|_{L^2
(\Omega)}^2 + C \left(\| e(\uk) \|_{L^4 (\Omega; \R^{d \times d})}^4
+\| e(u_{\eps\tau}^{k-1}) \|_{L^4 (\Omega; \R^{d \times d})}^4
\right)
\nonumber \\
& \leq \frac1{4\tau}\| \wk\|_{L^2 (\Omega)}^2 +
\frac{\tau}{4C_K^\gamma} \|\uk \|_{W^{1,\gamma}(\Omega;\R^d)}^\gamma
+ C \|u_{\eps\tau}^{k-1}\|_{W^{1,\gamma}(\Omega;\R^d)}^\gamma,
 \label{coerc31}
\end{align}
where we have used  H\"older inequality and the fact that
$\gamma>4$. Furthermore, relying on
\eqref{growthTheta}, and setting $p_{\omega}= 2 \omega/(\omega-1)$,
we find
\begin{align}
 I_2  & =\frac1{\tau}\int_{\Omega{\setminus} \GC}
  \Theta (\wk)\, \bbB{:}\left( e(\uk)
-e(u_{\eps\tau}^{k-1})\right) \wk\, \d x  \nonumber
 \\ & \leq \frac1{16\tau}\|
\wk\|_{L^2 (\Omega)}^2 + C \int_{\Omega}
\big|e(\uk)-e(u_{\eps\tau}^{k-1})\big|^2
\,\big|(\wk)^{2/\omega}\!+1\big|\, \d x \nonumber
\\
& \leq \frac1{8\tau}\| \wk\|_{L^2 (\Omega)}^2 + C\left( \| e(\uk)
\|_{L^{p_\omega} (\Omega;\R^{d \times d})}^{p_\omega} + \|
e(u_{\eps\tau}^{k-1}) \|_{L^{p_\omega} (\Omega;\R^{d \times
d})}^{p_\omega}+1 \right)
 \nonumber \\
& \leq \frac1{8\tau}\| \wk\|_{L^2 (\Omega)}^2 +
\frac{\tau}{8C_K^\gamma} \|\uk \|_{W^{1,\gamma}
(\Omega;\R^d)}^\gamma  + C\left(
\|u_{\eps\tau}^{k-1}\|_{W^{1,\gamma} (\Omega;\R^d)}^\gamma + 1\right),
\label{coerc32}
\end{align}
where we have successively used
H\"older's
and Young inequalities, and
that
$\gamma>p_{\omega}$ due to our assumption that
$\gamma>\max\{4,\frac{2\omega}{\omega-1}\}$.
 Besides, using that  $0 \leq \zk \leq 1 $ and
$ 0 \leq z_{\eps\tau}^{k-1} \leq 1$ a.e. in $\GC$, and the
continuous embedding~\eqref{e:contemb}, we find
\[
I_3 = \minus a_1\int_{\GC}\!\!\dt{\zk} \frac{ \wk|_{\GC}^+{+}
\wk|_{\GC}^-}2\, \d S \leq C \|\dt{\zk}\|_{L^2 (\GC)} \| \wk \|_{L^2
(\GC)} \leq \rho_1\|\wk\|_{W^{1,2}(\Omega)}^2\!+C_{\rho_1}
\]
where we choose $\rho_1>0$ ($C_{\rho_1}$ being some constant
depending on $\rho_1>0$) in such a way as to absorb $\| \wk
\|_{W^{1,2}(\Omega)}^2$ into the left-hand side of~\eqref{coerc3}.
Finally, we have
 \begin{align} I_4 & =
\int_\Omega\GRM_{\tau}^k\wk\,\d x+
\int_\Gamma\gRM_{\tau}^k\wk\,\d S
\leq \rho_2 \|
\wk \|_{W^{1,2}(\Omega)}^2 + C_{\rho_2} \left(\| \GRM_{\tau}^k
\|_{W^{1,2}(\Omega)^*}^2 + \| \gRM_{\tau}^k \|_{H^{1/2}(\partial
\Omega)^*}^2 \right)\,,\label{coerc333}
\end{align}
in which  we again choose a suitably small $\rho_2$.
Collecting \eqref{coerc1}--\eqref{coerc333}, we readily conclude an
estimate for $\| \uk \|_{W^{1,\gamma} (\Omega;\R^d)}$, $ \|
\zk\|_{L^\infty (\GC)}$ and $ \|\wk \|_{W^{1,2}(\Omega)} $.

Since  $\w_{{\eps}\tau}^k\in W^{1,2}(\Omega{{\setminus}}\GC)$, we
have that ${-}[\w_{{\eps}\tau}^k]^-\in
W^{1,2}(\Omega{{\setminus}}\GC)$ is a legal test function for
(\ref{GM3a}). Hence, we use recursively that
$\w_{{\eps}\tau}^{k-1}\ge0$ a.e. in $\Omega$ (starting from the
initial condition $\w_{{\eps}\tau}^0=\w_0\ge0$ a.e. in $\Omega$,
cf.~\eqref{wzero}), the fact that $\GRM_{\tau}^k \geq 0$ a.e. in
$\Omega $ and $\gRM_{\tau}^k\ge0$ a.e. in $\partial \Omega$ (cf.
with~\eqref{posg1}--\eqref{posg2}), and that $-a_1 \dt{\zk} \geq 0$
a.e. in $\GC$, and finally the property
$[\w_{{\eps}\tau}^k]^-\Theta(\w_{{\eps}\tau}^k)=0$  a.e. in $\Omega$,
 due to the fact that $\Theta $ in non-decreasing (cf.~(\ref{K-T})). Thus,
we conclude that $[\w_{{\eps}\tau}^k]^- =0 $ a.e.~in $\Omega$, whence
$\wk\ge0$ a.e.~in $\Omega$.

Finally, we prove~\eqref{posw} by adapting to the time-discrete
setting a comparison argument from~\cite[Sect.~4.2.1]{fpr09}.
Exploiting the fact that~$\GRM_{{\eps}\tau}^k\ge0$ a.e.\ in
$\Omega $,
we deduce from~\eqref{GM3a} that
\begin{align}
\dt \wk -\mathrm{div} \big(\mathcal{K}(\wk,e(\uk))\nabla\wk\big) &
\geq
 \frac{\mathsf{d}}4 \big|e(\dt \uk)\big|^2
- C |\Theta(\wk)|^2
\geq - C' |\wk|^2 \quad \text{in $\Omega{\setminus}\GC$}
\label{feiresl}
\end{align}
for any $k=1,\ldots,K_\tau$ and some $C'>0$ independent of $\tau$
and $\eps$, where \COL{$\mathsf{d}>0$ is the positive-definiteness
constant of $\bbD$, cf.~\eqref{posit}, and where} the last
inequality ensues from~\eqref{K-T} and~\eqref{30b}. We
compare~\eqref{feiresl} with the finite difference equation
\begin{equation}
\label{finite-difference} \dt \chi_{k} = -C'|\chi_k|^2\qquad\forall\,
k=1, \ldots, K_\tau,
\end{equation}
with $C'$ being the same constant as in~\eqref{feiresl}. \COL{In fact,
this is an implicit discretization of the ordinary-differential
equations $\DT{\chi}+C'|\chi|^2=0$ which, for
$\chi(0)=h_0(\theta^*)>0$ with $\theta^*$ from \eqref{strict-pos},
gives a sub-solution of the (continuous) heat equation. This initial-value
problem has the solution $\chi(t)=1/(C't+1/h_0(\theta^*))$ so that, in
particular, $\chi(\cdot)\ge1/(C'T+1/h_0(\theta^*))>0$ on $[0,T]$.
Now} we solve~\eqref{finite-difference} recursively starting from
the initial datum $\chi_0=h_0(\theta^*)>0$. In this way we
obtain an approximate solution to the mentioned initial-value
problem which, for $\tau\to0$, converges uniformly on the considered
finite interval $[0,T]$. In particular, for $\tau>0$ sufficiently
small, we may take for granted that, say,
\begin{equation}\label{pos}
\chi_k\ge\chi^* :=\frac1{C'T+1/h_0(\theta^*)+1}>0\qquad\forall\,k=1,
\ldots, K_\tau.
\end{equation}
For every $k=1,\ldots, K_\tau$ we subtract~\eqref{finite-difference}
from \eqref{feiresl} (the latter supplemented with the boundary
conditions~\eqref{BC1b}, \eqref{BC-contd}--\eqref{BC-conte}), and we
test the resulting inequality by $-(\wk-\chi_k)^-$. Thus, for all
$k=1, \ldots, K_\tau,$
\begin{equation}\label{almost-final}
\frac12\dt\big((\wk{-}\chi_k)^-\big)^2\le
-(\wk{-}\chi_k)^-\dt(\wk{-}\chi_k) \le C'(|\wk|^2-|\chi_k|^2
)(\wk-\chi_k)^{-}\le0
\end{equation}
in $\Omega{\setminus}\GC$;
the latter inequality also due to the previously proved positivity
$\wk \geq 0$ a.e. in $\Omega$. Summing~\eqref{almost-final} over
$k=1,\ldots,K_\tau$,
we easily conclude that $(\wk(x)-\chi_k)^-=0$ for almost
all $x \in \Omega$ and  for every $ k=0, \ldots, K_\tau$, whence
$\wk\ge\chi_k\ge\chi^*>0$ a.e. in $\Omega$. This concludes the proof
of~\eqref{posw}.
 $\hfill\Box$

\medskip
\par\noindent
\textbf{Approximate solutions.} In accordance with
notation~\eqref{not-interp}, for all $\tau>0$ we shall denote by
\begin{itemize}
\item   $\pwc u {{\eps}\tau}$, $\upwc u {{\eps}\tau}$,
 $\pwc {\w}{{\eps}\tau}$,
$\upwc {\w}{{\eps}\tau}$,
 and $\pwc
{z}{{\eps}\tau}$,  
the
piecewise constant interpolants of the elements
$\{\uk\}_{k=1}^{K_\tau} $, $\{\wk\}_{k=1}^{K_\tau} $, and
$\{\zk\}_{k=1}^{K_\tau} $;
\item by $\ut$, $\wt$, and $\zt$, the related piecewise linear
interpolants.
\end{itemize}

We shall now state the weak formulations of (the boundary value
problems for) equations~\eqref{GM1a}, \eqref{GM3a}, in terms of the
interpolants so far introduced by using   ``discrete test
functions''. Indeed, one verifies that  for every $K_\tau$-tuples
$\{\testu_\tau^k\}_{k=1}^{K_\tau} \subset
W_{\Gdir}^{1,2}(\Omega{\setminus} \GC; \R^d)$ and
$\{\testw_\tau^k\}_{k=1}^{K_\tau} \subset
W^{1,2}(\Omega{\setminus}\GC)$, the approximate solutions $(\pwc
u{{\eps}\tau},\pwc\w{{\eps}\tau}, \pwc z{{\eps}\tau},\pwl
u{{\eps}\tau},\pwl \w{{\eps}\tau}, \pwl z{{\eps}\tau})$ fulfil the
following:
\\
\textbf{the discrete (weak) momentum balance equation}
\begin{align}
 \nonumber \int_{Q}  & \big(\bbD e(
{\pwl{\DT{u}}{{\eps}\tau}})
  +\bbC e(\pwc u{{\eps}\tau})-\bbB\Theta(\pwc\w{{\eps}\tau})
+\tau |e(\pwc u {{\eps}\tau})|^{\gamma-2}e(\pwc u {{\eps}\tau})
\big){:} e(\pwc \testu{\tau})\, \d x\d t
\\ & \nonumber
+ \int_{\SC} \left( \dela{\pwc z {{\eps}\tau}}\JUMP{\pwc u
{{\eps}\tau}}{} +\yosd(\JUMP{\pwc u {{\eps}\tau}}{}) +\tau^\beta
\big(1{+}\big|\JUMP{\pwc u
{{\eps}\tau}}{}\big|^2\big)^{\mu/2-1}\JUMP{\pwc u {{\eps}\tau}}{}
\right){\cdot} \JUMP{\pwc\testu\tau}{} \, \d x\d t
\\ & \nonumber
-\int_{{\eps}\tau}^T \int_{\Omega}
\varrho{\pwl{\DT{u}}{{\eps}\tau}}(\cdot-\tau) {\cdot}
{\pwl{\DT{\testu}}{\tau}} \, \d x\d t
 + \varrho \int_{\Omega}
{\pwl{\DT{u}}{{\eps}\tau}}(T){\cdot}{\pwl{\testu}{\tau}} (T)\, \d x
\\& \label{e:discrmom}= \rho \int_{\Omega} \DT{u}_{0,\tau}{\cdot}
\pwl{\testu}{\tau}(\tau)\, \d x + \int_{Q}  \pwc \FRM{\tau} {\cdot}
\pwc{\testu}{\tau} \, \d x\d t + \int_{\Snew} \pwc\fRM{\tau} {\cdot}
 \pwc{\testu}\tau \,\d S  \d t
\end{align}
(where we have used the  notation $ \DT{u}_{0,\tau} =
\frac{u^0_{{\eps}\tau} - u^{-1}_{{\eps}\tau}}{\tau} = \DT{u}_0$),
 which can be obtained from \eqref{GM1a}, \eqref{BC0}, \eqref{BC1a}, \eqref{BC-contb},
 and \eqref{BC-contc} by
using a suitable discrete ``by-part'' summation formula
(cf.~\cite[Formula (4.49)]{tr1});
\\
\textbf{the discrete (weak) enthalpy equation}
\begin{align}
 \nonumber
&\int_\Omega \pwl\w{{\eps}\tau}(T)\pwl\testw{\tau}(T)\, \d x+
\int_{Q} \mathcal{K}(e(\pwc
u{{\eps}\tau}),\pwc\w{{\eps}\tau})\nabla\pwc\w{{\eps}\tau}{\cdot}\nabla
\pwc\testw{\tau} \, \d x\d t \\ \nonumber &
+\int_{\SC}\het({\JUMP{u_{\eps\tau}^{k-1}}{},} \pwc z{{\eps}\tau})
\JUMP{\Theta(\pwc\w{{\eps}\tau})}{}\JUMP{\pwc\testw{\tau}}{}\,\d S\d
t -\int_\tau^T\upwc\w{{\eps}\tau}{\pwl{\DT{\testw}}{\tau}}\, \d x\d
t
\\ & \nonumber  =\int_Q\!\left(\frac{2-\sqrt{\tau}}{2}\bbD
e(\pwl{\DT{u}}{{\eps}\tau}) {:}
e(\pwl{\DT{u}}{{\eps}\tau})+\Theta(\pwc\w{{\eps}\tau})\bbB {:}
e(\pwl{\DT{u}}{{\eps}\tau})\right)\pwc \testw{\tau}\,\d x\d t
+\int_{Q}\pwc\GRM{\tau}\pwc\testw{\tau}\,\d x \d t
\\ &
\minus \int_{\SC}\!\!a_1\pwl{\DT{z}}{{\eps}\tau}\frac{ \pwc
\testw{\tau}|_{\GC}^+{+} \pwc \testw{\tau}|_{\GC}^-}2 \, \d S\d t
+\int_\Omega\!\w_0  \pwc \testw{\tau}(\tau)\,\d x
\label{weak-heat-discr}
+\int_{\Sigma}\pwc\gRM{\tau} \pwc \testw{\tau}\,\d S \d t\,,
\end{align}
 again obtained from~\eqref{GM3a}, \eqref{BC1b},
\eqref{BC-contc}, \eqref{BC-contd}, and \eqref{BC-conte} by the use
of the  summation  (cf.~\cite[Formula (4.51)]{tr1});
\\
\textbf{the discrete flow rule of the delamination parameter} (cf.
\eqref{BC-conta})
\begin{align}
\label{discrez} \left. \begin{array}{ll}
&\hspace{-3.5em}\partial\ind_{(-\infty,0]}\big(\pwl
{\DT{z}}{{\eps}\tau}(t)\big) +\frac\dela{2}\big|\JUMP{\pwc
u{{\eps}\tau}(t)}{}\big|^2 +\tau^{\alpha} \pwc z{{\eps}\tau(t)} +
r(\pwc z{{\eps}\tau}(t))-a_0-a_1 \ni0 \vspace{1em}
\\
  &\hspace{-3.5em} r(\pwc z{{\eps}\tau}(t)) \in \partial \ind_{[0,1]}(\pwc z{{\eps}\tau}(t))
  \end{array}
  \right\}
 \quad \text{on $\SC$}\,.
\end{align}
\par\noindent
\textbf{A-priori estimates.}
 Like in~\cite[Lemma~4.1]{tr1}, we
derive some further energetic information on the approximate
solutions (see Lemma~\ref{lem-first} later on) by recurring to an
auxiliary minimization problem. With this aim, we first proceed to
the validation of a suitable (strict) semiconvexity property of the
stored energy.

We further introduce the short-hand notation for the regularized stored energy
\begin{subequations}\label{Phi-eps-tau}\begin{align}
&
\begin{aligned}
\Phi_{{\eps}\tau}(u,z):  =  & \int_{\Omega{\setminus}\GC}\!\Big(
\frac12\bbC e(u){:} e(u)  +\frac\tau\gamma|e(u)|^\gamma\Big)\,\d
x +\yosapp(\JUMP{u}{})\\ &  + \int_{\GC}\!\!\Big(\frac\dela 2
z\big|\JUMP{u}{}\big|^2 +\frac{\tau^{\alpha}}2|z|^2
+\frac{\tau^\beta}{\mu}\big(1{+}\big|\JUMP{u}{}\big|^2\big)^{\mu/2}+
\ind_{[0,1]}(z)-a_0z\Big)\,\d S.
\end{aligned}
\end{align}\end{subequations}
\begin{lemma}\label{l:semiconvex}
Under the assumptions of Theorem~\ref{th:4.1},
 suppose further  that the exponents $\mu \in (4,5)$,  $\alpha, \,
\beta  \in (0,1)$ in~\eqref{alpha-beta} comply with
\begin{equation}
\label{e:exponents}
 \alpha(\mu-2) + 2\beta < \frac{\mu-4}2.
\end{equation}
Then, for every $\delam>0$ there exists $\tau_\delam >0$ such that
for all $0 <\tau <\tau_\delam$ the function on $ W_{\Gdir}^{1,2}
(\Omega{\setminus} \GC; \R^d) \times L^\infty (\GC)$ given by
\begin{align}
 (u,z)
 \mapsto\Phi_{{\eps}\tau}(u,z) +\int_\Omega\!\frac{\bbD
e(u){:} e(u)}{2\sqrt\tau}\,\d x \  \ \ \text{is strictly convex.}
\label{e:impo}
\end{align}
\end{lemma}
\noindent{\it Proof.} 
%
Following the  calculations in \cite{book}, we prove \eqref{e:impo}
by investigating the monotonicity of the multivalued mapping
\[
W^{1,2}(\Omega{\setminus} \GC;\R^d)\times L^2(\GC) \rightrightarrows
W^{1,2}(\Omega{\setminus} \GC;\R^d)^*\times L^2(\GC)\, : \ \
(u,z)\mapsto \pl\Phi_{{\eps}\tau}(u,z)\,.
\]
To this goal, we have to estimate from below
\begin{align}\nonumber
 & \big\langle\pl\Phi_{{\eps}\tau}(u_1,z_1)
-\pl\Phi_{{\eps}\tau}(u_2,z_2),(u_1-u_2,z_1-z_2)\big\rangle
\\\nonumber &
=\int_{\Omega{{\setminus}}\GC}\!\!\!\big(\bbC e(u_1{-}u_2)
+\tau|e(u_1)|^{\gamma-2}e(u_1)-\tau|e(u_2)|^{\gamma-2}e(u_2)\big)
{:} e(u_1{-}u_2)\,\d x
\\\label{semi-convex2}
& +\int_{\GC}\!\mathcal{L}\big(\JUMP{u_1}{},\JUMP{u_2}{},
z_1,z_2\big)\,\d S\,.
\end{align}
As for the latter term, using the short-hand notation $s_i=
\JUMP{u_i}{}$  for $i=1,2$, and $r(z_i) \in \partial
\ind_{[0,1]}(z_i)$  as in~\eqref{BC-conta}, we can estimate the last term
in \eqref{semi-convex2} as
\begin{align}\nonumber
&\mathcal{L}(s_1,s_2,z_1,z_2)= \dela
(z_1s_1{-}z_2s_2){\cdot}(s_1{-}s_2) +\left(\yosd(s_1) {-} \yosd(s_2)
\right){\cdot}(s_1{-}s_2)
+\frac\dela{2}(z_1{-}z_2)(|s_1|^2{-}|s_2|^2)
\\ &\qquad +
(r(z_1){-}r(z_2))(z_1{-}z_2)
 + \tau^{\alpha}|z_1{-}z_2|^2
+\tau^{\beta}\big((1{+}|s_1|^2)^{\mu/2-1}s_1{-}(1{+}|s_2|^2)^{{\mu}/2-1}s_2\big){\cdot}(s_1{-}s_2)
\label{e:used-later}
\\&\quad
\geq \dela z_1|s_1{-}s_2|^2 +\frac\dela
{2}(z_1{-}z_2)(s_1{+}3s_2){\cdot}(s_1{-}s_2)
+\tau^{\alpha}(z_1{-}z_2)^2 \nonumber
\\ &\qquad
+\tau^{\beta}\big((1{+}|s_1|^2)^{{\mu}/2-1}s_1{-}(1{+}|s_2|^2)^{{\mu}/2-1}s_2\big){\cdot}(s_1{-}s_2)
\nonumber
\\ &\quad \geq \frac{\tau^{\alpha}}2(z_1{-}z_2)^2 - \frac{\dela^2}{8 \tau^\alpha}
|s_1{-}s_2|^2|s_1{+}3s_2|^2 \nonumber
+\tau^{\beta}\big((1{+}|s_1|^2)^{{\mu}/2-1}s_1{-}(1{+}|s_2|^2)^{{\mu}/2-1}s_2\big){\cdot}(s_1{-}s_2)
\nonumber
\\ &\quad
\geq \frac{\tau^{\alpha}}2\big(z_1{-}z_2\big)^2
-S_{\delam,\tau}\big|s_1{-}s_2\big|^2
 \label{eq6:delam-k-monotone1}
\end{align}
for some positive constant $-S_{\delam,\tau}$. Indeed, the first inequality follows
from the positivity (by monotonicity) of the second and fourth term on the
right-hand side of~\eqref{e:used-later}, and from simple algebraic
manipulations. So does the second inequality. To conclude the final
inequality \eqref{eq6:delam-k-monotone1} for some constant $S_{\delam,\tau}>0$
depending on $\delam$ and $\tau$, we have used that
(cf.~\cite[Lemma 5.2]{ThoMie09DNEM})
\[
\exists\, C_{\mu}>0\, :  \ \
(1{+}|s_1|^2)^{{\mu}/2-1}s_1{-}(1{+}|s_2|^2)^{{\mu}/2-1}s_2){\cdot}(s_1{-}s_2)\ge
C_{\mu}(|s_1|^{{\mu}-2}+|s_2|^{{\mu}-2})|s_1{-}s_2|^2
\]
 for all $ s_1,s_2 \in \R^d $
and that (since ${\mu}>4$)
\begin{equation}
\label{e:key} \forall\,\delam,\,  \tau>0 \ \
\exists\,S_{\delam,\tau}>0 \ \ \forall\, s_1,s_2 \in \R^d\,: \ \
c_{\mu} \tau^{\beta}(|s_1|^{{\mu}-2}+|s_2|^{{\mu}-2}) {-}
\frac{\delam^2}{8 \tau^\alpha} |s_1{+}3s_2|^2 \geq -
S_{\delam,\tau}\,.
\end{equation}
Combining~\eqref{semi-convex2} with~\eqref{eq6:delam-k-monotone1}
 the  boundedness of the jump
operator $u\mapsto\JUMP{u}{}$ from $
W_{\Gdir}^{1,2}(\Omega{\setminus} \GC;\R^d)$ onto $ L^2(\GC;\R^d)$,
as well as using Korn's inequality, we conclude that
\[
\begin{aligned}
\big\langle\pl\Phi_{{\eps}\tau}(u_1,z_1)
-\pl\Phi_{{\eps}\tau}(u_2,z_2),(u_1-u_2,z_1-z_2)\big\rangle
 \geq \frac{\tau^{\alpha}}2\|z_1{-}z_2\big\|_{L^2 (\GC)}^2
 {-} C
S_{\delam,\tau} \| e(u_1) {-} e(u_2)\|_{L^2 (\Omega;\R^{d\times
d})}^2
\end{aligned}
\]
with the constant $C$ depending on
\COL{the positive-definiteness constant of $\bbC$ (cf.~\eqref{posit}),}
on the norm of the trace operator from
$W^{1,2}(\Omega{{\setminus}}\GC;\R^d)$ to $ L^2(\GC;\R^d)$, and on
the constant in Korn's inequality~\eqref{e:korn}. Finally,
 the key observation is that, for $\delam>0$ fixed,
 the constant $S_{\delam,\tau}$
 in~\eqref{e:key} has the following qualitative behaviour
\[
S_{\delam,\tau} \sim
\frac{1}{\tau^{\alpha\frac{{\mu}-2}{{\mu}-4}+\frac{2\beta}{{\mu}-4}}}
\quad \text{as $\tau \to 0$.}
\]
 Thus,  using  condition~\eqref{e:exponents}, it can be verified
that for all $\delam>0$ there exists ${\tau}_\delam>0$ such that for
$0<\tau<\tau_\delam$ there holds $C S_{\delam,\tau} \leq
\mathsf{d}/{\sqrt{\tau}}$; \COL{again $\mathsf{d}>0$ is the
positive-definiteness constant of $\bbD$.}
This yields~\eqref{e:impo}. $\hfill\Box$ 

\begin{lemma}[First a priori information]\label{lem-first}
Under the assumptions of Theorem~\ref{th:4.1}, for all $\varrho \geq
0$ and for every $\eps>0$ {there is $\tau_\eps$ such that for all}
$0<\tau<\tau_\eps$ the approximate solutions $(\pwc u
{{\eps}\tau},\pwc\w{{\eps}\tau},\pwc z {{\eps}\tau},\pwl u
{{\eps}\tau},\pwl \w{{\eps}\tau},\pwl z {{\eps}\tau})$ fulfil the
following ``discrete mechanical energy'' inequality
\begin{align}\nonumber
&\!\!\!\!
T_\mathrm{kin}^{\varrho}\big(\pwl{\DT{u}}{{\eps}\tau}(t)\big)
+\Phi_{{\eps}\tau}\big(\pwc{u}{{\eps}\tau}(t),\pwc{z}{{\eps}\tau}(t)\big)
+\int_0^{\bar{\mathsf{t}}_{{\eps}\tau}(t)}\bigg(\int_{\Omega}\!\!
\frac{2{-}\sqrt\tau}2 \bbD
e\big(\pwl{\DT{u}}{{\eps}\tau}(s)\big){:}
e\big(\pwl{\DT{u}}{{\eps}\tau}(s)\big)\,\d x
 +\int_{\GC}\!\!\zeta_1\big(\pwl{\DT{z}}{{\eps}\tau}(s)\big)\,\d
S\bigg)\, \d s
\\
&\le  T_\mathrm{kin}^{\varrho}\big(\DT{u}_{0,\tau})
+\Phi_{{\eps}\tau}\big(u_{0,\tau},z_{0})\nonumber
\\
&+\int_0^{\bar{\mathsf{t}}_{{\eps}\tau}(t)}\bigg(\int_\Omega
\Theta(\pwc{\w}{{\eps}\tau}(s))\bbB{:}
e\big(\pwl{\DT{u}}{{\eps}\tau}(s)\big) \,\d x + \int_{\Omega} \pwc
\FRM{\tau}(s) {\cdot} \pwl{\DT{u}}{{\eps}\tau}(s)\,\d x +\int_{\Gnew}
\pwc \fRM{\tau}(s) {\cdot} \pwl{\DT{u}}{{\eps}\tau}(s)\,\d S \bigg)\,
\d s, \label{disc-energy0}
\end{align}
  as well as the
following ``discrete total energy'' inequality
\begin{align}\nonumber
& T_\mathrm{kin}^{\varrho}\big(\pwl{\DT{u}}{{\eps}\tau}(t)\big)
+\Phi_{{\eps}\tau}\big(\pwc{u}{{\eps}\tau}(t),\pwc{z}{{\eps}\tau}(t)\big)+\int_\Omega\pwc{\w}{{\eps}\tau}(t)\,\d
x
\le T_\mathrm{kin}^{\varrho}\big(\DT{u}_{0,\tau})
+\Phi_{{\eps}\tau}\big(u_{0,\tau},z_{0}) +\int_\Omega\w_0\,\d
x\nonumber
\\
& \hspace{1em} +\int_0^{\bar{\mathsf{t}}_{{\eps}\tau}(t)}
\bigg(\int_{\Omega} \pwc \FRM{\tau}(s) {\cdot}
\pwl{\DT{u}}{{\eps}\tau}(s)\,\d x +\int_{\Gnew}\!\!\pwc\fRM{\tau}(s)
{\cdot}\pwl{\DT{u}}{{\eps}\tau}(s)\,\d S + \int_{\Omega} \pwc
\GRM{\tau}(s) \, \d x +\int_{\partial \Omega} \pwc \gRM{\tau}(s)\,\d
S \bigg)\, \d s
 \label{disc-energy}
\end{align}
and also the ``discrete semistability'' for a.a. $t \in (0,T)$
(where $\Phi_{{\eps}\tau}$ is from~\eqref{Phi-eps-tau})
\begin{align}\label{disc-semistability}
&\Phi_{{\eps}\tau}\big(\pwc{u}{{\eps}\tau}(t),\pwc{z}{{\eps}\tau}(t)\big)
\le\Phi_{{\eps}\tau}\big(\pwc{u}{{\eps}\tau}(t),\tilde
z\big)+\calD\left(\tilde{z}- \pwc{z}{{\eps}\tau}(t) \right) \quad
\text{for all $\tilde z\in L^\infty(\GC)$.}
\end{align}
\end{lemma}

\noindent{\it Proof.}  Let us now fix a solution
$(u_{{\eps}\tau}^k,z_{{\eps}\tau}^k,\w_{{\eps}\tau}^k)$ of
Problem~\ref{probk}. Recall that such a triple exists  thanks
to~Lemma~\ref{lem-exist}. Let us consider an auxiliary minimization
problem, namely
\begin{align}\label{GMd-min1}
\hspace*{-.6em}\left.\begin{array}{ll}
\mbox{minimize}&\displaystyle{\int_\Omega\varrho
\dt^2u_{{\eps}\tau}^k {\cdot} u +\,\big(1{-}\sqrt\tau\,\big) \bbD
e\big(
\dt\uk
\big){:} e(u)
}
\\[0mm]&\quad
\displaystyle{ +\frac{\tau^{3/2}}2\bbD
e\big(\frac{u{-}u_{{\eps}\tau}^{k-1}}{\tau}\big) {:}
e\big(\frac{u{-}u_{{\eps}\tau}^{k-1}}{\tau}\big)
+\,
\Theta(\w_{{\eps}\tau}^k)\bbB{:} e(u)\,\d x}
\\[2mm]&\quad
\displaystyle{+\,\tau\int_{\GC}\zeta_1\Big(\frac{z{-}z_{{\eps}\tau}^{k-1}}{\tau}\Big)\,\d S
+\Phi_{{\eps}\tau}(u,z)- \int_{\Omega}\FRM_{\tau}^k {\cdot} u \, \d x-
\int_{\Gnew}\fRM_{\tau}^k }{\cdot}{u}\,\d S
\\[2mm]\mbox{subject to}&
(u,z)\in W_{\Gdir}^{1,\gamma}(\Omega{\setminus}\GC;\R^d)\times L^\infty(\GC)\,.
\end{array}\right\}\hspace*{-.4em}
\end{align}
By convexity  of $\Phi_{{\eps}\tau}$ and coercivity (cf. the
calculations developed in the proof of Lemma~\ref{lem-exist}), it is
immediate to check that the minimization problem (\ref{GMd-min1})
has a solution which we denote by $(\tilde u_{{\eps}\tau}^k,\tilde
z_{{\eps}\tau}^k)$; of course, it depends on the pair $(\uk,\wk)$ in
general. Writing optimality conditions for $(\tilde
u_{{\eps}\tau}^k,\tilde z_{{\eps}\tau}^k)$ gives, for all $\testu\in
W_{\Gdir}^{1,\gamma}(\Omega{\setminus} \GC;\R^d)$, that
\begin{subequations}\label{GMd-mod}
\begin{eqnarray}
\nonumber &&\int_\Omega \dt^2u_{{\eps}\tau}^k{\cdot} \testu
+\Big(\sqrt\tau\mathbb D e\big(\frac{\tilde
u_{{\eps}\tau}^k{-}u_{{\eps}\tau}^{k-1}}{\tau}\big) +\bbC e(\tilde
u_{{\eps}\tau}^k)+\tau\big|e(\tilde
u_{{\eps}\tau}^k)\big|^{\gamma-2} e(\tilde
u_{{\eps}\tau}^k)\Big){:}e(\testu)\,\d x
\\&&\nonumber
+\int_{\GC} \Big( {\dela}\tilde z_{{\eps}\tau}^k \JUMP{\tilde
u_{{\eps}\tau}^k}{} + \yosd(\JUMP{\tilde u_{{\eps}\tau}^k}{}) +
\tau^{\beta}\Big(1+ \big|\JUMP{\tilde u_{{\eps}\tau}^k}{}\big|^2
\Big)^{\mu/2-1} \JUMP{\tilde u_{{\eps}\tau}^k}{} \Big) {\cdot}
\JUMP{\testu}{}\,\d S
\\&&\label{GM1d-mod}
=\int_\Omega\big(1{-}\sqrt\tau\,\big) \mathbb De\big(\dt
u_{{\eps}\tau}^k\big){:} e(\testu)
-\Theta(\w_{{\eps}\tau}^k)\bbB{:}e(\testu) + \FRM_{\tau}^k
{\cdot} \testu  \,\d x +\int_{\Gnew} \fRM_{\tau}^k {\cdot} \testu \, \d
S
\end{eqnarray}
and for all $\tilde{z} \in L^\infty(\GC)$
\begin{align}
\nonumber
 & \!\!\!\!\! \!\!\!\!\! \int_{\GC} \zeta_1(\tilde z)  + \nabla\tilde
z_{{\eps}\tau}^k {\cdot} \nabla \left(\tilde z-\frac{\tilde
z_{{\eps}\tau}^k{-}z_{{\eps}\tau}^{k-1}}{\tau} \right)\, \d S  \hspace*{2em} \\
&  \!\!\!\!\! \!\!\!\!\! +  \int_{\GC}  \left(\tau^\alpha \tilde
z_{{\eps}\tau}^k + r(\tilde z_{{\eps}\tau}^k) + \frac\dela{2}\left|
\JUMP{\tilde u_{{\eps}\tau}^k}{} \right|^2-a_0
 \right)\left(\tilde z-\frac{\tilde
z_{{\eps}\tau}^k{-}z_{{\eps}\tau}^{k-1}}{\tau}\right)\,\d S
\ge\int_{\GC}\zeta_1 \Big(\frac{\tilde
z_{{\eps}\tau}^k{-}z_{{\eps}\tau}^{k-1}}{\tau}\Big)\,\d S\,.
\label{GM2d-mod}
\end{align}
\end{subequations}
Now,  we test the difference of (\ref{GM1a}) and (\ref{GM1d-mod}) by
$u_{{\eps}\tau}^k-\tilde u_{{\eps}\tau}^k$ and the difference of
\eqref{BC-conta} and (\ref{GM2d-mod}) by $z_{{\eps}\tau}^k-\tilde
z_{{\eps}\tau}^k$ and sum up the resulting relations. Using that the
underlying potential, namely the functional
 \[ (u,z)\mapsto\Phi_{{\eps}\tau}(u,z)
+\int_{\GC}\!\!
\zeta_1(z{-}z_{{\eps}\tau}^{k-1})\,\d S
+\int_\Omega\!\frac{\bbD e(u){:} e(u)}{2\sqrt\tau}\,\d x,
\]
is strictly convex on $W_{\Gdir}^{1,2}(\Omega {\setminus} \GC;\R^d)
\times L^\infty (\GC)$ by Lemma~\ref{l:semiconvex}, we conclude that
$ \uk=\tilde u_{{\eps}\tau}^k $, $\zk=\tilde z_{{\eps}\tau}^k.$
Then, the functional in (\ref{GMd-min1}) must have a bigger or equal
value on $(u_{{\eps}\tau}^{k-1},z_{{\eps}\tau}^{k-1})$ than on
$(\tilde u_{{\eps}\tau}^k,\tilde
z_{{\eps}\tau}^k)=(u_{{\eps}\tau}^k,z_{{\eps}\tau}^k)$, which gives
a discrete analog of \eqref{mech-energy}, namely
\begin{align}
 &T_\mathrm{kin}^{\varrho}\big(\dt u_{{\eps}\tau}^k
\big) +\Phi_{{\eps}\tau}\big(u_{{\eps}\tau}^k,z_{{\eps}\tau}^k\big)
+ \tau\int_{\GC}\!\! \zeta_1\big(\dt z_{{\eps}\tau}^k \big)\,\d S
+\frac{\tau}2\int_\Omega\big(2{-}\sqrt\tau\,\big)\bbD e\big(\dt
u_{{\eps}\tau}^k \big){:} e\big(\dt u_{{\eps}\tau}^k \big) \,\d x
\nonumber
\\
&\label{discrete-energy}
\le T_\mathrm{kin}^{\varrho}\big(\dt u_{{\eps}\tau}^{k-1} \big)
+\Phi_{{\eps}\tau}\big(u_{{\eps}\tau}^{k-1},z_{{\eps}\tau}^{k-1}\big)
-\tau\!\int_\Omega\! \Theta(\w_{{\eps}\tau}^k)\bbB{:}
e\big(\dt u_{{\eps}\tau}^k \big) + \FRM_{\tau}^k{\cdot}\dt
u_{{\eps}\tau}^k\d x
+\tau\!\int_{\Gnew}\!\!\!\fRM_{\tau}^k{\cdot}\dt u_{{\eps}\tau}^k \d S
\end{align}
when also employing  the algebraic inequality
\[
\dt^2u_{{\eps}\tau}^k{\cdot}\dt u_{{\eps}\tau}^k\ge\frac12|\dt
u_{{\eps}\tau}^k|^2 -\frac12|\dt u_{{\eps}\tau}^{k-1}|^2.
\]
 Upon summation
over $k$, we conclude  (\ref{disc-energy0}).

Now, to get (\ref{disc-energy}), we  add  to (\ref{discrete-energy})
the relation obtained testing (the weak formulation of) the boundary-value
problem (\ref{GM3a}, \ref{BC1b}, \ref{BC-contd},
\ref{BC-conte})  by $\tau$. Developing all calculations, one sees
that, thanks to our carefully designed discretization, the fourth
term on the left-hand side of~\eqref{discrete-energy} and the first
dissipative/adiabatic term on the  right-hand side of \eqref{GM3a}
mutually cancel out. So do the third term on the right-hand side
of~\eqref{discrete-energy} and the second right-hand-side term in \eqref{GM3a}.
Again upon summation over $k$, we arrive at~\eqref{disc-energy}.

Finally,  to check~(\ref{disc-semistability}), it just suffices to
realize that the functional minimized in (\ref{GMd-min1}) has a
lower value in $(u_{{\eps}\tau}^k,z_{{\eps}\tau}^k)$ than in
$(u_{{\eps}\tau}^k,\tilde z)$ for any $\tilde{z} \in L^\infty
(\GC)$, which gives
\begin{align}\nonumber
&\Phi_{{\eps}\tau}\big(u_{{\eps}\tau}^k,z_{{\eps}\tau}^k\big)+
\int_{\GC}\!\!
\tau\zeta_1\Big(\frac{z_{{\eps}\tau}^k{-}z_{{\eps}\tau}^{k-1}}{\tau}\Big)\,\d
S
\le\Phi_{{\eps}\tau}\big(u_{{\eps}\tau}^k,\tilde z\big)+
\int_{\GC}\!\!
\tau\zeta_1\Big(\frac{\tilde
z{-}z_{{\eps}\tau}^{k-1}}{\tau}\Big)\,\d S\,.
\end{align}
Then, by using that $\zeta_1$ is homogeneous degree $1$ and thus
satisfies the triangle inequality
$\zeta_1(\tilde{z}{-}z_{{\eps}\tau}^{k-1})\le \zeta_1(\tilde
z{-}z_{{\eps}\tau}^k)+\zeta_1(z_{{\eps}\tau}^k{-}z_{{\eps}\tau}^{k-1})$,
we find
\[
\begin{aligned}
\! \! \! \!
\Phi_{{\eps}\tau}\big(u_{{\eps}\tau}^k,z_{{\eps}\tau}^k\big)
\le \Phi_{{\eps}\tau}\big(u_{{\eps}\tau}^k,\tilde z\big)+
\int_{\GC}\!\!
\zeta_1\big(\tilde z{-}z_{{\eps}\tau}^{k-1}\big)
-\zeta_1\big(z_{{\eps}\tau}^k{-}z_{{\eps}\tau}^{k-1}\big)\,\d S
\le\Phi_{{\eps}\tau}\big(u_{{\eps}\tau}^k,\tilde z\big)+
\int_{\GC}\!\!
\zeta_1\big(\tilde z{-}z_{{\eps}\tau}^k\big)\,\d S\,.
\end{aligned}
\]
Being $k = 1,\ldots, K_\tau$ arbitrary, we
conclude~(\ref{disc-semistability}). $\hfill\Box$


\begin{lemma}[A priori estimates]
\label{prop:apriori}
 Under the  assumptions of Theorem~\ref{th:4.1}, there exist
 constants
$S_0$ and, for every $1\le r<\frac{d+2}{d+1}$, $S_r$ such that
for all $\varrho\geq 0$, $\eps
>0$ and for all $0 < \tau <\tau_\delam$
($\tau_\delam$ being as in Lemma~\ref{l:semiconvex}), for all
approximate solutions $(\pwc u {{\eps}\tau}, \pwc \w {{\eps}\tau},
\pwc z {{\eps}\tau}, \pwl u {{\eps}\tau}, \pwl \w {{\eps}\tau}, \pwl
z {{\eps}\tau})$ the following estimates  hold
\begin{subequations}\label{a-priori3}
\begin{align}
& \label{a30} \big\|\pwc
u{{\eps}\tau}\big\|_{L^{\infty}(0,T;W_{\Gdir}^{1,2}(\Omega;\R^d))}\le
S_0\,,
\\
& \label{a31} \big\|\pwl u{{\eps}\tau}\big\|_{
W^{1,2}(0,T;W_{\Gdir}^{1,2}(\Omega;\R^d))}\le S_0\,,
\\
& \label{a31bis} \varrho^{1/2} \big\|\pwl u{{\eps}\tau}\big\|_{
W^{1,\infty}(0,T;L^2(\Omega;\R^d)) }\le S_0\,,
\\
&  \label{a35} \big\|\pwc
u{{\eps}\tau}\big\|_{L^\infty(0,T;W_{\Gdir}^{1,\gamma}(\Omega;\R^d))}\le
\frac{S_0}{\sqrt[\gamma]{\tau}},
\\&
\label{a32} \big\|\pwc z{{\eps}\tau}\big\|_{L^\infty(\SC)}
\leq S_0\,,
\\
 & \label{a32bis}
 \big\|\pwl z{{\eps}\tau}\big\|_{\BV([0,T];L^1(\GC))} \leq S_0,
\\&
\label{a33} \big\|\pwc
\w{{\eps}\tau}\big\|_{L^\infty(0,T;L^1(\Omega))}
\le S_0,
\\
 & \label{a33bis}
\big\|\pwl \w{{\eps}\tau}\big\|_{ L^r(0,T;W^{1,r}(\Omega))} \leq
S_{{r}} \ \ \mbox{  for any $1\le r<\frac{d+2}{d+1}$},
\\&
\label{a34} \Big\|\pwl{\DT\w}{{\eps}\tau} \Big\|_{L^1(0,T;
W^{1,r'}(\Omega)^*)}\le S_0,
\\&
\label{a36} \varrho \,
\Big\|\pwl{\DT{u}}{{\eps}\tau}\Big\|_{\BV([0,T];W_{\Gdir}^{1,\gamma}(\Omega;\R^d)^*)}\le
S_0\,,
\end{align}
\end{subequations}
where $S_0$ and $S_r$ neither depend on $\eps$ nor on $\tau$.
Estimates \eqref{a32}, \eqref{a32bis}, \eqref{a33}, \eqref{a33bis}
respectively hold for $\pwl z{{\eps}\tau}$, $\pwc z{{\eps}\tau}$,
$\pwl \w{{\eps}\tau}$ and $\pwc \w{{\eps}\tau}$, as well.
\end{lemma}
\noindent{\it Proof.} Some of the calculations we shall develop
hereafter are analogous to the ones in the
 proof of~\cite[Prop.~4.2]{tr1}, to which we shall systematically refer.

First of all, we use the ``discrete total energy''
balance~\eqref{disc-energy}. Indeed, on the one hand,  by
definition~\eqref{Phi-eps-tau} of $\Phi_{{\eps}\tau}$, the second
term on the left-hand side of~\eqref{disc-energy}  is non-negative
and, thanks to
\COL{positive-definiteness of $\bbC$} and Korn's
inequality~\eqref{e:korn}, it provides a  bound for $\|\pwc u
{{\eps}\tau}(t) \|_{W^{1,2}(\Omega;\R^d)}^2$ and for $\tau\| \pwc u
{{\eps}\tau}(t) \|_{W^{1,\gamma}(\Omega;\R^d)}^\gamma $ uniformly
w.r.t. $t \in [0,T]$. Further,
 being $\pwc\w{{\eps}\tau} \geq
0$ a.e. in $\Omega$ thanks to~\eqref{posw}, the third term on the
left-hand side  of~\eqref{disc-energy} estimates $\|
\pwc\w{{\eps}\tau} \|_{L^\infty(0,T;L^1(\Omega))}$. To estimate the
right-hand side of~\eqref{disc-energy}, we employ the
 discrete ``by-part'' summation \cite[Formula~(4.51)]{tr1},
 to the effect that
\begin{align}
\nonumber &\!\!\!
\!\int_0^{\bar{\mathsf{t}}_{\tau}(t)}\!\!\int_{\Gnew} \pwc \fRM{\tau}(s)
{\cdot}\pwl{\DT{u}}{{\eps}\tau}(s)\,\d S\d s
 =  \int_{\Gnew}\pwc\fRM{\tau}(t){\cdot}\pwc
u{{\eps}\tau}(t) \, \d S - \int_{\Gnew}  \pwc \fRM{\tau}(\tau){\cdot}
u_{0,\tau}\, \d S \!\!
  - \int_\tau^{\bar{\mathsf{t}}_{\tau}(t)}\!\!\int_{\Gnew}\pwl
{\DT{\fRM}}{\tau}(s){\cdot}\upwc{u}{{\eps}\tau}(s) \, \d S \d s
\\ &\nonumber \qquad
\leq \rho_3 \|\pwc u{{\eps}\tau}(t) \|_{W^{1,2}(\Omega;\R^d)}^2 +
 C_{\rho_3} \left( \|u_{0,\tau} \|_{W^{1,2}(\Omega;\R^d)}^2 +
\| \pwc \fRM{\tau}\|_{L^\infty (0,T; L^{4/3}(\Gnew))}^2 \right)
\\\label{est:02} &\qquad\qquad
 + C\int_0^{\bar{\mathsf{t}}_{\tau}(t)} \|\pwl
{\DT{\fRM}}{\tau}(s)\|_{L^{4/3}(\Gnew)} \| \pwc{u}{{\eps}\tau}(s)
\|_{W^{1,2}(\Omega;\R^d)}\, \d s\,,
\end{align}
where inequality~\eqref{est:02} is also due to the continuous
embedding~\eqref{e:contemb} and $\rho_3$ is chosen in such a way as
to absorb the first term on the right-hand side into the term
$\|\pwc u {{\eps}\tau}(t) \|_{W^{1,2}(\Omega;\R^d)}^2$ on the
left-hand side of~\eqref{disc-energy}.
 Furthermore, in the case $\varrho >0$ we estimate the fourth term
 on
the right-hand side of~\eqref{disc-energy} by
\begin{equation}
\label{est:01} \int_0^{\bar{\mathsf{t}}_{\tau}(t)} \int_{\Omega}\pwc
\FRM{\tau}(s) {\cdot}  \pwl{\DT{u}}{{\eps}\tau}(s) \, \d x  \, \d s
\leq \int_0^{\bar{\mathsf{t}}_{\tau}(t)} \|\pwc\FRM{\tau}(s) \|_{L^2
(\Omega;\R^d)} \,\| \pwl{\DT{u}}{{\eps}\tau}(s)\|_{L^2
(\Omega;\R^d)}\, \d s\,.
\end{equation}
We then combine~\eqref{disc-energy}, \eqref{est:01},
and~\eqref{est:02}, and
use~\eqref{data-converg}--\eqref{data-converg-bis} for $\pwc
\FRM{\tau},\,  \pwc \fRM{\tau},\, \pwl {\DT{\fRM}}{\tau}, \, \pwc
\GRM{\tau},$ and $\pwc \gRM{\tau}$.
 Applying the Gronwall Lemma, we
conclude estimates \eqref{a30}, \eqref{a31bis}, \eqref{a35},
\eqref{a32},   and \eqref{a33} (the estimates for $\pwl
{z}{{\eps}\tau}$ and $\pwl {\w}{{\eps}\tau}$ following from the
bounds for $\pwc {z}{{\eps}\tau}$ and $\pwc {\w}{{\eps}\tau}$, and
from~\eqref{elementary1}). In the case $\varrho=0$, the only change
in the above calculations is that, under the additional
assumption~\eqref{effebis}, we estimate the fourth term
 on
the right-hand side of~\eqref{disc-energy} by use of the
aforementioned discrete by-part summation formula. Namely, on
account of the  Sobolev embedding~\eqref{e:contemb}
\begin{align}
\nonumber \!  \!  \!  \!  \!  \!
\int_0^{\bar{\mathsf{t}}_{\tau}(t)}&\!\!\int_{\Omega} \pwc
\FRM{\tau}(s){\cdot}\pwl{\DT{u}}{{\eps}\tau}(s)  \, \d x \d s =
\int_{\Omega} \pwc\FRM{\tau}(t){\cdot}\pwc u{{\eps}\tau}(t)\, \d x -
\int_{\Omega}  \pwc \FRM{\tau}(\tau){\cdot}u_{0,\tau}\, \d x -
\int_\tau^{\bar{\mathsf{t}}_{\tau}(t)}\!\!\int_{\Omega} \pwl
{\DT{\FRM}}{\tau}(s){\cdot}\upwc{u}{{\eps}\tau}(s)  \, \d x \d s
\\ &
\label{est:02-bis}
\begin{aligned}
\leq \rho_4 & \|\pwc u{{\eps}\tau}(t) \|_{W^{1,2}(\Omega;\R^d)}^2 +
 C_{\rho_4} \left( \|u_{0,\tau} \|_{W^{1,2}(\Omega;\R^d)}^2 +
\| \pwc \FRM{\tau}\|_{L^\infty (0,T; L^{6/5}(\Omega;\R^d))}^2
\right)
\\ &
 + C\int_0^{\bar{\mathsf{t}}_{\tau}(t)} \|\pwl
{\DT{\FRM}}{\tau}(s)\|_{L^{6/5}(\Omega;\R^d)} \|
\pwc{u}{{\eps}\tau}(s) \|_{W^{1,2}(\Omega;\R^d)}\, \d s\,,
\end{aligned}
\end{align}
where again the  positive constant $\rho_4$ is such that the first
term on the right-hand side of \eqref{est:02-bis} is controlled by
$\|\pwc u {{\eps}\tau}(t) \|_{W^{1,2}(\Omega;\R^d)}^2$ on the
left-hand side of~\eqref{disc-energy}.

  Secondly, again  arguing as for~\cite[Prop.~4.2]{tr1}, we make use
  of the technique by \textsc{Boccardo \& Gallou\"et}~\cite{boccardo-gallouet1}, with the simplification devised
  in~\cite{feireisl-malek}. Hence, we test  the heat equation~\eqref{weak-heat-discr}
by $\Pi (\pwc {\w} {{\eps}\tau})$, where $\Pi : [0,+\infty) \to
[0,1]$ is the map
\[
w \mapsto \Pi (w) = 1-\frac1{(1{+}w)^\varsigma},\ \ \ \ \varsigma>0;
\]
note that $\Pi (\pwc {\w} {{\eps}\tau})\in W^{1,2}(\Omega{\setminus}
\GC )$, because $\Pi$ is Lipschitz continuous. With the same
calculations as in~\cite{tr1}, taking into account~\eqref{30dprimo}
we find
\begin{align}\nonumber
&\hspace*{-1em}
\varsigma\, \mathsf{k}\,
\int_Q
\frac{|\nabla \pwc {\w} {{\eps}\tau}|^2}{(1+\pwc {\w}
{{\eps}\tau})^{1+\varsigma}}\,\d x\d t \le
\int_Q\mathcal{K}(e(\pwc {u} {{\eps}\tau}),\pwc {\w} {{\eps}\tau})
\nabla\pwc {\w} {{\eps}\tau}{\cdot}\nabla\Pi(\pwc {\w}
{{\eps}\tau})\,\d x\d t
\\\nonumber &\qquad+
\int_{\SC} \eta({\underline u_{\eps\tau},}\pwc z{{\eps}\tau})
\JUMP{\Theta(\pwc\w{{\eps}\tau})}{}\JUMP{\Pi(\pwc {\w}
{{\eps}\tau})}{}\,\d S\d t
 +\int_\Omega\widehat{\Pi}(\pwc {\w} {{\eps}\tau}(T,{\cdot}))\,\d x
\end{align}
\begin{align}\nonumber
 &\qquad \le\int_\Omega\!\widehat{\Pi}(\w_0)\,\d x + \| \pwc\GRM{\tau}
\|_{L^1(Q)} + \| \pwc\gRM{\tau}  \|_{L^1 (\Sigma))}
\\ &\label{est:2}\qquad
 +
C \big( \| \bbD e(\pwl{\DT{u}}{{\eps}\tau})  {:}
e(\pwl{\DT{u}}{{\eps}\tau}) \|_{L^1 (Q)} +
\|\Theta(\pwc\w{{\eps}\tau})\bbB{:}e(\pwl{\DT{u}}{{\eps}\tau})
\|_{L^1 (Q)} + \| \zeta_1 (\pwl {\DT{z}}{{\eps}\tau})\|_{L^1 (\SC)}
\big)
\end{align}
where $\widehat{\Pi}$ is the primitive function of $\Pi$ such that
$\widehat{\Pi}(0)=0$. Note that inequality~\eqref{est:2} follows
from the  fact that $\eta({\underline u_{\eps\tau},}\pwc
z{{\eps}\tau}) \JUMP{\Theta(\pwc\w{{\eps}\tau})}{}\JUMP{\Pi(\pwc
{\w} {{\eps}\tau})}{} \geq 0$ a.e. in $\SC$ (by the positivity of
$\eta$
and the monotonicity of $\Theta$ and $\Pi$), from the ``discrete
chain rule'' \cite[Formula~(4.30)]{tr1} for $\widehat{\Pi}$, and
from  $0 \leq \Pi(\pwc {\w} {{\eps}\tau}) \leq 1$.
Combining~\eqref{est:2} with the Gagliardo-Nirenberg inequality,
 we find, for all $1\leq r <{(d{+}2)}/{(d{+}1)}$, that
\begin{align}
 \big\|\nabla \pwc {\w}
{{\eps}\tau} \big\|_{L^r(Q;\R^d)}^r
 \leq  C_r\big(&1+\| \bbD
e(\pwl{\DT{u}}{{\eps}\tau}){:}e(\pwl{\DT{u}}{{\eps}\tau})
\|_{L^1 (Q)}\!
\label{est:3}
 + \|\Theta(\pwc\w{{\eps}\tau})\bbB{:}e(\pwl{\DT{u}}{{\eps}\tau})
\|_{L^1 (Q)}\!+ \| \zeta_1 (\pwl {\DT{z}}{{\eps}\tau})\|_{L^1
(\SC))}\big)
\end{align}
for some positive constant $C_r$,  depending on $r$ and also on  the
function $\eta$, cf.\ \eqref{eta-affine}.

Then, we multiply \eqref{est:3} by a  constant  $\rho_5>0$ and add
it to~\eqref{disc-energy0} (in which we set $t=T$). Now, by
\COL{positive-definiteness of $\bbC$}, the third term on the
left-hand side of~\eqref{disc-energy0}
 is bounded from below by
 $\mathsf{d}/2 \|e(\pwl{\DT{u}}{{\eps}\tau}) \|_{L^2(Q;\R^{d\times d})}^2 $, whereas the
 fourth term controls $ \| \zeta_1 (\pwl {\DT{z}}{{\eps}\tau})\|_{L^1 (\SC)}$.
  Thus,
we choose $\rho_5$ small enough in such a way to absorb the first
and the third term on the right-hand side of~\eqref{est:3} into the
left-hand side of~\eqref{disc-energy0}. Hence, we find
\begin{align}
\frac{\mathsf{d}}4  &
\|e(\pwl{\DT{u}}{{\eps}\tau}) \|_{L^2 (Q; \R^{d\times d})}^2
+ (1{-}\rho_5)
\|\zeta_1 (\pwl {\DT{z}}{{\eps}\tau})\|_{L^1 (\SC)}
+\rho_5
\big\|\nabla \pwc {\w} {{\eps}\tau} \big\|_{L^r(Q;\R^d)}^r
\leq   T_\mathrm{kin}^{\varrho}\big(\DT{u}_{0,\tau})
+\Phi_{{\eps}\tau}\big(u_{0,\tau},z_{0,\tau}) \nonumber\\ & +
\int_0^T\int_{\Omega}\pwc \FRM{\tau}{\cdot} \pwl{\DT{u}}{{\eps}\tau}
\, \d x \, \d t
+\int_0^T\int_{\Gnew} \pwc \fRM{\tau} {\cdot} \pwl{\DT{u}}{{\eps}\tau}
\,  \d S \d t
 + (\rho_5 C_r{+}1)
\|\Theta(\pwc\w{{\eps}\tau})\bbB e(\pwl{\DT{u}}{{\eps}\tau})
\|_{L^1(Q)}.
\label{disc-energy01}
\end{align}
The first two terms on the right-hand side of~\eqref{disc-energy01}
are estimated in view of
\eqref{data-converg}--\eqref{data-converg-bis}
and~\eqref{est-init-data}, whereas, taking into account the Sobolev
embedding~\eqref{e:contemb} and Korn's inequality~\eqref{e:korn}, we
have
\[
\begin{aligned} &
\int_0^T\int_{\Omega}\pwc \FRM{\tau}{\cdot} \pwl{\DT{u}}{{\eps}\tau}
\, \d x  \d t \leq C \int_0^T \| \pwc
\FRM{\tau}(s)\|_{L^{6/5}(\Omega;\R^d)}^2 + \frac{\mathsf{d}}{16}
\int_0^T \|e(\pwl{\DT{u}}{{\eps}\tau}) \|_{L^2 (\Omega; \R^{d\times
d})}^2\, \d t,
\\
 &
\int_0^T\int_{\Gnew} \pwc \fRM{\tau} {\cdot} \pwl{\DT{u}}{{\eps}\tau}
\,\d S \d t \leq C \int_0^T \| \pwc
\fRM{\tau}(s)\|_{L^{4/3}(\Gnew;\R^d)}^2 + \frac{\mathsf{d}}{16}
\int_0^T \|e(\pwl{\DT{u}}{{\eps}\tau}) \|_{L^2 (\Omega; \R^{d\times
d})}^2\, \d t;
\end{aligned}
\]
\COL{here again $\mathsf{d}>0$ is the positive-definiteness constant
of $\bbD$.} To estimate the last summand, we use
\COL{the positive-definiteness of $\bbC$} and~\eqref{growthTheta},
finding
\begin{align}
(\rho_5 C_r{+}1)\|\Theta(\pwc\w{{\eps}\tau})\bbB{:}
e(\pwl{\DT{u}}{{\eps}\tau}) \|_{L^1(Q) }
& \leq \rho_6
|\bbD|\,\| e(\pwl{\DT{u}}{{\eps}\tau}) \|_{L^2(Q;
\R^{d\times d}))}^2 + C_{\rho_6} \|\Theta(\pwc\w{{\eps}\tau})
\|_{L^2(Q)}^2 \nonumber\\ & \leq \rho_6
|\bbD|\,\|e(\pwl{\DT{u}}{{\eps}\tau})\|_{L^2(Q; \R^{d\times d}))}^2 +
C_{\rho_6}\big(\|\pwc\w{{\eps}\tau}
\|_{L^{2/\omega}(Q)}^{2/\omega}\!+ 1\big) \label{est:4}
\end{align}
in which we choose  the positive constant $\rho_6$ small enough,
again to absorb the first term on the right-hand side
of~\eqref{est:4} into the left-hand side of~\eqref{disc-energy01}.
In order to estimate $\|\pwc\w{{\eps}\tau} \|_{L^{2/\omega}(Q)}$, we
again employ the Gagliardo-Nirenberg inequality. Indeed, with the
same calculations as throughout~\cite[Formulae~(4.39)--(4.43)]{tr1},
and relying on the restriction of $\omega$ in \eqref{30b} and on the
bound for $\| \pwc \w {{\eps}\tau}\|_{L^\infty (0,T;L^1 (\Omega))}$,
we conclude
\begin{equation}
\label{e:plugin} \int_0^T \|\pwc\w{{\eps}\tau}
\|_{L^{2/\omega}(\Omega))}^{2/\omega} \leq  \rho_7 \int_0^T
\big\|\nabla \pwc {\w} {{\eps}\tau} \big\|_{L^r(\Omega;\R^d)}^r +
C_{\rho_7}
\end{equation}
for a suitably small $\rho_7>0$. Then, we plug~\eqref{e:plugin}
into~\eqref{est:4}, and the latter into~\eqref{disc-energy01}, and
choose $\rho_6$ in such a way as to absorb $\|\nabla \pwc {\w}
{{\eps}\tau} \big\|_{L^r(Q;\R^d)}^r$ into the left-hand side
of~\eqref{disc-energy01}. Thus, we conclude estimate~\eqref{a31}, as
well as an estimate for
$\|\zeta(\pwl{\DT{z}}{{\eps}\tau})\|_{L^1(\SC)}$
(yielding~\eqref{a32bis}),  and a bound for $\nabla \pwc {\w}
{{\eps}\tau}$ in $L^r(Q;\R^d)$. Combining the latter information
with the estimate for $\pwc {\w} {{\eps}\tau}$ in $L^\infty (0,T;
L^1(\Omega))$, we infer~\eqref{a33bis} (the estimate for $\pwl
{\w}{{\eps}\tau}$  due to the bound for $\pwc \w{{\eps}\tau}$ and
to~\eqref{elementary1}). As a by-product of the above calculations,
we find
\begin{equation}
\label{e:later} \|\pwl\Lambda {{\eps}\tau}\|_{L^1(Q) }\leq C\,, \ \
\text{with} \ \ \pwl\Lambda {{\eps}\tau}: =
\frac{2-\sqrt{\tau}}{2}\bbD e(\pwl{\DT{u}}{{\eps}\tau}){:}
e(\pwl{\DT{u}}{{\eps}\tau}) +
\Theta(\pwc\w{{\eps}\tau})\bbB{:}e(\pwl{\DT{u}}{{\eps}\tau})\,.
\end{equation}
For later convenience, we also remark that \eqref{a33bis} yields
\begin{equation}
\label{e:later2} \| \Theta(\pwc\w{{\eps}\tau})\|_{L^{\omega r}(0,T;
L^{\omega q}(\GC))} \leq C \quad \text{for all $1\leq q \leq
\frac{dr}{d-r}$}\,,
\end{equation}
where we have also used  the continuous embedding $W^{1,r}(\Omega)
\subset L^q (\GC)$ for $q$ ranging in the above-mentioned index
interval, as well as the growth restriction~\eqref{growthTheta}
imposed on $\Theta$.

To prove  \eqref{a34}, we argue by comparison
in~\eqref{weak-heat-discr}, to the effect that
\[
\| \pwl {\DT{\w}}{{\eps}\tau}\|_{L^1 (0,T; W^{1,r'}(\Omega)^*)} =
\sup_{
\|\testw\|_{L^\infty (0,T; W^{1,r'}(\Omega))} \leq 1}\left(I_5+ I_6
+ I_7 + I_8 \right), \quad  \text{where}
\]
\[
I_5 =  \int_Q \pwl\Lambda {{\eps}\tau} \testw \leq \|\pwl\Lambda
{{\eps}\tau}\|_{L^1(Q)} \| \testw\|_{L^\infty(Q)} \leq C
\|\testw\|_{L^\infty (0,T;W^{1,r'}(\Omega))}
\]
thanks to \eqref{e:later} and the continuous embedding
$W^{1,r'}(\Omega) \subset L^\infty (\Omega)$ (since $r'>d+2$), while
\[
\begin{aligned}
I_6 =  - \int_Q \mathcal{K}(e(\pwc
u{{\eps}\tau}),\pwc\w{{\eps}\tau})\nabla\pwc\w{{\eps}\tau}{\cdot}\nabla
\testw \leq C_{\mathcal{K}} \|\nabla \pwc\w{{\eps}\tau}\|_{L^r (Q)}
\| \nabla \testw\|_{L^{r'} (Q)}\leq C
\|\testw\|_{L^\infty(0,T;W^{1,r'}(\Omega))}
\end{aligned}
\]
due to~\eqref{growthKappa} and~\eqref{a33}. Further,
\[
\begin{aligned}
I_7 &  =  -\int_{\SC} \Big(\pwc z{{\eps}\tau}
\JUMP{\Theta(\pwc\w{{\eps}\tau})}{}\JUMP{\testw}{}
 +
 a_1\pwl{\DT{z}}{{\eps}\tau}\frac{
\testw|_{\GC}^+{+}  \testw|_{\GC}^-}2 \Big)
 \,\d S\d t
\\ &  \leq \left(  \| \pwc z{{\eps}\tau} \|_{L^\infty(\SC)}  \,  \|
\Theta(\pwc\w{{\eps}\tau})\|_{L^{\omega r}(0,T; L^{\omega q}(\GC))}
\, + \| \zeta_1 (\pwl{\DT{z}}{{\eps}\tau})\|_{L^1 (\SC)}\right)\, \|
\testw\|_{L^\infty(\SC)}
\leq C \|\testw\|_{L^\infty(0,T;W^{1,r'}(\Omega))}
\end{aligned}
\]
thanks to~\eqref{a32}, 
\eqref{e:later2}, and  the
continuous embedding $W^{1,r'}(\Omega) \subset L^\infty (\GC)$, and,
finally,
\[
\begin{aligned}
I_8   =  \int_{Q}  \pwc\GRM{\tau}  \testw\,\d x \d t +
  \int_{\Sigma} \pwc\gRM{\tau} \testw\,\d S \d t  \leq
  \left( \| \pwc\GRM{\tau}
\|_{L^1(Q)} + \| \pwc\gRM{\tau}  \|_{L^1 (\Sigma)}
\right)\|\testw\|_{L^\infty (0,T;W^{1,r'}(\Omega))}\,.
\end{aligned}
\]
Collecting the above calculations, we conclude~\eqref{a34}.

Finally, for~\eqref{a36} we use that $\pwl{\DDT{u}}{{\eps}\tau}$ is
a measure on $[0,T]$, supported at the jumps of
$\pwl{\DT{u}}{{\eps}\tau}$, and we estimate  the norm
 $\varrho \|\pwl{\DDT{u}}{{\eps}\tau}\|_{\mathrm{M}(0,T;W_{\Gdir}^{1,\gamma}(\Omega{{{\setminus}}}\GC;\R^d)^*)}$,
where $\mathrm{M}(0,T;W_{\Gdir}^{1,\gamma}(\Omega{{{\setminus}}}\GC;\R^d)^*)$
denotes the space of Radon measures on $[0,T]$ with values in
 $W_{\Gdir}^{1,\gamma}(\Omega{{{\setminus}}}\GC;\R^d)^*$,
arguing by comparison in~\eqref{e:discrmom}; see the proof
of~\cite[Prop.~4.2]{tr1}, where similar calculations were carried
out.
 $\hfill\Box$

\section{Limit passage with
$\tau\to0$ and proof of Theorem~\ref{th:4.1}} \label{ss:4.4}
Throughout this section, we shall keep $\eps>0$ fixed, and let $\tau
\to 0$. We shall  develop a proof of the passage to the limit
unifying  the cases $\varrho>0$ and $\varrho=0$.
\medskip
\par\noindent {\textbf{Step $0$: selection of convergent subsequences.}}
First of all, it follows from estimates~\eqref{a31}, \eqref{a31bis},
and \eqref{a36}, from the Banach selection principle, and from the
Aubin-Lions theorem (see, e.g., \cite[Thm.~5,Cor.~4]{simon86}
and~\cite[Cor.~7.9]{NPDE_roubicek} for the generalization to the
case of time derivatives as measures), that there exist a (not
relabeled) sequence $\tau \to 0$
 and a  limit function
$\ue \in W^{1,2}(0,T;W_{\Gdir}^{1,2}(\Omega{{\setminus}} \GC;\R^d))$
such that the following  weak,  weak$^*$, and strong  convergences
hold as $\tau\to 0$:
\begin{subequations}
\label{e:convutau}
\begin{align}
 \label{e:convutau1} & \pwl {u}{{\eps}\tau} \weakto
\ue  \ \  \text{ in $
W^{1,2}(0,T;W_{\Gdir}^{1,2}(\Omega{{\setminus}} \GC;\R^d)),$}
\\
 \label{e:convutau2}
 & \pwl {u}{{\eps}\tau} \to
\ue  \ \ \text{ in $
\mathrm{C}^{0}([0,T];W_{\Gdir}^{1{-}\epsilon,2}(\Omega{{\setminus}}
\GC;\R^d))$} \ \ \forall\, \epsilon \in
 (0,1],
 \\
\label{e:convutau3} &  \text{if $\varrho>0$,} \  \ \pwl
{u}{{\eps}\tau} \weaksto \ue \ \  \text{ in $
W^{1,\infty}(0,T;L^2(\Omega;\R^d))$.}
\end{align}
In the case $\varrho>0$ we also have  $\pwl {u}{{\eps}\tau} \to
u_\eps$    in
$W^{1,2}(0,T;W_{\Gdir}^{1{-}\epsilon,2}(\Omega{{\setminus}}
\GC;\R^d))  \cap W^{1,q}(0,T;L^2(\Omega;\R^d)) $ for all $\epsilon
\in (0,1]$ and $1 \leq q <\infty$.  Furthermore,
estimate~\eqref{a36} and a generalization of Helly's principle (see
\cite{BarbuPrecupanu86} as well as \cite[Thm.~6.1]{MieThe04RIHM})
yield that $\DT{u}_\eps \in \BV
([0,T];W_{\Gdir}^{1,\gamma}(\Omega{{\setminus}} \GC;\R^d)^*)$ and,
in addition, $ \pwl {\DT{u}}{\eps \tau}(t) \weakto \DT{u}_{\eps}(t)$
in $W_{\Gdir}^{1,\gamma}(\Omega{{\setminus}} \GC;\R^d)^*$ for all $t
\in [0,T]$. By virtue  of estimate~\eqref{a31bis} and of a trivial
compactness argument, this pointwise weak convergence improves to
\begin{equation}
\label{pointiwise-for-u} \pwl {\DT{u}}{\eps \tau}(t) \weakto
\DT{u}_{\eps}(t) \quad \text{in $L^2(\Omega;\R^d)$ for all $t \in
[0,T]$, in the case $\varrho>0$.}
\end{equation}
Combining \eqref{e:convutau1} and \eqref{e:convutau2} with the
general inequality~\eqref{elementary2}, we conclude that, up to the
extraction of a further subsequence, for all $\epsilon \in(0,1]$,
\begin{equation}
\label{e:convutau5}
\begin{aligned}
& \pwc {u}{{\eps}\tau} \weaksto u_\eps \ \ \text{ in
 $L^{\infty}(0,T;W_{\Gdir}^{1,2}(\Omega{{\setminus}} \GC;\R^d))$,}
 \quad
 \pwc {u}{{\eps}\tau} \to u_\eps  \ \ \text{ in
$L^{\infty}(0,T;W_{\Gdir}^{1{-}\epsilon,2}(\Omega{{\setminus}}
\GC;\R^d))$,}
 \\ &
\pwc {u}{{\eps}\tau}(t) \to \ue(t) \ \ \text{in
$W_{\Gdir}^{1{-}\epsilon,2}(\Omega{{\setminus}} \GC;\R^d))$ for all
$t \in [0,T]$,}
 \end{aligned}
\end{equation}
\end{subequations}
the latter pointwise convergence due to~\eqref{e:convutau2}
and~\eqref{elementary2}.

With the aforementioned compactness results, we deduce
from estimates \eqref{a32},  and \eqref{a32bis} that there exists a
function $\ze \in L^\infty(\SC) \cap \BV ([0,T];
\mathcal{Z})$, ($\mathcal{Z}$ being some reflexive space  such that
$L^1(\GC) \subset \mathcal{Z}$ with a continuous embedding, for
example $\mathcal{Z} = W^{1,2+\epsilon}(\GC)^*$ for some
$\epsilon>0$), such that (possibly along a subsequence)
\begin{subequations}
\label{e:convztau}
\begin{align}
\label{e:convztau1}
    \pwc {z}{{\eps}\tau}
\weaksto z_\eps  \ \ \text{ in $L^{\infty}(\SC)$,}
\end{align}
and, again by~\cite[Thm.~6.1, Prop.~6.2]{MieThe04RIHM}, $ \pwc
{z}{{\eps}\tau}(t)\weakto \ze(t)$ in $\mathcal{Z}$ for all $t \in
[0,T]$. In view of~\eqref{a32}, we indeed have pointwise weak$^*$
convergence in $L^\infty (\GC)$, i.e.
\begin{equation}
\label{e:convztau2} \pwc {z}{{\eps}\tau}(t)\weaksto \ze(t) \ \
\text{ in $L^\infty (\GC)$ for all $t \in [0,T]$.}
\end{equation}

 Finally,
 arguing as in the proof of~\cite[Thm.~6.1]{MieThe04RIHM} with
Helly's selection principle and taking into
account~\eqref{e:convztau2},  we conclude that for all $0 \leq s
\leq t \leq T$
\begin{equation}
\label{e:finally} \mathrm{Var}_{\mathcal{R}}(z_\eps;[s,t]) \leq
\lim_{\tau \to 0} \int_s^t\!\int_{\GC}\!\zeta_1
\big(\pwl{\DT{z}}{{\eps}\tau}(r)\big)\, \d S \d r\,,
\end{equation}
which ultimately yields that  there exists $S_0'>0$ such that
\begin{equation}
\label{bvl1}
 \ze \in \BV ([0,T]; L^1 (\GC))\ \ \ \text{and}\  \ \ \|\ze \|_{\BV ([0,T]; L^1
 (\GC))} \leq S_0' \quad \text{for all $\eps>0$}.
\end{equation}
\end{subequations}

Thirdly, by the same tokens we conclude from estimates \eqref{a33},
\eqref{a33bis}, and \eqref{a34} that there exists $\we \in
L^{r}(0,T;W^{1,r}(\Omega{{\setminus}} \GC)) \cap \BV ([0,T];
W^{1,r'}(\Omega{{\setminus}} \GC)^*)$ such that
\begin{subequations}
\label{e:convwtau}
\begin{align}
 \label{e:convwtau1} & \pwc {\w}{{\eps}\tau}, \, \pwl {\w}{{\eps}\tau}
\rightharpoonup \w_\eps \ \ \text{ in $
L^{r}(0,T;W^{1,r}(\Omega{{\setminus}} \GC))$,}
\\
 \label{e:convwtau2} &
\pwc {\w}{{{\eps}\tau}}, \,  \pwl {\w}{{\eps}\tau} \to \w_\eps \ \
\text{ in $ L^{r}(0,T;W^{1,r-\epsilon}(\Omega{{\setminus}} \GC))
\cap L^q (0,T; L^1 (\Omega))$}\,,
\end{align}
for all $\epsilon \in
 (0,r-1]$ and $1 \leq q <\infty$,  as well as
\begin{equation}
\label{e:poinwtiwise-w}
 \pwc {\w}{{{\eps}\tau}}(t), \, \pwl {\w}{{\eps}\tau}(t)
\weakto \w_{\eps}(t) \quad \text{in $W^{1,r'}(\Omega{{\setminus}}
\GC)^*$ for all $t \in [0,T]$.}
\end{equation}
Notice that, under condition \eqref{strict-pos} on $\theta_0$,
convergence~\eqref{e:convwtau2} and~\eqref{posw}
yield~\eqref{poswe}. It also follows from
\cite[Thm.~6.1]{MieThe04RIHM} that $ \mathrm{Var}_{\|
\cdot\|_{W^{1,r'}(\Omega)^*}}(\we;[s,t]) \leq \lim_{\tau \to 0}
\mathrm{Var}_{\| \cdot\|_{W^{1,r'}(\Omega)^*}}(\pwl
\w{{\eps}\tau};[s,t])\,, $ with $\mathrm{Var}_{\|
\cdot\|_{W^{1,r'}(\Omega)^*}}$ denoting the total variation w.r.t.
the norm $\| \cdot\|_{W^{1,r'}(\Omega)^*}$. This entails that
\begin{equation}
\label{bvw1}
 \|\we \|_{\BV ([0,T];W^{1,r'}(\Omega)^*)} \leq S_0' \quad \text{for all $\eps >0$}.
\end{equation}
For later purposes, we also point out that, in view of
estimate~\eqref{a33} and of~\eqref{e:convwtau2}, there holds
\begin{equation}
\label{useful-later} \| \w_\eps\|_{L^\infty (0,T; L^1 (\Omega))}\leq
S_0 \qquad \quad \text{for all $\eps>0$},
\end{equation}
$S_0$ being the same constant as in estimates~\eqref{a-priori3}.
\end{subequations}

Besides, \eqref{a35} yields that
\begin{subequations}
 \label{pass-limi-u}
\begin{align}
\label{conv1} \tau\big\||e(\pwc u {{\eps}\tau})|^{\gamma-2}e(\pwc u
{{\eps}\tau})\big\|_{L^{\gamma/{(\gamma-1)}}(Q; \R^{d\times d})} \le
S_0{\tau^{1/\gamma}} \to 0 \quad \text{ as $\ \tau\to0$.}
\end{align}
In view of~\eqref{e:contemb} and the second of~\eqref{e:convutau5},
it is not difficult to verify that, for all $\epsilon \in (0,3]$,
\begin{align}
\label{conv2} \left.
\begin{array}{lll}
 \JUMP{\pwc u{{\eps}\tau}}{} \to \JUMP{u_{\eps}}{} & \text{in $L^\infty
(0,T; L^{4-\epsilon} (\GC;\R^d))$,} \\   \JUMP{\pwc
u{{\eps}\tau}(t)}{} \to \JUMP{u_{\eps}(t)}{}  &  \text{in
$L^{4-\epsilon} (\GC;\R^d)$ for all $t \in [0,T]$.}
\end{array}
 \right\}
\end{align}
Furthermore, using that $\yosd$ is given by \eqref{yos-repre}, and
recalling~\eqref{e:contraction}, from \eqref{a31}  we  easily infer
that
\begin{equation}
\label{later-useful} \exists\, S_1= S_1(\eps)>0 \  \  \forall
\tau>0:\quad \big\|\yosd\big(\JUMP{\pwc u{{\eps}\tau}}{}\big)
\big\|_{L^\infty (0,T;L^2(\GC;\R^d))}\leq S_1,
\end{equation}
with $S_1(\eps)\to \infty $ as $\eps\to0$; more specifically, due
to~\eqref{e:contraction} we have
$S_1(\eps)=\mathscr{O}(1/\sqrt\eps)$. Combining \eqref{conv2} with
the strong-weak closedness  of the graph of the operator $\yosd$, up
to the extraction of a further subsequence we find that
\begin{equation}
\label{e:sweak} \yosd\big(\JUMP{\pwc u{{\eps}\tau}}{}\big)\weaksto
\yosd\big(\JUMP{u_{\eps}}{}\big)\quad\text{in
$L^\infty(0,T;L^2(\GC;\R^d))$.}
\end{equation}
 Moreover, using that
$ \big|\big(1+\big|\JUMP{\pwc
u{{\eps}\tau}}{}\big|^2\big)^{\frac\mu2-1} \JUMP{\pwc
u{{\eps}\tau}}{}\big|\le 2^{\frac\mu2 -2}\big(\big|\JUMP{\pwc
u{{\eps}\tau}}{}\big|+ \big|\JUMP{\pwc
u{{\eps}\tau}}{}\big|^{\mu-1}\big)$  a.e. in $\SC$, as well as
estimate~\eqref{a31}, one sees that the sequence  $\{ (1+ |
\JUMP{\pwc u{{\eps}\tau}}{} |^2 )^{\frac\mu2-1} \JUMP{\pwc
u{{\eps}\tau}}{}\}$ is bounded in $L^\infty (0,T;
L^{4/{(\mu-1)}}(\GC; \R^d))$. Thus,
\begin{equation}
\label{conv4} \tau^{\beta/2}\big(1+ \big| \JUMP{\pwc
u{{\eps}\tau}}{} \big|^2\big)^{\frac\mu2-1}\JUMP{\pwc
u{{\eps}\tau}}{}\to0\quad \text{in $L^\infty (0,T;
L^{4/{(\mu-1)}}(\GC; \R^d))$.}
\end{equation}
\end{subequations}

 Next,  let us point out that, in the case the space
dimension is $d=3$, \eqref{a33bis} holds for all $1 \leq r<5/{4}$,
so that \eqref{e:convwtau2} yields by interpolation
\begin{subequations}
 \label{pass-limi-w}
\begin{equation}
\label{e:new-label} \pwc {\w}{{\eps}\tau} \to \w_\eps \qquad \text{
in $L^{15/7-\epsilon} (Q)\ $ for all }\ \epsilon \in
\big(0,\mbox{$\frac87$}\big].
\end{equation}
 In particular, $\Theta(\pwc {\w}{{\eps}\tau}) \to \Theta (\w_\eps)$ a.e. in $Q$. Combining this
information with~\eqref{growthTheta} (note that, by~\eqref{30b},
$\omega
> \frac65$ for $d=3$),
 it is immediate to deduce from \eqref{e:new-label}
that, for example,
\begin{equation}
\label{conv3} \Theta(\pwc {\w}{{\eps}\tau})\to\Theta(\w_\eps)\qquad
\text{ in $L^{2} (Q)$}.
\end{equation}
Furthermore,  using  standard trace theorems we also deduce from
\eqref{e:convwtau2}  that for all $\epsilon \in  (0,\frac37]$
\[
 \pwc {\w}{{\eps}\tau}^+|_{\GC}  \to  \w_\eps^+|_{\GC}  \quad \text{and} \quad
  \pwc {\w}{{\eps}\tau}^-|_{\GC} \to  \w_\eps^-|_{\GC} \qquad\text{ in $
L^{r}(0,T;L^{10/7-\epsilon}(\GC))$,}
\]
so that, again by~\eqref{growthTheta}, for $d=3$,  using that
$\omega
> 6/5$, we conclude  that
\begin{equation}
\label{e:interesting}
 \JUMP{\Theta(\pwc {\w}{{\eps}\tau})}{}  \to  \JUMP{\Theta(\w_\eps)}{}\qquad
\text{ in $L^{r}(0,T;L^{12/7-\epsilon}(\GC)) $}  \ \ \forall\, \epsilon \in
\big(0,\mbox{$\frac57$}\big]\,.
\end{equation}
Similar calculations leading to~\eqref{conv3}
and~\eqref{e:interesting} can be performed in the case $d=2$.
\end{subequations}

In the end, we are now going to  show that
\begin{equation}
\label{e:altogether} \Phi_{\eps}(u_{\eps}(t), {z}_\eps(t)) \leq
\liminf_{\tau \to 0} \Phi_{{\eps}\tau}(\pwc u{{\eps}\tau}(t), \pwc
z{{\eps}\tau}(t)) \quad \text{for all $t \in [0,T]$}\,.
\end{equation}
 Indeed, taking into account~\eqref{e:convutau5} it is not difficult to deduce that for all $t
\in [0,T]$
\[
\liminf_{\tau \to 0} \int_{\Omega{{\setminus}}\GC} \frac12\bbC
e\big(\pwc u{{\eps}\tau}(t)\big){:}e\big(\pwc
u{{\eps}\tau}(t)\big)\,\d x\ge\int_{\Omega{{\setminus}}\GC}
\frac12\bbC e\big(\ue(t)\big){:}e\big(\ue(t)\big) \,\d x.
\]
Combining~\eqref{e:convztau2} with~\eqref{conv2},  we have
\[
\liminf_{\tau \to 0} \int_{\GC} \frac{\delam}2
\pwc{z}{{\eps}\tau}(t)\big|\JUMP{\pwc u{{\eps}\tau}(t)}{}\big|^2\,
\d S \geq
 \int_{\GC} \frac{\delam}2
\ze(t)\big|\JUMP{\ue(t)}{}\big|^2\, \d S \quad \text{for all $t \in
[0,T]$}\,.
\]
Besides, taking into account that $\yosapp$ is lower semicontinuous
on $L^2 (\GC;\R^d)$ (cf.~\eqref{F-eps}),
we immediately conclude
\[
\liminf_{\tau \to 0}
 \yosapp(\JUMP{\pwc u{{\eps}\tau}(t)}{}) \geq
 \yosapp(\JUMP{\ue(t)}{}) \quad \text{for all $t \in [0,T]$}\,.
\]
Collecting the above inequalities and also relying
on~\eqref{e:convztau2}, we infer~\eqref{e:altogether}.
\par\noindent \textbf{Step $1$: passage to the limit in the momentum equation.}
As a first step, we shall take the limit as $\tau \to 0$ of the
discrete momentum equation~\eqref{e:discrmom} and of the discrete
heat equation~\eqref{weak-heat-discr} with more regular  test
functions, which,  for technical reasons, we shall need to
approximate carefully. More precisely,
 for the momentum balance
equation~\eqref{e:weak-momentum} we shall use  test functions
\begin{subequations}
\begin{equation}
\label{testu} \testu \in
 L^{2}
(0,T;W_{\Gdir}^{1,2+\upsilon}(\Omega{{\setminus}} \GC;\R^d))  \cap
W^{1,1}(0,T; L^2(\Omega;\R^d)) \quad \text{for some $\upsilon>0$,}
\end{equation} and we shall approximate them with
  discrete approximations $\{
\testu_{\tau}^k \}$, such that the related piecewise constant and
linear interpolants fulfil as $\tau \to 0$,
\begin{equation}
\label{conve-testu} \left.
\begin{array}{ll}
\pwl \testu{\tau} \to \testu & \text{in $ W^{1,1}(0,T; L^2(\Omega;
\R^d))$,}   \\ \pwc \testu{\tau} \to \testu & \text{in $L^2 (0,T;
W^{1,2{+}\upsilon}(\Omega{{\setminus}} \GC;\R^d)$ for some
$\upsilon>0$,}
\\
\tau^{1/\gamma} e(\pwc \testu{\tau}) \to 0 & \text{in $
L^{\gamma}(Q;\R^{d\times d})$,}  \\
   \tau^{\beta/2}\| \pwc \testu{\tau} \|_{L^1 (0,T;
L^{4/(5{-}\mu)}(\GC))} \leq C\,. &
\end{array}
 \right\}
\end{equation}
\end{subequations}
 {Now,} combining \eqref{e:convutau1} and
\eqref{conv3}  with the second of \eqref{conve-testu}, and
\eqref{conv1} with the third of \eqref{conve-testu},
 we pass to the limit as $\tau \to 0$ in
the first integral term on the left-hand side of~\eqref{e:discrmom}.
Secondly,  \eqref{e:convztau1} and \eqref{conv2} yield
\[
 \dela{\pwc z {{\eps}\tau}}\JUMP{\pwc u {{\eps}\tau}}{} \weaksto  \dela
 \ze\JUMP{\ue}{}\qquad \text{in $L^\infty (0,T; L^{4-\epsilon}(\GC))$ for all $\epsilon \in [0,3)$,}
\]
which we combine with the second of~\eqref{conve-testu}. Also taking
into account \eqref{e:sweak} and \eqref{conv4}, together with the
fourth of \eqref{conve-testu},  we take the limit of the second
integral term on the left-hand side of~\eqref{e:discrmom}. In the
case $\varrho>0$, we take the limit of the third and fourth terms on
the left-hand side, and of the first term on the right-hand side of
\eqref{e:discrmom} by means of \eqref{e:convutau3} (combined with
the first of \eqref{conve-testu}),  of \eqref{pointiwise-for-u}, and
of \eqref{est-init-data}. Finally, using \eqref{data-converg-1} and
\eqref{data-converg-2}  we handle the second and third
right-hand-side terms
 in~\eqref{e:discrmom}. We thus conclude   that the
triple $(\ue,\ze,\we)$ fulfils equation~\eqref{e:weak-momentum},
first with test functions  as in~\eqref{testu} and ultimately, by a
density argument, with test functions $\testu \in L^2
(0,T;W_{\Gdir}^{1,2}(\Omega{{\setminus}} \GC;\R^d)) \cap
W^{1,1}(0,T;L^2
(\Omega;\R^d))$. 
\par\noindent
\textbf{Step $2$: passage to the limit in the semistability
condition.}
We consider  a subset $\mathcal{N} \subset (0,T)$ of full
measure such that  for all $t \in \mathcal{N}$  the approximate
stability condition \eqref{disc-semistability} holds independently
of  $\tau \to 0$. Then we fix $t \in \mathcal{N}$ and $\tilde{z} \in
L^\infty (\GC)$. We may suppose without loss of generality that
$\mathcal{R}(\tilde{z} - z_\eps(t))<+\infty$, hence
\begin{equation}
\label{e:use} \tilde{z}(x) \leq z_\eps(t,x)\qquad \foraa\,x \in \GC\,.
\end{equation}
We then construct the following recovery sequence
\begin{align}\label{eq5:recov-seq}
   \tilde z_{{\eps}\tau}(t,x):=\begin{cases}
   \displaystyle \pwc z{{\eps}\tau}(t,x)\frac{\tilde{z}(x)}{z_\eps(t,x)}
    & \text{where $z_\eps(t,x)>0$},
\\
   \displaystyle 0 & \text{where $z_\eps(t,x)=0$}.\end{cases}
\end{align}
Now, using~\eqref{e:use} and~\eqref{e:convztau2} one immediately
sees that
\begin{equation}
\label{e:use2} \tilde z_{{\eps}\tau}(\cdot,t) \leq \pwc
z{{\eps}\tau}(\cdot,t) \quad \aein \,\GC, \qquad \tilde
z_{{\eps}\tau}(t) \weaksto \tilde{z} \ \ \text{in $L^\infty (\GC)$.}
\end{equation}
Plugging $\tilde z_{{\eps}\tau}$ in~\eqref{disc-semistability}, we find
\begin{align}
\nonumber 0 & \leq \limsup_{\tau \to 0} \Big( \Phi_{{\eps}\tau}(\pwc
u{{\eps}\tau}(t), \tilde z_{{\eps}\tau}(t)){+} \mathcal{R}(\tilde
z_{{\eps}\tau}(t) - \pwc z{{\eps}\tau}(t) ){-}\Phi_{{\eps}\tau}(\pwc
u{{\eps}\tau}(t), \pwc z{{\eps}\tau}(t))\Big)
\\\nonumber&
\begin{aligned}=\limsup_{\tau \to 0}\int_{\GC}\!\!\Big(\frac{\tau^\alpha}{2}
\left(|\tilde z_{{\eps}\tau}(t)|^2  {-} | \pwc z{{\eps}\tau}(t)|^2
\right) & + \frac{\dela}{2} \big|\JUMP{\pwc
u{{\eps}\tau}(t)}{}\big|^2 \left( \tilde z_{{\eps}\tau}(t) -\pwc
z{{\eps}\tau}(t)\right)
-(a_0{+}a_1)\big(\tilde z_{{\eps}\tau}(t) {-} \pwc
z{{\eps}\tau}(t)\big) \Big)\, \d S
\end{aligned}
\\\nonumber&
\begin{aligned}
\leq \int_{\GC} \Big(  \frac{\dela}{2}
\big|\JUMP{u_{\eps}(t)}{}\big|^2 \left( \tilde {z}(t)
-z_{\eps}(t)\right)  {-}(a_0{+}a_1)\big(\tilde z_{{\eps}}(t) {-}
\pwc z{{\eps}}(t)\big)\Big)\, \d S
\end{aligned}
\\\label{eq6:delam-recovery}
& = \Phi_{\eps}(u_{\eps}(t), \tilde {z}(t)){+} \mathcal{R}(\tilde
{z}(t) - z_{\eps}(t) ){-}\Phi_{\eps}(u_{\eps}(t), {z}_\eps(t))\,,
\end{align}
where the
\COL{second inequality} ensues from~\eqref{e:use2}
and~\eqref{conv2}. 
\par\noindent
\textbf{Step $3$: passage to the limit in the  mechanical and total
energy inequalities.} Using \eqref{e:convutau1},
\eqref{pointiwise-for-u}, \eqref{e:finally}, and
\eqref{e:altogether},
 we pass to the limit in the left-hand side of the
discrete mechanical energy inequality \eqref{disc-energy0}  by weak
lower semicontinuity. To take the limit of the right-hand side, we
employ~\eqref{est-init-data}, the weak convergence
\eqref{e:convutau1} and the strong convergence~\eqref{conv3}, which
yield
\begin{align}
\label{conv-adiab} \Theta (\pwc \w{{\eps}\tau}) \bbB{:}
e(\pwl{\DT{u}}{{\eps}\tau}) \rightharpoonup \Theta (\w_{\eps}) \bbB
{:}e({\DT{u}}_{\eps}) \quad \text{weakly in $L^1 (Q)$.}
\end{align}
 Also using
\eqref{data-converg-1}--\eqref{data-converg-2},  we conclude that
the triple $(\ue,\ze,\we)$ complies  for all $t\in [0,T]$ with
\begin{align}\nonumber
&\! \! \! \! \! \! \! \!
T_\mathrm{kin}^{\varrho}\big(\DT{u}_\eps(t)\big)
+\Phi_{\eps}\big(\ue(t),\ze(t)\big) +\int_0^t \int_{\Omega}\!\!\!
\bbD e\big(\DT{u}_\eps(s)\big){:} e\big(\DT{u}_\eps(s)\big)\,\d
x\d s + \mathrm{Var}_{\mathcal{R}}(z_\eps;[0,t]) \le
T_\mathrm{kin}^{\varrho}\big(\DT{u}_{0})
\\
&\! \! \! \! \! \! \! \!  +\Phi_{\eps}\big(u_{0},z_{0})
+\int_0^{t}\bigg(\int_\Omega \Theta(\we(s))\bbB
{:}e\big(\DT{u}_\eps (s) \big) \,\d x + \int_{\Omega}\FRM(s)
{\cdot}\DT{u}_\eps(s)\, \d x   +\int_{\Gnew} \fRM(s) {\cdot}
\DT{u}_\eps(s) \, \d S \bigg) \d s\,. \label{disc-energy0-lim}
\end{align}

We also pass to the limit in the discrete total energy
inequality~\eqref{disc-energy}. Indeed, one tackles the left-hand
side by  the above-mentioned lower-semicontinuity arguments (also
using~\eqref{e:poinwtiwise-w}), and passes to the limit in the
right-hand side by convergences \eqref{est-init-data} and
\eqref{data-converg}--\eqref{data-converg-bis}. Thus, the total
energy inequality for the $\eps$-approximate problem holds for all
$t \in [0,T]$.

\par\noindent
 \textbf{Step $4$: mechanical energy equality.}
Like in~\cite{tr1}, we prove
 that, in the limit, the mechanical energy inequality~\eqref{disc-energy0-lim}  in fact holds as an
 equality, obtaining
\begin{align}\nonumber&
T_\mathrm{kin}^{\varrho}\big(\DT{u}_\eps(t)\big)
+\Phi_{\eps}\big(u_\eps(t),z_\eps(t)\big) +\int_0^{t}\int_{\Omega}
\bbD e\big(\DT{u}_{\eps}\big){:} e\big({\DT{u}}_{\eps}\big)\, \d
x \d t + \mathrm{Var}_{\mathcal{R}}(z_\eps;[0,t])
= T_\mathrm{kin}^{\varrho}\big(\DT{u}_{0})
+\Phi_{\eps}\big(u_{0},z_{0})
\\\label{mecheq} &\quad
+\int_0^{t}\bigg(\int_\Omega\Theta({\w}_{\eps})\bbB{:}
e\big(\DT{u}_{\eps}\big)\, \d x + \int_{\Omega} \FRM {\cdot}
{\DT{u}}_{\eps}\, \d x+\int_{\Gnew} \fRM(s){\cdot}\DT{u}_\eps(s)
\, \d S\bigg)\, \d s \quad \text{for all $t \in [0,T]$}.
\end{align}
To this aim, we develop the same calculations  as
 throughout~\cite[Formulae~(4.69)-(4.76)]{tr1}. The  first step  of
 the argument is
 a
sophisticated trick based on the  previously proved
semistability condition (see also, e.g.,  \cite{DMFraToa05,MieFra06}
for the use of such a technique in a rate-independent context),
which allows us to prove the following inequality for all $t \in
[0,T]$
\begin{equation}
\label{e:faith}
\begin{aligned}
\Phi_{\eps} & \big(u_\eps(t),z_\eps(t)\big)  -
\Phi_{\eps}\big(u_{0},z_{0})  +
\mathrm{Var}_{\mathcal{R}}(z_\eps;[0,t])  \geq \int_0^t
\pairing{}{}{(\Phi_\eps)_u'(\ue,\ze)}{\DT{u}_\eps}\,
 \d s
\\ &
= \int_0^t \int_{\Omega{\setminus}\GC} \bbC e(\ue) {:}
e(\DT{u}_\eps)\, \d x \d s  + \int_{0}^t \int_{\GC} \left(\delam \ze
\JUMP{\ue}{}{\cdot} \JUMP{\DT{u}_\eps}{} + \yosd(\JUMP{\ue}{})
{\cdot}\JUMP{\DT{u}_\eps}{}    \right)\, \d S\d s\,,
\end{aligned}
\end{equation}
where $(\Phi_\eps)_u'$ denotes the partial G\^ateaux-derivative with
respect to $u$ of the functional
$\Phi_\eps:W^{1,2}(\Omega{\setminus}\GC)\times L^\infty(\GC)\to\R$, and the
\COL{equality} follows from the definition~\eqref{8-1eps} of
$\Phi_\eps$. The second step consists of  testing of the momentum balance
equation~\eqref{e:weak-momentum} by $\DT{u}_\eps$. In the case $\varrho=0$,
$\DT{u}_\eps \in L^2 (0,T; W_{\Gdir}^{1,2}(\Omega{{\setminus}} \GC;\R^d))$ is an
admissible test function for~\eqref{e:weak-momentum}. In the case
$\varrho>0$, the test by $\DT{u}_\eps$ may be performed after proving that
$\DDT{u}_\eps \in L^2 (0,T;W_{\Gdir}^{1,2}(\Omega{{\setminus}} \GC;\R^d)^*)$,
cf.~Remark~\ref{rem:relax-test-delam}. In fact, a comparison
argument in~\eqref{e:weak-momentum} readily yields that
$\DDT{u}_\eps \in
L^2(0,T;W_{\Gdir}^{1,2}(\Omega{{\setminus}}\GC;\R^d)^*)$. Choosing
$\DT{u}_\eps$ as a test function in~\eqref{e:weak-momentum} and
integrating on $(0,t)$ for all $t \in [0,T]$ leads, after an
integration by parts, to
\begin{align}\nonumber
\frac{\varrho}2 & \int_\Omega |\DT{u}_\eps (t)|^2 \,\d x  + \int_0^t
\int_{\Omega}\bbD e(\DT{u}_\eps){:}  e(\DT{u}_\eps)\, \d x \d s
 +  \int_0^t \int_{\Omega} \bbC e(\ue) {:} e(\DT{u}_\eps)\, \d x
\d s
\\\nonumber &
+ \int_{0}^t \int_{\GC} \left(\delam \ze \JUMP{\ue}{} +
\yosd(\JUMP{\ue}{}) {\cdot}\JUMP{\DT{u}_\eps}{}    \right)\, \d S\d s
\\ &\label{e:moravia} =\frac\varrho2 \int_\Omega |\DT u_0|^2\, \d x
+\int_0^t\left( \int_\Omega \Theta({\w}_{\eps})\bbB{:}
e\big(\DT{u}_{\eps}\big)\,\d x + \int_\Omega\!\FRM{\cdot}\DT{u}_\eps
\, \d x +\int_{\Gnew}\!\!\fRM(s){\cdot}\DT{u}_\eps(s) \,
\d S\right) \, \d s\,.
\end{align}
Combining \eqref{e:faith} with \eqref{e:moravia}, we obtain the
reverse inequality in \eqref{disc-energy0-lim} and thus conclude
\eqref{mecheq}.
 \par\noindent
\textbf{Step $5$: passage to the limit in the enthalpy equation.}
First of all,  we observe that the following chain of inequalities
holds for all $t\in[0,T]$:
\begin{align}
\nonumber
 &\hspace*{-2em}\mathrm{Var}_{\mathcal{R}}(z_\eps;[0,t]) +
\int_0^{t}\int_{\Omega} \bbD e(\DT{u}_{\eps}) {:}
e(\DT{u}_{\eps})\, \d x \d t
\\& \nonumber
  \leq \liminf_{\tau \to 0} \int_0^t
\int_{\GC}\zeta_1 \left(\pwl{\DT{z}}{{\eps}\tau} \right)\, \d S \d t
+ \int_0^{t}\int_{\Omega} \bbD e(\pwl {\DT{u}}{{\eps}\tau}) {:}
e(\pwl {\DT{u}}{{\eps}\tau})\, \d x \d t
 \\& \nonumber
\le\limsup_{\tau \to 0}
T_\mathrm{kin}^{\varrho}(\DT{u}_{0,\tau})
+\Phi_{{\eps}\tau}(u_{0,\tau},z_{0}) -
T_\mathrm{kin}^{\varrho}(\pwl{\DT{u}}{{\eps}\tau}(t))
-\Phi_{{\eps}\tau}(\pwc
u{{\eps}\tau}(t), \pwc z{{\eps}\tau}(t))
\\\nonumber
&\qquad+\int_0^{t}\left(\int_\Omega
\Theta(\pwc{\w}{{\eps}\tau})
\bbB{:}e\big(\pwl{\DT{u}}{{\eps}\tau}\big) + \int_{\Omega}
\pwc \FRM{\tau} {\cdot} \pwl{\DT{u}}{{\eps}\tau}\, \d x
+\int_{\Gnew}\pwc \fRM{\tau} {\cdot} \pwl{\DT{u}}{{\eps}\tau} \, \d S
\right)\, \d t
\\ \nonumber
& \leq  T_\mathrm{kin}^{\varrho}(\DT{u}_{0}) +\Phi_{\eps}(u_{0},z_{0})
- T_\mathrm{kin}^{\varrho}({\DT{u}}_{\eps}(t)) -\Phi_{\eps}(
u_{\eps}(t), z_{\eps}(t))
+\int_0^{t}\bigg(
\int_\Omega\Theta({\w}_{\eps})\bbB{:} e\big({\DT{u}}_{\eps}\big)
+\FRM {\cdot} \DT{u}_\eps\,\d x
\\ &\qquad   + \int_{\Gnew}\!\!\fRM
{\cdot} \DT{u}_\eps \, \d S\bigg)\, \d t
\label{e:finally-bis}
 =
\mathrm{Var}_{\mathcal{R}}(z_\eps;[0,t]) + \int_0^{t}\int_{\Omega}
\bbD e(\DT{u}_{\eps}) {:} e(\DT{u}_{\eps})\, \d x \d t.
\end{align}
Indeed,  the first
\COL{inequality} ensues from~\eqref{e:convutau1}
and~\eqref{e:finally}, the second one from the discrete mechanical
energy inequality~\eqref{disc-energy0}, the third one from
\eqref{est-init-data}, \eqref{pointiwise-for-u},
\eqref{e:altogether}, \eqref{conv-adiab}, and from
\eqref{data-converg-1}--\eqref{data-converg-2}, cf.\ also Step 3.
Finally, the last equality ensues  from~\eqref{mecheq}. Thus, all of
the above inequalities turn out to hold with an equality sign. By a
standard $\liminf/\limsup$ argument,  this entails that
\begin{equation}
\label{e:stop0} \left.
\begin{array}{ll}
 \lim_{\tau \to 0}
T_\mathrm{kin}^{\varrho}(\pwl{\DT{u}}{{\eps}\tau}(t))=
T_\mathrm{kin}^{\varrho}({\DT{u}}_{\eps}(t)),
 \\
 \lim_{\tau \to 0}\Phi_{{\eps}\tau}(\pwc u{{\eps}\tau}(t), \pwc z{{\eps}\tau}(t))=\Phi_{\eps}(
u_{\eps}(t), z_{\eps}(t))\end{array} \right\}\quad \text{for all $t
\in [0,T]$,}
\end{equation}
 as well as
\begin{equation}
\label{e:stop2} \bbD e(\pwl {\DT{u}}{{\eps}\tau}) {:} e(\pwl
{\DT{u}}{{\eps}\tau}) \to \bbD e({\DT{u}}_{\eps}) {:}
e({\DT{u}}_{\eps}) \quad \text{strongly in $L^1(Q)$.}
\end{equation}
 Furthermore, arguing as in~\cite{tr1}, from
\eqref{e:finally-bis} holding as an equality we conclude that
\begin{equation}
 \label{e:stop} \zeta_1 \big(\pwl{\DT{z}}{{\eps}\tau}\big)
\weaksto\,\varmeaps \quad \text{in measure on $\overline{\SC}$,}
\end{equation}
with $\varmeaps $ being the measure introduced in~\eqref{meash}.

We are now in the position of taking the limit
of~\eqref{weak-heat-discr},  where we shall use test functions
\begin{equation}
\label{testw} \testw \in
\mathrm{C}^0(0,T;W^{1,r'+\varsigma}(\Omega{{\setminus}}\GC))\cap
W^{1,r'}(0,T;L^{r'}(\Omega))\qquad \text{for some $\varsigma>0$,}
\end{equation}
and we shall approximate them with discrete approximations
$\{ \testw_{\tau}^k \}$, such that, $\tau\to0$, the related interpolants fulfil
as $\pwc\testw{\tau} \to \testw  $ in $\mathrm{C}^0
(0,T; W^{1,r'+\varsigma}(\Omega{{\setminus}} \GC))$ for some
$\varsigma>0$, and $ \pwl \testw{\tau} \to \testw $ in
$W^{1,r'}(0,T; L^{r'}(\Omega))$.
 Then, we pass to the limit in the
first integral term on the left-hand side by exploiting
\eqref{e:poinwtiwise-w} and the aforementioned convergence for the
test functions $\pwc {\testw}{{\eps}\tau}$. To deal with the second
term we observe that, due to \eqref{e:convutau5}, to
\eqref{e:convwtau2}, and to the boundedness of the function
$\mathcal{K}: \R^{d\times d} \times \R \to \R^{d\times d}$, there
holds $ \mathcal{K}(e(\pwc
u{{\eps}\tau}),\pwc\w{{\eps}\tau})\to\mathcal{K}(e(\ue),\we)$ in
$L^q (Q)$ for all $1\leq q<\infty$, which we combine with the weak
convergence~\eqref{e:convwtau1} for $\pwc\w{{\eps}\tau}$ and with
the convergence for $\pwc {\testw}{{\eps}\tau}$. It follows
from~\eqref{eta-affine} and from convergences~\eqref{e:convztau1}
and~\eqref{conv2} that $ \eta({\JUMP{\underline{u}_{\eps\tau}}{},}
\pwc z{{\eps}\tau}) \weaksto \eta({\JUMP{u_\eps}{},}\ze)$ in
$L^\infty (0,T;L^{3{-}\epsilon}(\GC))$ for all $\epsilon \in (0,2]$,
which we exploit with \eqref{e:interesting} to take the limit of the
third integral term. The passage to the limit in the fourth term
results from \eqref{e:convwtau1} and the convergence for
$\pwl\w{{\eps}\tau}$. As for the right-hand side of
\eqref{weak-heat-discr}, to deal with the first integral term  we
exploit \eqref{e:stop2} and the convergence for $\pwc
{\testw}{{\eps}\tau}$, which in particular yields $ \pwc
{\testw}{{\eps}\tau} \to \testw $ in $\mathrm{C}^0 (\overline{Q})$.
Relying on this convergence and on \eqref{e:stop}, we also infer
\[
\lim_{\tau \to 0}
\bigg({-}\int_{\SC}\!\!a_1\pwl{\DT{z}}{{\eps}\tau}\frac{\pwc
\testw{\tau}|_{\GC}^+{+} \pwc\testw{\tau}|_{\GC}^-}2 \, \d S\d
t\bigg)=\int_{\SC}\!\!\frac{ \testw|_{\GC}^+{+}
\testw|_{\GC}^-}2\,\varmeaps (\d S\d t).
\]
Finally, employing \eqref{data-converg-bis}, one takes the limit of
the last three terms on the right-hand side of
\eqref{weak-heat-discr}, thus finding that the triple
$(\ue,\ze,\we)$ fulfils the weak formulation~\eqref{weak-heat} of
the enthalpy equation for all test functions as in~\eqref{testw}.
Again by a density argument, we conclude that $(\ue,\ze,\we)$
fulfil~\eqref{weak-heat} with test functions $ \testw \in
\mathrm{C}^0 ([0,T];W^{1,r'}(\Omega{{\setminus}} \GC)) \cap
W^{1,r'}(0,T; L^{r'}(\Omega))$.

This concludes the proof that solves $(\ue,\ze,\we)$
the approximate problem, i.e.\
Theorem~\ref{th:4.1}. $\hfill\Box$
\begin{remark}[Strong convergence]
\upshape
 Let us observe that, in the case $\varrho>0$,
 by a classical argument based on a Korn-type inequality and
  the uniform convexity of the space
$W^{1,2}(0,T;W_{\Gdir}^{1,2}(\Omega{{\setminus}} \GC;\R^d))$,
convergence~\eqref{e:stop2} joint with~\eqref{pointiwise-for-u}
allows us to conclude that $ \pwl {u}{{\eps}\tau} \to u_\eps$ in
$W^{1,2}(0,T;W_{\Gdir}^{1,2}(\Omega{{\setminus}}\GC;\R^d))$.
Likewise, if $\varrho>0$ \eqref{e:stop0} gives
$\pwl{\DT{u}}{{\eps}\tau}(t)\to{\DT{u}}_{\eps}(t)$ in
$L^2(\Omega;\R^d)$ for all $t \in [0,T]$.
\end{remark}

\begin{remark}[Numerics]
\upshape We point out that our method of proof may yield some
strategy for numerical analysis after making  a spatial
discretization, although the non-variational structure of \eq{GMa}
(preserved if discretized in space) would still require some
iterative procedure for numerical solution, cf.
\cite{bartels-roubicek}.
\end{remark}

\section{Limit passage with $\eps\to0$ and proof of Theorem~\ref{th:3.0}}\label{s:5}
\noindent In passing to the limit in the $\eps$-approximate problem
as $\eps \to 0$,  we shall follow the steps of the proof of
Theorem~\ref{th:4.1}. Thus, we shall sketch most of the arguments,
referring to the detailed calculations developed in
Section~\ref{ss:4.4}, and dwell with some detail only on the
passages to the limit as $\eps \to 0$ in the momentum equation, and
on the proof of the mechanical energy equality.
\par\noindent
\textbf{Step $0$: a-priori estimates and compactness.} The sequence
$(\ue,\ze,\we)_\eps$ inherits the a priori estimates of
Lemma~\ref{prop:apriori}, i.e.~now
\begin{subequations}\label{a-priori4}
\begin{align}
& \label{a41} \! \! \! \! \! \!  \! \! \big\|\ue \big\|_{
W^{1,2}(0,T;W_{\Gdir}^{1,2}(\Omega;\R^d))} + \varrho^{1/2} \big\|\ue
\big\|_{ W^{1,\infty}(0,T;L^2(\Omega;\R^d)) } +  \varrho \|
\DDT{u}_\eps\|_{L^2 (0,T;
 W_{\GC}^{1,2}(\Omega;\R^d)^*)} \le S\,,
 \\
\label{a42}  & \! \! \! \! \! \!  \! \!
\big\|\ze\big\|_{L^\infty(\SC)} +\big\|\ze \big\|_{\BV([0,T];L^1(\GC))} \leq
S\,,
\\
& \label{a43} \! \! \! \! \! \!  \! \!
\big\|\we\big\|_{L^\infty(0,T;L^1(\Omega))} + \big\|\we
\big\|_{L^r(0,T;W^{1,r}(\Omega))}+\big\|\we\big\|_{\BV([0,T];W^{1,r'}(\Omega)^*)}
\le S_r' \ \mbox{  for any $1\le r<\frac{d+2}{d+1}$\,,}
 \\
& \label{bound-energies} \! \! \! \! \! \!  \! \!  \sup_{t\in [0,T]}
\Phi_{\eps} \left(\ue(t), \ze(t) \right) \leq S.
\end{align}
\end{subequations}
for some $S>0$ and $S_r'>0$ depending on $1\le r<\frac{d+2}{d+1}$.
Indeed, the first two estimates
in~\eqref{a41},
 the first of~\eqref{a42}, and the second of~\eqref{a43} respectively follow
from \eqref{a31}, \eqref{a31bis}, \eqref{a32}, and \eqref{a33bis}
via lower-semicontinuity. The second of~\eqref{a42}, and \eqref{a43}
have been proved throughout Section~\ref{ss:4.4}, cf.~with
\eqref{bvl1}, \eqref{bvw1}, and \eqref{useful-later}. The third of
estimates~\eqref{a41} follows by testing~\eqref{e:weak-momentum} by
functions $\testu\in L^2(0,T;W_{\GC}^{1,2}(\Omega{{\setminus}}\GC;\R^d))
\cap W^{1,1}(0,T; L^2(\Omega;\R^d))$ and taking into account \eqref{a41},
the second of~\eqref{a43}. Finally, \eqref{bound-energies} is
a direct consequence of the total energy inequality.

By the Banach and the Helly selection principles, there is a
subsequence (for simplicity, denoted by the same indexes) and
$(\ude,\zde,\wde)$ such that the following
convergences hold: \nopagebreak
\begin{subequations}
\label{convej}
\begin{align}
 \label{e:convej1} & \uej \weakto
\ude  \ \  \text{ in $
W^{1,2}(0,T;W_{\Gdir}^{1,2}(\Omega{{\setminus}} \GC;\R^d)),$}
\\
 \label{e:convej2}
 & \uej \to
\ude  \ \ \text{ in $
\mathrm{C}^{0}([0,T];W_{\Gdir}^{1,2-\epsilon}(\Omega{{\setminus}}
\GC;\R^d))$} \ \ \forall\, \epsilon \in
 (0,1],
\\
\label{e:convej3} &  \text{if $\varrho>0$,} \  \ \uej{\weaksto}\ude
\text{ in $ W^{1,\infty}(0,T;L^2(\Omega{{\setminus}} \GC;\R^d))
\cap W^{2,2}(0,T;
 W_{\GC}^{1,2}(\Omega{{\setminus}} \GC;\R^d)^*)$,}
\\
& \label{e:convej4}
\begin{aligned}
\text{if $\varrho>0$,} \  \ \uej{\to}\ude   & \text{
in $W^{1,2}(0,T;W_{\Gdir}^{1,2-\epsilon}(\Omega{{\setminus}}
\GC;\R^d))
 \cap W^{1,q}(0,T;L^2(\Omega;\R^d))$}
   \ \forall\, 0{<}\epsilon{\le}1, \,
1{\le}q{<}\infty,\!\!\!
 \end{aligned}
\\
& \label{e:convej5} \text{if $\varrho>0$,} \  \ \DT{u}_{\eps} (t)
\weakto \DT{u}(t) \quad \text{in $L^2(\Omega;\R^d)$ for all $t \in
[0,T]$,}
\\
& \label{e:convej6}
 \zej \weaksto \zde  \ \ \text{ in $L^{\infty}(\SC)$,}
\\
& \label{e:convej7}
 \zej(t) \weaksto \zde(t) \ \  \text{ in $L^\infty(\GC)$ for
all $t \in [0,T]$,}
\\
 &
\label{e:convej8}
  \wej
\rightharpoonup \wde \ \ \text{ in $
L^{r}(0,T;W^{1,r}(\Omega{{\setminus}} \GC))$,}
\\
 \label{e:convej9} &
\wej  \to \wde \ \ \text{ in $
L^{r}(0,T;W^{1,r-\epsilon}(\Omega{{\setminus}} \GC)) \cap L^q (0,T;
L^1 (\Omega))$}
 \ \ \forall\, \epsilon \in
 (0,r{-}1], \, 1{\le}q{<}\infty,
\\
 \label{e:convej10} &
 \wej(t)
\weakto \wde (t) \quad \text{in $W^{1,r'}(\Omega{{\setminus}}\GC)^*$
for all $t \in [0,T]$,}
\\
 \label{e:convej11} &
\!\!\!\!\left.\begin{array}{ll} \JUMP{\uej}{} \to \JUMP{\ude}{} &
\text{in $L^\infty (0,T; L^{4-\epsilon} (\GC;\R^d))$,} \\[.3em]
\JUMP{\uej(t)}{} \to \JUMP{\ude(t)}{}
 &  \text{in $L^{4-\epsilon} (\GC;\R^d)$ for any $t \in [0,T]$,}
 \end{array}
 \right\}\text{ for all $\epsilon\in(0,3]$}
\\
 \label{e:convej13} &
\Theta(\wej) \to \Theta (\wde) \text{ in $L^2(Q)$,}
\\
\label{e:convej14} &
 \JUMP{\Theta(\wej)}{}  \to  \JUMP{\Theta(\wde)}{}  \text{ in $
L^{r}(0,T;L^{3/2}(\GC))$\,.}
\end{align}
\end{subequations}
Convergences \eqref{convej} can be deduced from
estimates~\eqref{a-priori4} arguing in the very same way as
throughout~\eqref{e:convutau1}--\eqref{e:interesting} in
Section~\ref{ss:4.4}.
\par\noindent
\textbf{Step $1$: passage to the limit in the momentum equation.}
First of all, notice that~\eqref{bound-energies},
\eqref{e:convej11}, and \eqref{liminf-later} yield
$S\ge\liminf_{\eps_j \to 0} \yosappej\big(\JUMP{\uej(t)}{}\big)
\ge\ind_{\cone}\big(\JUMP{\ude(t)}{}\big)$ for all $t\in [0,T]$,
whence~\eqref{constraints-delam}. We now exploit (\ref{convej}a-f,k-l) to
pass to the limit in \eqref{e:weak-momentum} with $\varrho\ge0$, tested by
$\testu-\ue$, for any $\testu$ in
$L^2(0,T;W_{\Gdir}^{1,2}(\Omega{\setminus}\GC;\R^d))$ (and in
$W^{1,1}(0,T;L^2(\Omega;\R^d))$ if $\varrho>0$), with $\JUMP{\testu}{}\GE0$
on~$\SC$, i.e. for any test function in~\eqref{e:weak-momentum-variational}.
Using the mentioned convergences, we find
\begin{align}\nonumber
 & \!  \! \! \!\limsup_{\eps\to 0} \int_{\SC}\yosdej\big(
\JUMP{\uej}{}\big){\cdot} \JUMP{\uej{-}v}{}\,\d S\d t
=\limsup_{\eps\to0}\int_\Omega\varrho\DT u_0{\cdot}
\big(\uej(0){-}\testu(0)\big)
-\varrho\DT u_\eps(T){\cdot} (\uej(T) {-}\testu(T))\,\d x
\\\nonumber
&\qquad
+\int_{\Snew}\!\!\fRM{\cdot}(\uej{-}\testu)\,\d S\d t
-\int_{\SC}\!\!
\dela z_\eps\JUMP{u_\eps}{}{\cdot}\JUMP{\uej{-}\testu}{} \, \d S\d t
\\\nonumber &\qquad
 + \int_Q\!\FRM{\cdot}(\uej{-}v)
-\big(\bbD e(\DT{u}_{\eps_j}){+}\bbC
e(\uej){-}\bbB\Theta(\w_{\eps_j})\big){:}e(\uej{-}v)
+\varrho\DT{\uej}{\cdot}(\DT{\uej}{-}v) \,\d x\d t
\\\nonumber&\le\int_\Omega\varrho\DT u_0{\cdot} (u_0{-}\testu(0))
-\varrho\DT u(T){\cdot}(u(T){-}\testu(T))
\,\d x
+\int_{\Snew}\!\!\fRM{\cdot}(u{-}\testu)\,\d S\d t -\int_{\SC}\!\!
\dela z\JUMP{u}{}{\cdot}\JUMP{u{-}\testu}{} \, \d S\d t
\\ &\qquad
+\int_Q\!\FRM{\cdot}(u{-}\testu)-\big(\bbD
e(\DT{u}){+}\bbC e(u){-}\bbB\Theta(\w)\big){:}e(u{-}v)
+\varrho\DT{u}{\cdot}(\DT{u}{-}v) \,\d x\d t.
\label{e:step1}
\end{align}
On the other hand, recalling formula~\eqref{yos-repre} for   the
Yosida regularization $\yosd$, we see that
\begin{align}\nonumber
&\liminf_{\eps\to0}\int_{\SC}\!\!\!\yosdej\big(\JUMP{\uej}{} \big)
{\cdot} \big( \JUMP{\uej}{} {-} \JUMP{\testu}{} \big)\d S \d t
\geq \liminf_{\eps \to 0} \frac1{\eps} \int_{\SC}\!\!\!
\big(\JUMP{\uej}{}{-}\mathrm{P}_{\cone}\big(\JUMP{\uej}{} \big)
\big) {\cdot} \big( \mathrm{P}_{\cone}\big(\JUMP{\uej}{}  {-}
\JUMP{\testu}{}
 \big)\big)\d S \d t
\\ &\hspace{9em}
+\liminf_{\eps\to 0}\frac1{\eps}\int_{\SC}\!\!\!\left(\JUMP{\uej}{}
- \mathrm{P}_{\cone}\big(\JUMP{\uej}{}\big) \right) {\cdot}
\big(\JUMP{\uej}{} {-} \mathrm{P}_{\cone}\big(\JUMP{\uej}{}\big)
\big) \, \d S \d t\ge0 \label{projection-inequality}
\end{align}
the latter inequality holding due to the properties of the
projection operator and the fact that $\testu\GE0$ on $\SC$.
Combining \eqref{e:step1} and \eqref{projection-inequality}, and
rearranging some terms, we readily conclude the weak
formulation~\eqref{e:weak-momentum-variational} of the momentum
inclusion.
\par\noindent
\textbf{Step $2$: passage to the limit in the semistability
condition.} It can be performed  by the  very same   recovery
sequence trick devised in Step $2$ of the proof of
Theorem~\ref{th:4.1}.
\par\noindent
\textbf{Step $3$: passage to the limit in the  mechanical and total
energy inequalities.} It follows from \eqref{e:convej1},
\eqref{e:convej6},
 \eqref{e:convej11}, and
\eqref{liminf-later} that
\begin{equation}
\label{e:altogether-delam} \Phi(\ude(t), \zde(t)) \leq
\liminf_{\eps_j \to 0} \Phi_{\eps_j}(\uej(t), \zej(t)) \quad
\text{for all $t \in [0,T]$.}
\end{equation}
Combining~\eqref{e:altogether-delam} with convergences
\eqref{convej} and arguing exactly  like in Step $3$  of the proof
of Theorem~\ref{th:4.1}, we pass to the limit by
lower-semicontinuity in conclude that $(\ude,\zde,\wde)$ complies
for all $t \in [0,T]$ with the mechanical energy
inequality~\eqref{disc-energy0-lim}, with $\Phi$ in place of
$\Phi_\eps$. Likewise, we conclude the total energy
inequality~\eqref{total-energy-brittle}, with $\Phi$ in place of
$\Phi_\eps$.
\par\noindent
 \textbf{Step $4$: mechanical energy equality.}
Arguing like  in Section~\ref{ss:4.4}
(cf.~\cite[Formulae~(4.69)-(4.76)]{tr1}), we first of all prove that
for all $t \in [0,T]$
\begin{align}\nonumber
\Phi & \big(\ude(t),\zde(t)\big)  - \Phi\big(u_{0},z_{0}) +
\mathrm{Var}_{\mathcal{R}}(z;[0,t]) \geq \int_0^t
\pairing{}{}{\lambda}{\DT{u}}
 \, \d s
 \\\label{e:faith-delam}
 & \text{for any
$\lambda\in L^2(0,T;W^{1,2}(\Omega{\setminus}\GC;\R^d)^*)$ with $\lambda(t)\in
 \partial_u \Phi (u(t),z(t))$ for a.a. $t \in (0,T)$,}
\end{align}
where $\partial_u \Phi: W^{1,2}(\Omega{{\setminus}} \GC;\R^d)
\rightrightarrows W^{1,2}(\Omega{{\setminus}} \GC;\R^d)^*$ denotes
the subdifferential w.r.t. $u$ of the functional $\Phi:
W^{1,2}(\Omega{{\setminus}} \GC;\R^d)\times L^\infty (\GC) \to \R $.
It follows from definition~\eqref{8-1-k} that the operator
$\partial_u \Phi$ is given by
\begin{align}\nonumber
\lambda  \in \partial_u \Phi(u,z)
 \ \ \text{if and only if} \ \ \exists\, \ell \in
 \indabs(u);\ \
\forall\, v \in W^{1,2}(\Omega{\setminus} \GC;\R^d):
\\\label{subdif-repre}
\pairing{}{}{\lambda}{v} =\int_{\Omega{\setminus}\GC}\!\!\!\bbC e(u)
{:} e(v)\, \d x + \int_{\GC}\!\!\delam z \JUMP{u}{}  {\cdot}
\JUMP{v}{} \, \d S + \pairing{}{}{\ell}{v},
\end{align}
where we have introduced for  notational convenience the functional
$\abs: W^{1,2}(\Omega{\setminus}\GC;\R^d) \to [0,+\infty]$  defined
by $ \abs(u)= \ind_{\cone}(\JUMP{u}{})$ its subdifferential
$\indabs: W^{1,2}(\Omega{\setminus}\GC;\R^d) \rightrightarrows
W^{1,2}(\Omega{\setminus}\GC;\R^d)^* $. Notice that $\indabs =
\jum^* \circ \partial I_{\cone} \circ \jum$, \COLLL{where $J$
denotes} the jump operator $\jum(u)=\JUMP{u}{}$ and $\jum^*$ its
adjoint. Now, let us observe that
\begin{equation}
\label{chain-rule}
\begin{aligned}
& \int_{0}^t\!\!\big\langle\ell,\DT{u}\big\rangle\,\d s
=\abs(\ude(t))-\abs(\ude(0))=  \ind_{\cone}\big( \JUMP{\ude(t)}{}\big) -
\ind_{\cone}\big( \JUMP{\ude(0)}{}\big)=0
\\
 &\qquad\text{for all $\ell \in L^2(0,T;W^{1,2}(\Omega{\setminus}\GC;\R^d)^*)$
such that $\ell(s)\in\partial
\COLL{\abs}\big(\JUMP{\ude(s)}{}\big) \ \foraa\, s \in (0,T)$,}
\end{aligned}
\end{equation}
by the chain rule for the convex functional  $\abs$
(cf.~\cite[Prop.~XI.4.11]{visintin96}), and by~\eqref{uzero}
and~\eqref{constraints-delam}. Combining
\eqref{e:faith-delam}--\eqref{chain-rule}, we conclude the following
inequality for all $t\in[0,T]$
\begin{align}
\label{e:other-crucial-inequality}
\Phi\big(u(t),z(t)\big) & - \Phi\big(u_{0},z_{0}) +
\mathrm{Var}_{\mathcal{R}}(z;[0,t])
\ge\int_0^t\!\bigg(\int_{\Omega} \bbC e(\ude){:}e(\DT{u})\, \d
x \d s +
\int_{\GC}\!\!\!\delam\zde\JUMP{\ude}{}{\cdot}\JUMP{\DT{u}}{}\,\d
S\bigg)\d s\,.
\end{align}
We now distinguish the two cases $\varrho=0$ and~$\varrho>0$.
\\ \textbf{Case $\varrho>0$:} arguing by
comparison in~\eqref{e:weak-momentum-variational} we may readily
check that $\DDT{u}\in L^2
(0,T;W_{\cone}^{1,2}(\Omega{{\setminus}} \GC;\R^d)^*)$,
which is in duality with
$\DT{u}\in L^2(0,T; W_{\cone}^{1,2}(\Omega{{\setminus}}
\GC;\R^d)).$ As pointed out in Remark~\ref{rem:relax-test-delam},
this entails that $\DT{u}$ is an admissible test function for the
momentum balance inclusion~\eqref{e:weak-momentum-variational}.
Then, upon proceeding with such a test and again
using~\eqref{chain-rule} we conclude for all $t \in [0,T]$ that
\begin{align}\nonumber
\frac{\varrho}2\int_\Omega |\DT{u} (t)|^2 \,\d x  &+ \int_0^t\!\!
\int_{\Omega}\bbD e(\DT{u}){:} e(\DT{u})\, \d x \d s
+  \int_0^t\!\!\int_{\Omega} \bbC e(\ude) {:}
e(\DT{u})\, \d x \d s + \int_{0}^t\!\!\int_{\GC} \delam \zde
\JUMP{\ude}{}{\cdot} \JUMP{\DT{u}}{}\, \d S\d s
\\\label{e:moravia-delam}
& =\frac\varrho2 \int_\Omega |\DT u_0|^2\, \d x + \int_0^t\!\!\left(\,
 \int_\Omega \Theta(\wde)\bbB{:}
e\big(\DT{u}\big)\, \d x  + \int_\Omega\!\FRM{\cdot}\DT{u} \, \d x +
\int_{\Gnew}\!\!\!\fRM{\cdot}\DT{u}\, \d S \right) \, \d s\,.
\end{align}
Combining~\eqref{e:moravia-delam}
with~\eqref{e:other-crucial-inequality},   we get the reverse of the
mechanical energy inequality, which leads to the  desired mechanical
energy equality. \\
\textbf{Case $\varrho=0$:}
From \eqref{e:weak-momentum-variational} one infers that the functional
\begin{align}\label{def-of-ell}
\ell:v\mapsto \int_Q\!\big(\bbD e(\DT{u}){+}\bbC
e(u){-}\bbB\Theta(\w)\big) {:}e(v)
\,\d x\d t +\int_{\SC}\!\!\!
\dela z\JUMP{u}{}
{\cdot}\JUMP{v}{} \, \d S\d t
-\int_{Q}\!\FRM{\cdot}v\, \d x\d t-
\int_{\Snew}\!\!\!\fRM{\cdot}v\,\d S\d t.
\end{align}
is in $L^2(0,T;W^{1,2}(\Omega{\setminus}\GC;\R^d)^*)$, and fulfils
\begin{align}
&
\int_{\SC}\!\!
I_{\cone}\big(\JUMP{v}{}\big)\d S\d t
\label{use-of-ell}
\ge
\int_{\SC}\!\!
I_{\cone}\big(\JUMP{u}{}\big)\,\d S\d t + \int_0^T
\big\langle\ell,v{-}u\big\rangle \, \d t.
\end{align}
Hence, $\ell(t)\in\partial{\abs}\big(\JUMP{\ude(t)}{}\big)$ for
almost all $t \in (0,T)$. Thus, \eqref{chain-rule} yields $\int_0^t
\pairing{}{}{\ell}{\DT{u}}\, \d s =0$ for all $t \in [0,T]$, which
is just relation~\eqref{e:moravia-delam} with $\varrho=0$.  Again,
we combine the latter with~\eqref{e:other-crucial-inequality}, and
conclude the mechanical energy equality.
\par\noindent
\textbf{Step $5$: passage to the limit in the enthalpy equation.} It
can be developed  in the very same way as in the proof  of
Theorem~\ref{th:4.1}. We point out that, if~\eqref{strict-pos}
holds, convergence~\eqref{e:convej9} and~\eqref{poswe} yield for
almost all $(t,x)\in Q$~the strict positivity of
$\theta(t,x)=\Theta(\w(t,x))$. $\hfill\Box$ 
\begin{remark}[Convergence of the reaction force]
\upshape \COLLL{Notice that,  if $\varrho=0$,} there exists $S'>0$
such that, for all $\eps>0$,
\begin{equation}\nonumber
\Big\|J^*\!\circ\yosd(\JUMP{u_\eps}{})\Big\|_{L^2(0,T; W^{1,2}(\Omega{\setminus}\GC;\R^d)^*)}
=\sup_{\testu \in L^2(0,T; W^{1,2}(\Omega{\setminus}\GC;\R^d))}
\bigg|\int_{\SC}\!\!\yosd\big(\JUMP{u_\eps}{}\big){\cdot}\JUMP{\testu}{}\,\d S\d
t\bigg|\le S',
\end{equation}
which
ensues from a comparison in~\eqref{e:weak-momentum} by using
\eqref{a-priori4} (in spite of the blow-up of the reaction force
$\yosd(\JUMP{u_\eps}{})$ in $L^\infty(0,T;L^2(\GC;\R^d))$,
cf.~\eqref{e:sweak}).
Passing to the limit in \eqref{e:weak-momentum}, we can see that,
for the subsequence selected in Step~0,
\COLLL{$J^*\!\circ\yosd(\JUMP{u_\eps}{})$ weakly   converges in
$L^2(0,T; W^{1,2}(\Omega{\setminus}\GC;\R^d)^*)$  to  the function
$\ell$} defined in \eqref{def-of-ell} .
\end{remark}

\bigskip

\baselineskip=10pt {\small \noindent{\it Acknowledgements}: This
research was initiated during a visit of R.R. at the
Charles University in Prague,
supported by the ``Ne\v cas
center for mathematical modeling'' LC 06052 (M\v SMT \v CR), and
partially by grants from PRIN 2008 project ``Optimal mass
transportation, geometric and functional inequalities and
applications". T.R.~acknowledges the
hospitality of the University of Brescia, as well as partial support
from the grants A~100750802 (GA~AV~\v CR),
201/09/0917 and 201/10/0357 (GA \v CR), and MSM~21620839 (M\v SMT \v
CR), and from the research plan AV0Z20760514 ``Complex dynamical systems
in thermodynamics, mechanics of fluids and solids'' (\v CR).}


\begin{thebibliography}{00}

\vspace{-.3em}\bibitem{AKRS02ODTC}
K.T.Andrews, K.L.Kuttler, M.Rochdi, M.Shillor: One-dimensional
dynamic thermoviscoelastic contact with damage. \emph{J.\ Math.\
Anal.\ Appl.} {\bf 272} (2002), 249--275.

\vspace{-.3em}\bibitem{attouch}
H.~Attouch: \newblock
      \emph{Variational Convergence for Functions and Operators}.
      Pitman,
      Boston, 1984.


\vspace{-.3em}\bibitem{barbu76}
V.~Barbu:
\newblock {\em Nonlinear Semigroups and Differential Equations in Banach
  Spaces}.
\newblock Noordhoff, Leyden, 1976.


\vspace{-.3em}\bibitem{BarbuPrecupanu86} V. Barbu, T. Precupanu:
\newblock  \emph{Convexity and Optimization in Banach
Spaces}. \newblock  D. Reidel Pub. Co., Dordrecht, 2nd edition,
1986.

\vspace{-.3em}\bibitem{bartels-roubicek} S. Bartels, T. Roub\'\i\v cek:
Thermo-visco-elasticity with rate-independent plasticity in
isotropic materials undergoing thermal expansion.
(Preprint no.486, SFB 611, Univ. Bonn, 2009) Math. Modelling Anal. Numer.
(submitted).

\vspace{-.3em}\bibitem{bedford}
A. Bedford: {\it Hamilton's Principle in Continuum Mechanics.} Pitman,
Boston, 1985.

\vspace{-.3em}\bibitem{boccardo-gallouet1}  L. Boccardo,  T. Gallou\"et: \newblock
Non-linear elliptic and parabolic equations involving measure data.
\newblock {\it J. Funct. Anal.} {\bf 87} (1989), 149--169.

\vspace{-.3em}\bibitem{BBR1}
E.Bonetti, G.Bonfanti,  R.Rossi:  \newblock  Global existence
for a contact problem with adhesion. \newblock  {\it Math. Methods
Appl. Sci.} {\bf 31} (2008), 1029--1064.

\vspace{-.3em}\bibitem{BBR2}
E.Bonetti, G.Bonfanti,  R.Rossi: \newblock   Well-posedness and
long-time behaviour for a model of contact with adhesion. \newblock
{\it Indiana Univ. Math. J.} {\bf 56} (2007), 2787–-2820.

\vspace{-.3em}\bibitem{BBR3} E. Bonetti, G. Bonfanti, R. Rossi: \newblock
  Thermal effects in adhesive contact: modelling and
analysis. \newblock  \emph{Nonlinearity} \textbf{22} (2009),
2697--2731.

\vspace{-.3em}\bibitem{BBR4} E. Bonetti, G. Bonfanti, R. Rossi: \newblock
Long-time behaviour of a thermomechanical model for adhesive
contact. \newblock (Preprint arXiv:0909.2493) \newblock
\emph{Discrete Contin. Dyn. Syst. Ser. S}, in print (2010).

\vspace{-.3em}\bibitem{brezis73}
H.~Br\'ezis: \newblock {\em Op\'erateurs Maximaux Monotones et
Semi-groupes de Contractions dans les Espaces de Hilbert}.
\newblock
North-Holland, Amsterdam,  1973.

\vspace{-.3em}\bibitem{ChGiPo08CIBM}
A. Chambolle, A. Giacomini,  M. Ponsiglione: \newblock   Crack
initiation in brittle materials. \newblock  {\it Arch. Rational Mech.
Anal.} {\bf 188} (2008), 309--349.

\vspace{-.3em}\bibitem{CooGor64MCCP}
J. Cook,  J.~E. Gordon, C.~C. Evans, D.~M. Marsh: \newblock  A
mechanism for the control of crack propagation in all-brittle
systems. \newblock   {\it Proc. R. Soc. Lond. Ser.A} {\bf 282}
(1964), 508--520.

\vspace{-.3em}\bibitem{DMFraToa05} G. Dal Maso, G. Francfort,  R. Toader: \newblock
Quasistatic crack growth in nonlinear elasticity. \newblock  {\it
Arch. Rational Mech. Anal.} {\bf 176} (2005), 165--225.

\vspace{-.3em}\bibitem{eck0}
C. Eck: \newblock  Existence of solutions to a thermo-viscoelastic
contact problem with Coulomb friction. {\it Math. Models Methods
Appl. Sci.} {\bf 12} (2002), 1491--1511.

\vspace{-.3em}\bibitem{eck}
C. Eck, J. Jaru\v sek,  M. Krbec: \newblock  {\it Unilateral
Contact Problems.}  \newblock  Chapman \& Hall/CRC, Boca Raton,
2005.

\vspace{-.3em}\bibitem{feireisl-malek}  E. Feireisl, J. M\'alek: \newblock  On the
Navier-Stokes equations with temperature-dependent transport
coefficients.  \newblock  {\it Diff. Equations Nonlin. Mech.}
(2006), 14pp. (electronic), Art.ID 90616.

\vspace{-.3em}\bibitem{fpr09} E. Feireisl, H. Petzeltov\'a, E. Rocca:
\newblock Existence of solutions to a phase transition model with microscopic
    movements. \newblock
    {\it Math. Methods Appl. Sci.} \textbf{32} (2009),  1345--1369.

\vspace{-.3em}\bibitem{fremlin3} D.~H. Fremlin:
\newblock
 {\it Measure theory. Vol. 3. Measure algebras.} \newblock
  Corrected second printing of the 2002 original. \newblock Torres Fremlin, Colchester, 2004.

\vspace{-.3em}\bibitem{Fre82} M. Fr\'emond: \newblock Equilibre
des structures qui adh\`erent \`a leur support. \newblock {\it
Comptes Rendus
Acad. Sci. Paris}, \textbf{295} (1982), 913--916.

\vspace{-.3em}\bibitem{Fre87} M. Fr\'emond: \newblock Adh\`erence des solides. \newblock {\it
Journal de M\'echanique Th\'eorique et Appliqu\'ee}, \textbf{6}
(1987), 383--407.


\vspace{-.3em}\bibitem{Grif20PRFS}
A. A. Griffith:  \newblock  The phenomena of rupture and flow in
solids.
\newblock   {\it Philos. Trans. Royal Soc. London Ser. A. Math. Phys.
Eng. Sci.} {\bf 221} (1921), 163--198.

\vspace{-.3em}\bibitem{MieFra06} G. Francfort,  A. Mielke: \newblock
Existence results for a class of rate-independent material
                  models with nonconvex elastic energies. \newblock
{\it J. reine angew. Math.} {\bf 595} (2006), 55--91.


\vspace{-.3em}\bibitem{GiaPon06GCAS}
A. Giacomini, M. Ponsiglione: \newblock   A
{${\Gamma}$}-convergence approach to stability of unilateral
minimality properties in fracture mechanics and applications.
\newblock   {\it Arch. Rational Mech. Anal.} {\bf 180} (2006),
399--447.

\vspace{-.3em}\bibitem{KoMiRo}
M. Ko\v cvara, A. Mielke,  T. Roub\'\i\v cek:  \newblock  A
rate-independent approach to the delamination problem. \newblock
 {\it Mathematics and Mechanics of Solids} {\bf 11} (2006), 423--447.

\vspace{-.3em}\bibitem{Legu02STCC}
D. Leguillon: Strength or toughness? {A} criterion for crack onset
at a notch.  \newblock  {\it Euro. J. Mechanics A/Solids} {\bf 21}
(2002), 61--72.

\vspace{-.3em}\bibitem{Mant08ICOC}
V. Manti\v{c}: \newblock  Interface crack onset at a circular
cylindrical inclusion under a remote transverse tension.
{\it Int. J. Solids Struct.} {\bf 46} (2009), 1287--1304.

\vspace{-.3em}\bibitem{Miel05ERIS}
A. Mielke.: \newblock   Evolution in rate-independent systems.
\newblock   In: {\it Handbook of Differential Equations,
Evolutionary Equations, 2} (Eds.: Dafermos, C.M., Feireisl,
E.), Elsevier, Amsterdam, 2005, pp.~461--559.

\vspace{-.3em}\bibitem{MieRoudamage} A. Mielke,  T. Roub\'\i\v cek:
\newblock {Rate-independent damage processes
in nonlinear elasticity}. \newblock M$^3\!$AS Math. Models Methods
Appl. Sci. \textbf{16} (2006), 177--209.

\vspace{-.3em}\bibitem{book}
A. Mielke, T. Roub\'\i\v cek: \emph{Rate-Independent Systems,
Theory and Application.} In preparation.

\vspace{-.3em}\bibitem{MiRoTh??DDNE}
A. Mielke, T. Roub\'\i\v cek, M.Thomas:
From damage to delamination in nonlinearly elastic materials at small
strain. In preparation.

\vspace{-.3em}\bibitem{MieThe99MMRI} A. Mielke, F. Theil:  \newblock A mathematical model for rate-independent phase
                  transformations with hysteresis. \newblock
Proceedings of the Workshop on ``Models of Continuum
                  Mechanics in Analysis and Engineering'' (Eds.: Alber, H.-D.  Balean, R.M., and Farwig,
                  R.),
Shaker-Verlag, Aachen, 1999, pp.~117-129.

\vspace{-.3em}\bibitem{MiThLe02VFRI} A. Mielke, F. Theil, V. Levitas:
\newblock A variational formulation of rate--independent
   phase transformations using an extremum principle. \newblock
{\it Arch. Rational Mech. Anal.} \textbf{162} (2002), 137--177.

\vspace{-.3em}\bibitem{MieThe04RIHM}
A. Mielke, F. Theil: \newblock   On rate-independent hysteresis
models.
\newblock  {\it Nonlin. Diff. Equations Appl.} \textbf{11} (2004),
151--189. (Accepted July 2001).

\vspace{-.3em}\bibitem{ThoMie09DNEM}
A. Mielke,  M. Thomas: \newblock  Damage of nonlinearly elastic
materials at small strain: existence and regularity results.
{\it ZAMM} \textbf{90} (2010), 88--112.

\vspace{-.3em}\bibitem{petrov-schatzman2009} A.~Petrov, M. Schatzman:
\newblock Mathematical results on existence for viscoelastodynamic
problems with unilateral constraints. \newblock
  {\it SIAM J. Math. Anal.} \textbf{40} (2009), 1882--1904.

\vspace{-.3em}\bibitem{petrov-schatzman2010} A.~Petrov, M. Schatzman: \newblock
A pseudodifferential linear complementarity problem related to a one
dimensional viscoelastic model with Signorini conditions.
\newblock {\it Arch. Rational Mech. Anal.}, in print (2010).

\vspace{-.3em}\bibitem{Point}
N. Point:
\newblock Unilateral contact with adherence.
\newblock {\em Math. Methods Appl. Sci.} \textbf{10} (1988), 367--381.

\vspace{-.3em}\bibitem{Raous}
M.~Raous, L.~Cang\'emi, M.~Cocu:
\newblock A consistent model coupling adhesion, friction, and unilateral contact.
\newblock {\em Comput. Methods Appl. Mech. Eng.} \textbf{177} (1999), 383--399.

\vspace{-.3em}\bibitem{NPDE_roubicek}
T. Roub\'\i\v cek: {\it Nonlinear Partial Differential Equations
with Applications}. Birkh\"auser, Basel, 2005.

\vspace{-.3em}\bibitem{tr0}
T. Roub\'\i\v cek: Thermo-visco-elasticity at small strains with
$L^1$-data. \emph{Quarterly Appl. Math.} {\bf 67} (2009), 47--71.

\vspace{-.3em}\bibitem{tr2}
T. Roub\'\i\v cek: Rate independent processes in viscous solids at
small strains. {\it Math. Methods Appl. Sci.} {\bf 32} (2009),
825--862.

\vspace{-.3em}\bibitem{tr1}
T. Roub\'\i\v cek: Thermodynamics of rate independent processes in
viscous solids at small strains.  \emph{SIAM J. Math. Anal.}
{\bf 40} (2010), 256-297.

\vspace{-.3em}\bibitem{tr-LS-CZ}
T. Roub\'\i\v cek, L. Scardia,  C. Zanini: Quasistatic
delamination problem.
{\it Cont. Mech. Thermodynam.} {\bf 21} (2009), 223-235.

\vspace{-.3em}\bibitem{SadStu10MTCC}
P.Sadowski, S.Stupkiewicz:
A model of thermal contact conductance at high real contact area fraction.
{\it Wear} \textbf{268} (2010), 77-85.

\vspace{-.3em}\bibitem{simon86}  J.~Simon: \newblock
  Compact sets in the space {$L^p(0,T;B)$}. \newblock
  \emph{Ann. Mat.  Pura Appl.},
  {\bf 146}  (1987),   65--96.

\vspace{-.3em}\bibitem{SoHaSh06AACP}
M. Sofonea, W. Han, M. Shillor: {\it Analysis and Approximation of
Contact Problems with Adhesion or Damage}. Chapman \& Hall/CRC, Boca
Raton, FL, 2006.

\vspace{-.3em}\bibitem{visintin96}
A. Visintin: {\it Models of Phase Transitions}.  Birkh\"auser
Boston, Inc., Boston,  1996.




\end{thebibliography}
\end{document}